\documentclass{amsart}[12pt]
\usepackage[colorlinks=true,          
            linkcolor=blue,
            citecolor=red,
            urlcolor=blue]{hyperref}
\usepackage{amssymb,amsfonts,amsmath,hyperref,mathrsfs,latexsym,stmaryrd,graphicx,manfnt}
\DeclareMathAlphabet{\mathpzc}{OT1}{pzc}{m}{it}

\DeclareMathAlphabet{\mathpzc}{OT1}{pzc}{m}{it}

\newtheorem{defn}{Definition}[section]
 \newtheorem{thm}{Theorem}[section]

 \newtheorem{Prop}{Proposition}[section]
  \theoremstyle{definition}\newtheorem{Ex}{Example}[section]

  \numberwithin{equation}{section}

\usepackage {amsfonts,amssymb}
\usepackage{latexsym}
\usepackage{mathrsfs}
\input xy
\xyoption{all}

\newcommand {\CC}{\mathbb{C}}

\newcommand{\NN}{\mathbb{N}}
\newcommand {\PP}{\mathbb{P}}

\newcommand {\RR}{\mathbb{R}}

\newcommand {\HH}{\mathbb{H}}

\newcommand{\ZZ}{\mathbb{Z}}

\newcommand {\bph}{\boldsymbol{\phi}}

\newcommand {\bdk}{\boldsymbol{k}}
\newcommand {\bdE}{\boldsymbol{E}}
\newcommand {\bdU}{\boldsymbol{U}}
\newcommand {\bdP}{\boldsymbol{P}}

\newcommand {\bk}{{\bf k}}

\newcommand {\bB}{{\bf B}}

\newcommand {\bU}{{\bf U}}

\newcommand{\bun}{\textrm{Bun}}
\newcommand{\rk}{\textrm{rk}}

\newcommand{\cA}{\mathcal{A}}
\newcommand{\cB}{\mathcal{B}}
\newcommand{\cC}{\mathcal{C}}
\newcommand{\cD}{\mathcal{D}}
\newcommand{\cE}{\mathcal{E}}
\newcommand{\cF}{\mathcal{F}}

\newcommand{\cH}{\mathcal{H}}
\newcommand{\cI}{\mathcal{I}}

\newcommand{\cK}{\mathcal{K}}
\newcommand{\cL}{\mathcal{L}}

\newcommand {\cO}{\mathcal{O}}
\newcommand{\cP}{\mathcal{P}}
\newcommand{\cQ}{\mathcal{Q}}
\newcommand{\cR}{\mathcal{R}}
\newcommand{\cS}{\mathcal{S}}
\newcommand{\cT}{\mathcal{T}}
\newcommand{\cU}{\mathcal{U}}

\newcommand{\cX}{\mathcal{X}}

\newcommand{\cih}{\mathpzc{h}}
\newcommand{\cx}{\mathpzc{x}}
\newcommand{\cy}{\mathpzc{y}}

\newcommand{\scB}{\mathscr{B}}
\newcommand{\scC}{\mathscr{C}}

\newcommand{\scE}{\mathscr{E}}
\newcommand{\scF}{\mathscr{F}}

\newcommand{\scP}{\mathscr{P}}

\newcommand{\scV}{\mathscr{V}}

\newcommand{\scZ}{\scZ}

\newcommand{\fa}{\mathfrak{a}}
\newcommand{\fb}{\mathfrak{b}}

\newcommand{\fg}{\mathfrak{g}}

\newcommand{\fl}{\mathfrak{l}}

\newcommand{\fo}{\mathfrak{o}}

\newcommand{\fs}{\mathfrak{s}}

\newcommand{\ft}{\mathfrak{t}}
\newcommand{\slt}{\mathfrak{sl}(2,\CC)}

\newcommand{\fv}{\mathfrak{v}}

\newcommand{\fz}{\mathfrak{z}}

\newcommand{\fD}{\mathfrak{D}}

\newcommand{\sym}{\textrm{Sym}}

\newcommand {\dbar}{\overline{\partial}}

\newcommand{\mhom}{\textrm{Hom}}
\newcommand {\mend}{\textrm{End}}
\newcommand {\misom}{\textrm{Isom}}
\newcommand {\maut}{\textrm{Aut}}
\newcommand{\pr}{\textrm{pr}}

\newcommand {\send}{\underline{End} }

\newcommand {\ad}{\textrm{ad} }

\newcommand{\img}{\textrm{Im }}
\newcommand{\spec}{\textrm{Spec }}

\newcommand{\gspec}{\underline{\textrm{Spec }}}
\newcommand {\cok}{\textrm{coker}}
\newcommand{\tot}{\textrm{tot }}

\newcommand{\ctimes}{\otimes_\CC}

\newcommand{\pic}{\textrm{Pic}}
\newcommand{\tr}{\textrm{tr }}

\newcommand  {\eps}{\varepsilon}

\newcommand {\fii}{\varphi}

\newcommand{\Higgs}{{\bf Higgs}}
\newcommand{\Bun}{{\bf Bun}}

\newcommand{\Prym}{{\bf Prym}}

\newcommand{\rts}{{\sf root}}
\newcommand{\wts}{{\sf weight}}
\newcommand{\crts}{{\sf coroot}}
\newcommand{\cwts}{{\sf coweight}}
\newcommand{\chr}{{\sf char}}
\newcommand{\cchr}{{\sf cochar}}

\newcommand{\hookr}{\hookrightarrow}

\newcommand{\les}[9]{
\xymatrix{
 0 \ar[r] & {#1} \ar[r]  &  {#2} \ar[r]  &  {#3}
\ar@{->}`r/10pt[d] `[l] `^dl[dlll]  `^dr/10pt[dll]    [dll] \\
 &  {#4} \ar[r] & {#5} \ar[r] & {#6}
\ar@{->}`r/10pt[d] `[l] `^dl[dlll]  `^dr/10pt[dll]    [dll] \\
 & {#7} \ar[r]  & {#8} \ar[r] & {#9}
\ar@{->}`r/10pt[d] `[l] `^dl[dlll]  `^dr/10pt[dll]    [dll] \\
 & 0 \ar[r] & \cdots & }
}

\newcommand{\lestwo}[9]{
\xymatrix{     
 0 \ar[r] & {#1} \ar[r]  &  {#2} \ar[r]  &  {#3} 
\ar@{->}`r/10pt[d] `[l] `^dl[dlll]  `^dr/10pt[dll]    [dll] \\
 &  {#4} \ar[r] & {#5} \ar[r] & {#6} 
\ar@{->}`r/10pt[d] `[l] `^dl[dlll]  `^dr/10pt[dll]    [dll] \\
 & {#7} \ar[r]  & {#8} \ar[r] & {#9} }
}

\newcommand{\lesthree}[5]{
\xymatrix{     
 0 \ar[r] & {#1} \ar[r]  &  {#2} \ar[r]  &  {#3} 
\ar@{->}`r/10pt[d] `[l] `^dl[dlll]  `^dr/10pt[dll]    [dll] \\
 &  {#4} \ar[r] & {#5} & }
}

\newcommand{\lesfour}[8]{
\xymatrix{     
 0 \ar[r] & {#1} \ar[r]  &  {#2} \ar[r]  &  {#3} 
\ar@{->}`r/10pt[d] `[l] `^dl[dlll]  `^dr/10pt[dll]    [dll] \\
 &  {#4} \ar[r]^-{#8} & {#5} \ar[r] & {#6} 
\ar@{->}`r/10pt[d] `[l] `^dl[dlll]  `^dr/10pt[dll]    [dll] \\
 & {#7} \ar[r]  & \cdots  &  }
}

\include{biblio}
 \date{\today} 
\title{ Lectures on Higgs Moduli and Abelianisation }
\author{Peter Dalakov}
\address{Institute of Mathematics and Informatics, Bulgarian Academy of Sciences, acad.G.Bontchev str., bl.8, 1113 Sofia Bulgaria}

\begin{document}

\maketitle
\setcounter{tocdepth}{1}
\tableofcontents

\eject
	\section{Introduction}
	      \subsection{}
Almost thirty years ago N.Hitchin introduced, in his seminal works  \cite{hitchin_sd}, \cite{hitchin_sb} 
the notion of a \emph{Higgs field}  -- a
 ``twisted  endomorphism'' of a holomorphic vector bundle $E$ over a compact Riemann surface $X$.
This   ``twisted endomorphism'' is, more precisely,  a sheaf homomorphism
$\theta: E\to E\otimes \Omega^1_X$, i.e., a global section $\theta\in H^0(X,\send E\otimes \Omega^1_X)$. Hitchin's motivation was rooted 
in gauge theory. Its methods of  mathematical exploration of Yang--Mills equations on Riemann surfaces, 
self-duality equations on $\RR^4$ and  magnetic monopoles and instantons on $\RR^3$ have lead to outstanding progress in
these areas 
during the late 1970-ies and early 1980-ies. 
Hitchin considered self-dual connections on (trivial)  $SU(2)$- and $SO(3)$-bundles on  $\RR^4$, and imposed the condition that these
connections be  invariant under translations along $\RR^2\subset \RR^4$. The dimensionally reduced self-duality equations 
turned out to be conformally invariant and hence could be studied on an arbitrary Riemann surface. These  \emph{self-duality equations on Riemann surface},
 now called \emph{Hitchin's equations}, are non-linear elliptic PDE,  involving the data of a Higgs field and
 a  unitary connection on a vector bundle.

Hitchin constructed analytically a moduli space of solutions of the self-duality equations and endowed it with the structure of
an algebraic variety. He showed that the connected component of the moduli space, corresponding to topologically trivial bundles
  is in fact hyper-K\"ahler (orbifold).
Its twistor family (over $S^2$) contains  two non-isomorphic complex structures, both of which admit modular interpretation. One 
of them corresponds to the coarse moduli space of (semi-)stable Higgs bundles of fixed rank and degree zero, while the other is a coarse moduli space
for (semi-) simple local systems of fixed rank.

The hyper-K\"ahler geometry behind the self-duality equations is instrumental in the non-abelian Hodge theory
developed by C.Simpson (\cite{hbls},\cite{simpson_ICM}). This  rich and beautiful theory is completely outside
of the scope of the present lectures, but we direct the reader to the surveys \cite{simpson_hodge},  \cite{alberto_rayan}
and \cite{boalch_naht}, and to the article \cite{KW} for some recent applications to physics.

One of the key insights in Hitchin's seminal works   was the idea
to replace the  data of a Higgs field $\theta$ and
a \emph{vector bundle} $E$
    by the  data
of a \emph{line bundle} $L$ over 
a ramified (``spectral'') cover $S\to X$.
  For a reasonable theory  
the canonical (line) bundle $K_X=\Omega_X^1$
 must  have global sections, i.e., the genus $g= \dim H^0(X,K_X)$ must be at least one.  To obtain a coarse moduli space
of such pairs  Hitchin introduced a notion of stability (extending the notion of stability for vector bundles) and
it turned out that 
 stable pairs  exist only when $g\geq 2$.

The spectral curve $S$ encodes the spectrum of $\theta$, while  the line bundle $L$  encodes the eigenspaces of $\theta$.
The origins of this idea can be traced to an earlier work of Hitchin's on the construction of monopoles in $\RR^3$ (\cite[\S3]{hitchin_monopoles}).
In fact, various kinds of spectral curves over rational and elliptic bases have been extensively used by the integrable systems 
and mathematical physics communities in the context of Lax pairs, solitons and monopoles -- see \cite{adlervm_cis},
\cite{hitchin_monopoles}, \cite{griffiths_flows},
 \cite{adams_harnad_previato_isospectral_1}, \cite{adams_harnad_hurtubise_isospectral_2}, \cite{hurtubise_kk} and the references therein.
One should also recall here 
Atiyah's abelianisation program for the geometric quantisation of the WZW model
(\cite{atiyah_knots}, \cite{hitchin_flat_GQ}, \cite{witten_GQCS}).
The idea to replace non-abelian theta-functions  by  abelian ones 
is still  an active and exciting area of research in mathematics  and physics.

      The most general  version of this  ``spectral correspondence''  relating Higgs  and spectral data
is described in \cite{don-gaits}. It  generalises  many 
works in algebraic geometry,  most notably   \cite{hitchin_sb},
 \cite{faltings},  \cite{donagi_spectral_covers} and \cite{scognamillo_elem}.
The introduction of ``meromorphic'' Higgs fields allows one to consider spectral data on  curves  of low genus, thus connecting
with the earlier work on algebraic integrability.
It seems that the merger of the two flows of research finally occurred after the appearance of
  \cite{markman_thesis}, \cite{bottacin}, \cite{hurtubise_markman_rk2} and \cite{markman_sw}.

In these lectures we begin by a simple motivational setup in which we consider spectral covers arising from
families of matrices  and proceed with considering basic properties of Higgs bundles,
spectral and cameral curves. We then continue with a potpourri of specialised topics: principal $\slt$-subalgebras, the
Kostant and Hitchin sections, special K\"ahler geometry, the Donagi--Markman cubic and the $G_2$ Hitchin system.

Most of these results are classical, except for some of the topics in Sections \ref{flow}, \ref{dm_cubic} and \ref{G2}, where we also include some
recent material and work in progress.

We work exclusively with Higgs and principal bundles whose structure groups are \emph{complex} (simple or reductive) Lie groups.
Hitchin (\cite[\S 10]{hitchin_sd}, \cite{hitchin_teich}) initiated  the study of Higgs bundles and their moduli for \emph{real forms} of
complex reductive Lie groups. This rich and still developing area is yet another sin of omission that we admit to have made. 
We direct the reader to the theses \cite{ana_thesis}, \cite{schaposnik_thesis}, the survey \cite{schaposnik_singapore}, the preprint
 \cite{ana_cameral_su} and the references therein for
 a discussion of spectral and cameral data in such a context.

	      \subsection{Notation and conventions}
Except for the occurence of the group $GL_n$ in Sections \ref{toy_version} and \ref{abel_lin}, we are going to use
principal bundles with structure group $G$ which is assumed to be a simple complex (affine) Lie group.
We are going to denote by $T$ and $B$ Cartan and Borel subgroups of $G$, and by lowercase fraktur letters ($\ft$, $\fb$ and $\fg$)
the corresponding Lie algebras. The root system will be denoted by $\cR$ and the Weyl group by $W$.

We denote by $\ad$ (respectively, $\textrm{Ad}$) the adjoint representation of $\fg$ (respectively, $G$)
in $\mend(\fg)$, respectively, $GL(\fg)$. For $x\in \fg$ we use $\ad x$ and $\ad_x$ interchangeably and denote
by $\zeta_x$ or $\zeta(x)$ the centraliser $\ker\ad x$.
The adjoint group of $G$ is denoted by
$G^{ad}=G/Z(G)\subset GL(\fg)$. 
The adjoint bundle of a principal bundle $P$ will be denoted by $\ad P$.
We deemphasise the difference between a (holomorphic or algebraic) vector bundle and its sheaf of sections.
If we need to underline this difference, we write $\tot$ for the total space
$\tot E:= \gspec\sym^\bullet E^\vee$.

We denote by $\scV$ the vanishing scheme of a function or section of a vector bundle.

       When studying spectral covers, one can work either with $\CC$-schemes (locally of finite type) or with complex-analytic spaces. 
However, when dealing with principal bundles
 it should be kept in mind  that algebraic $G$-bundles need not be Zariski locally trivial.
We are going to work predominantly with holomorphic principal bundles, usually over compact Riemann surfaces.
For a line bundle $L=\cO(D)$ we denote by $|L|$ the corresponding linear system.
 
  The Hitchin base will be denoted by $\cB$ (or $\cB_\fg$ and $\cB_G$ if the structure group must be specified).
The locus of non-singular cameral curves with simple ramification will be denoted by $\scB$.

	      \subsection{Acknowledgements}
This  work is dedicated to Ugo Bruzzo on the occasion of his sixtieth birthday. I would like to thank him for
being a great collaborator and a wonderful friend. The last section of this work is based on my talk at
the conference \emph{Interactions in geometry and physics} held in Guaruj\'a, Brazil in honour of Ugo's birthday.
I would like to thank the organising committee:  F.Sala, M.Jardim, A.Henni, V.Lanza, F.Perroni, J.Scalise and P.Tortella
for the impeccable organisation and endless tolerance.

The first part of these notes is an extended version of a mini-course on abelianisation and
Higgs bundles that I gave at L'Universit\'e d'Angers, France in October 2014. I am grateful to   V.Roubtsov for the invitation
to give the lectures and would like to thank him,
 M.Cafasso, G.Powell, I.Reider and LAREMA for the hospitality. I would also like to thank Volodya for encouraging
me to write up these notes and for his nearly endless patience with regard to this matter.

	\section{Spectral covers from families of matrices}\label{toy_version}
The  spectral construction is rooted in an elementary observation from linear algebra:
a regular semi-simple endomorphism of $\CC^n$ is uniquely determined by the collection of its eigenvalues and
their respective eigenspaces. Consider then $\mend (\CC^n)$ and its  open subset $\mend^{r,ss}(\CC^n)$ consisting of 
regular semi-simple endomorphisms. 
The totality of their spectral data can be described by  a line bundle over an  $n$-fold \'{e}tale cover of $\mend^{r,ss}(\CC^n)$,
thus leading to a toy version of Hitchin's construction. 
We begin by  spelling out the details of this story, which  
 I have learned  from T.Pantev. 

		  \subsection{Regular and Semisimple Endomorphisms}

Let $\fg$ be a reductive
Lie algebra over $\CC$. Recall that an element $x\in \fg$ is  \emph{regular} if its centraliser
$\fz_x=\ker\ad x\subset \fg$ has the lowest possible dimension, i.e., if
$\dim \ker\ad x= \rk \fg$.  In these notes
we shall be dealing exclusively 
 with two cases: $\fg = \fg\fl_n(\CC)$
and $\fg$ a simple complex Lie algebra. 
In the former case, the rank is $n$ and an element is regular precisely when it has a single Jordan block per eigenvalue.
 In the latter case the rank is the dimension of a Cartan subalgebra.

To  begin with the general linear case  we fix a
 complex vector space $V$, $\dim V =n>0$, and 
show that   $\mend^{r,ss}(V)\subset \mend (V)$ is 
a  $PGL(V)$-invariant Zariski open set,  in fact an affine variety.

Indeed, under the adjoint representation 
\[
 \textrm{ad}: \mend (V)\longrightarrow \mend (\mend (V)),\ \phi \longmapsto [\phi,\ ]
\]
an element $\phi$ with eigenvalues $\{\lambda_i\}$ is mapped to 
 $\ad\phi$, whose  eigenvalues are $(\lambda_i-\lambda_j)$.
Among these the  eigenvalue zero appears at least $n$ times and the $n(n-1)$ differences $\lambda_i-\lambda_j$, $i\neq j$,
if nonzero, come in pairs with opposite signs.
The characteristic polynomial of $\ad\phi$ is then
\[
 \det (\lambda \boldsymbol{1} - \ad\phi)=   \sum_{k=n}^{n^2} \lambda^k P_k(\phi) = \lambda^n
\left\{ P_n(\phi)+ \lambda P_{n+1}(\phi) + \ldots \right\}
\]
for some $PGL(V)$-invariant polynomials  $P_k$, $\deg P_k = n^2-k$.
The first of these is the \emph{discriminant}:
$\fD=P_n\in \sym^{n(n-1)}(\mend (V)^\vee)$,
  \[
\fD(\phi)= P_n(\phi)= \pm \prod_{i\neq j} (\lambda_i-\lambda_j)^2.
\] 
We denote by $m_{\ad \phi}$ the multiplicity of $\lambda$ in $\det (\lambda \boldsymbol{1} - \ad\phi)$.
The discriminant of $\phi$ vanishes precisely when some of the differences $\lambda_i-\lambda_j$ vanishes
(for some $i\neq j$), i.e., when $m_{\ad \phi}$ is greater than its minimal value, $n$. 
One sees that, up to a (non-zero) constant multiple,
$\fD$ coincides with the discriminant $\cD$ of the characteristic polynomial $\det (\lambda\boldsymbol{1}- \phi)$ of $\phi$.
The divisor
\[
 \scV(\fD)= \left\{\phi\left| m_{\ad\phi}> n \right.\right\}\subset \mend V
\]
 turns out to be precisely 
 the  complement of the regularly semisimple locus.
      \begin{Prop}
Let  $V$ be a complex vector space of dimension  $n$ and discriminant $\fD\in \sym^{n(n-1)}(\mend (V)^\vee)$. Then
\[
 \scV(\fD)  = \mend (V) \backslash \mend^{r,ss} (V)\subset \mend (V).
\]
      \end{Prop}

\emph{Proof:}
The complement of $\mend^{r,ss}V$ in $\mend V$ is the union of the sets of non-regular and of non-semisimple
endomorphisms.
Since $n\leq \dim\ker\ad \phi\leq m_{\ad \phi}$, $\phi$-nonregular implies $\phi\in \scV(\fD)$.
We show next that any non-semisimple element is contained in $\scV(\fD)$. 
Indeed,   
 $\phi\in \mend(V)$ is semisimple or nilpotent precisely when $\ad\phi$ is so. 
 Consider the Jordan decomposition $\ad\phi = (\ad\phi)^{ss}+ (\ad \phi)^{nlp}$. 
 If  $\phi\notin \mend^{ss} V$ then   $(\ad\phi)^{nlp}\neq 0$ and so $(\ad\phi)^{ss}\in \mend (\mend V)$ is non-regular, leading to
 $m_{(\ad \phi)^{ss}}>n$. Since the characteristic polynomials of $\ad\phi$ and $(\ad\phi)^{ss}$ coincide, 
we get $m_{\ad\phi}>n$, i.e., $\phi\in\scV(\fD)$. This shows that $ \mend (V) \backslash \mend^{r,ss} (V)\subset\scV(\fD)$.

For the opposite inclusion observe that $\fD(\phi)=0$ if and only if $m_{\ad\phi}>n$. The latter inequality,  together with  
$m_{\ad\phi}\geq \dim\ker\ad\phi\geq n$ implies that either $\ad \phi$ (and hence $\phi$) is non-semisimple, or
$\dim\ker\ad\phi>n$, i.e., $\phi$ is non-regular.

\qed

More generally, given a  simple complex Lie algebra $\fg$ and an element $x\in\fg$ we can consider the characteristic polynomial
\[
\det\left(\lambda\boldsymbol{1}- \ad x\right) = \sum_{k=\rk \fg}^{\dim \fg} \lambda^k P_k(x)
\]
of  $\ad x$.   The discriminant $\fD= P_{\rk \fg}\in \sym^{\dim\fg-\rk \fg}(\fg^\vee)$ is a $G$-invariant polynomial,  and in fact
$\fD= \prod \alpha$,  the product over all roots. Essentially the same argument as above shows that 
$\fg^{r,ss}=\fg\backslash \scV(\fD)$.
We shall return to this situation in Section \ref{cameral}.

      \begin{Ex}
Let $\fg=\slt$ and $x\in \fg$. Then $  \det\left(\lambda\boldsymbol{1}- \ad x\right) = \lambda^3 + \lambda 4\det x$ 
and $\det \left(\lambda\boldsymbol{1}-  x\right)= \lambda^2- \det x$.
      \end{Ex}

	    \subsection{The spectral cover}
An endomorphism $\phi\in\mend^{r,ss}(V)$
is determined by  its spectrum and the collection of  eigenspaces (eigenlines), i.e., by: 
    \begin{enumerate}
     \item $\{\lambda_1,\ldots, \lambda_n\}$, $\lambda_j\in\CC$, distinct
     \item $\{ L_1,  \ldots,  L_n\}$,  $L_i\in \PP(V)$ distinct, $V= L_1\oplus \ldots\oplus L_n$.
    \end{enumerate}
A crucial  point -- subsumed  in the notation -- is that we have  fixed a bijection  between the sets $\{\lambda_i\}_i$ and $\{L_i\}_i$,
so that $L_i$ is the $\lambda_i$-eigenspace of $\phi$. 

We now let these data move in a family: suppose $S$ is a complex algebraic (or analytic) variety 
and $\Phi: S\to \mend V$ a morphism.
 These determine a collection of endomorphisms 
$\{\phi_s=\Phi(s)\}_{s\in S}$
which  can be assembled   in  a ramified $n$-to-one cover  $p:\widetilde{S}\to S$, cut out by the characteristic
polynomial(s):
  \[
   \xymatrix{\widetilde{S}=\{\left. (s,\lambda)\right| \det (\lambda\boldsymbol{1}-\phi_s)=0\}\ar[drr]_-{p}\ar@{^{(}->}[rr]<-0.5ex> &  & S\times \CC\ar[d]_-{\textrm{pr}_1}\\
									      &         &S
}.
  \]
The fibre  $p^{-1}(s)\subset \widetilde{S}$  is the spectrum of  $\phi_s$ and thereby $\widetilde{S}$
 is called \emph{the spectral cover of $S$, corresponding to the family $\Phi$}.
There are different ways of thinking about $\widetilde{S}$,  each of them carrying the
germs of possible generalisations. We discuss these below.
		  \subsection{Global Spectra}
We have defined $\widetilde{S}\subset S\times \CC$ as a hypersurface, cut out by a specific equation. 
Now $S\times\CC$ is the total space
$\tot\cO_S=\gspec\sym^\bullet\cO_S$ and 
$\textrm{pr}_{S\ast}\cO_{S\times\CC} =\sym^\bullet\cO_S\simeq \CC[\lambda]\ctimes\cO_S$.
Then $\Phi\in \mend\left(V\ctimes\cO_S\right)$
and $\widetilde{S}$ 
is   the closed subscheme (complex subspace)
\[
\widetilde{S}= \scV\left(\det\left(\lambda\boldsymbol{1} - \textrm{pr}_S^\ast\Phi\right)\right) \subset \tot\cO_S.
\]
The spectral cover  is also a   ``global spec'' of a sheaf of
$\cO_S$-algebras.
Indeed, let $\cI$ be the
ideal sheaf, generated by the image of the $\cO_S$-module homomorphism
\[
 \left(a_n(\Phi),\ldots, a_1(\Phi),1,0,0\ldots\right)^t: \cO_S\longrightarrow \sym^\bullet\cO_S.
\]
Then there is an $\cO_S$-module isomorphism
$\sym^\bullet\cO_S/\cI\simeq_{\cO_S} \cO_S^{\oplus n}$ and
\[
\widetilde{S}=\gspec\sym^\bullet\cO_S/\cI.
\]

Equivalently, the characteristic polynomial of $\Phi$ can be thought of as a morphism 
(of $S$-varieties)
$\tot\cO_S\to\tot\cO_S$, and $\widetilde{S}$ is the preimage of the $0$-section  
$S\subset \tot\cO_S$.

		  \subsection{The Adjoint Quotient as a proto-Hitchin map}

We start by describing
the cover $\widetilde{S}\to S$
as a pullback (via  $\Phi$) of a certain ``universal'' cover of $\mend (V)$.

Consider first the standard
 $n$-to-one cover $\pi:\widetilde{\CC^n}\to \CC^n$, defined by
\[
\widetilde{\CC^n}=\left\{ (a_1,\ldots, a_n; \lambda)\left| \lambda^n+a_1\lambda^{n-1}+\ldots +a_n=0 \right. \right\}\subset \CC^n\times \CC,
\]
with  $\pi=\left.\textrm{pr}_1\right|_{\widetilde{\CC^n}}$. The fibre  $\pi^{-1}(\underline{a})$ consists of the roots of
the unique monic, degree-$n$ polynomial, having  coefficients $\underline{a}=(a_1,\ldots, a_n)$.
 The total space $\widetilde{\CC^n}\subset \CC^n\times\CC$ is a (smooth) hypersurface and the fibres of $\pi$ are complete intersections, 
so $\pi$  is flat  by the relations criterion for flatness.

Define a  morphism
 $h:\mend (V)\to \CC^n$ by  $h(\phi)=(a_1(\phi),\ldots a_n(\phi))$, where 
 $a_i(\phi) = (-1)^i\textrm{tr}\Lambda^i\phi$.
By the very definition of the spectral cover we have an isomorphism $\widetilde{S}\simeq (h\circ \Phi)^\ast \widetilde{\CC^n}$
and a commutative diagram
\[
 \xymatrix{
	  \widetilde{S}\ar[d]_-{p}\ar[r]&h^\ast\widetilde{\CC^n}\ar[r]\ar[d] &\widetilde{\CC^n}\ar[d]_-{\pi}\\
          S\ar[r]^-{\Phi}          &\mend V\ar[r]^-{h}& \CC^n\\
}.
\]
Since $\pi$ is flat,  $p$ must be  flat as well.

We note in passing that there is another, $n!$-to-one, cover of $\CC^n$:  the quotient morphism for  the standard action 
of the symmetric group $S_n$ on $\CC^n$. Pulling that cover via $h\circ\Phi$ would lead us to the notion of a \emph{cameral cover}, 
to be discussed in section \ref{cameral}.

The fibre of $h$ over $\underline{a}$ consists
of all endomorphisms  whose characteristic polynomial has coefficients $\underline{a}=(a_1,\ldots, a_n)$, i.e., 
 $h^{-1}(\underline{a})\subset \mend (V)$
is the
closure  of an adjoint $GL(V)$-orbit. 
We can decompose any   $\phi \in\mend V$ into a semi-simple and nilpotent part, $\phi= \phi^{ss}+\phi^{nlp}$, and for 
any $t\in\CC^\times$, $\phi^{ss}+ t\phi^{nlp}$ is conjugate to $\phi$, so  $\phi^{ss}\in \overline{GL(V)\cdot \phi}$.
 In general, the closure
   $\overline{GL(V)\cdot \phi}$ contains several orbits,  among which there is a unique closed orbit
(that of  $\phi^{ss}$), and a unique open orbit (that of a regular element with the same spectrum as $\phi$).

The space of  orbit closures  is the GIT quotient
\[
  \mend V\sslash GL(V)\simeq \mend^{ss}V\slash GL(V)\simeq \mend^{reg}V/GL(V) \simeq \CC^n.
\]
 The map $h$ is the adjoint quotient map and a  precursor of the Hitchin map.
 Again, we shall return to this discussion and its analogues for
arbitrary simple groups $G$ in Section \ref{cameral}.

It is easy to give equations for the 
 cover $h^\ast\widetilde{\CC^n}$: it is  the hypersurface 
\[
 \mend(V)\times\CC \supset h^\ast\widetilde{\CC^n}=\left\{ (\phi,\lambda)\left| \lambda^n +\sum_{i=1}^n a_i(\phi)\lambda^{n-i}=0 \right. \right\}
\]
and its ramification  locus is
\[
\mend(V)\times\CC \supset
 \left|
	\begin{array}{l}
	 \lambda^n +\sum_{i=1}^n a_i(\phi)\lambda^{n-i}=0\\
	 n\lambda^{n-1}+\sum_{i=1}^n (n-i)a_i(\phi)\lambda^{n-i-1}=0\\
	\end{array} \right. .
\]
The branch locus, i.e., the discriminant locus of the morphism $h^\ast\widetilde{\CC^n}\to \mend(V)$,
is cut out by
 the discriminant of the characteristic polynomial, i.e., by the  singular  hypersurface  
\[
\scV(\fD)= \scV(\cD) = \left\{\phi\left| \fD(\phi)=0 \right.\right\} \subset \mend V,
\]  
away from which we have 
an \'etale $n$-to-one cover $\left. h^\ast\widetilde{\CC^n}\right|_{\mend^{r,ss}V}\to \mend^{r,ss}V$.
The hypersurface   $h^\ast\widetilde{\CC^n}\subset\mend V\times\CC$ is singular. 
Its  singular locus is contained in its  ramification locus   and
is cut by the additional equation $\sum_{i=1}^n \lambda^{n-i}da_i=0$.

      \begin{Ex}
If $V=\CC^2$, the double cover $h^\ast\widetilde{\CC^2}\to \textrm{Mat}_2(\CC)$ is given by
\[
 \left\{\left. (A,\lambda)\right| \lambda^2 -\lambda \tr A  + \det A =0\right\}\subset \textrm{Mat}_2(\CC)\times\CC.
\]
It is branched  over the vanishing  locus of  $\cD =\tr^2-4\det$, which is nothing but the set of $2\times 2$ matrices with
coinciding eigenvalues. The singular locus  of $h^\ast\widetilde{\CC^2}$ is   the line
\[
 \left\{\left. (z\boldsymbol{1},z)\right| z\in\CC\right\}\subset h^\ast\widetilde{\CC^2}
\]
lying over the line of scalar matrices  in $\textrm{Mat}_2(\CC)$. 
      \end{Ex}
Once again, there is another approach that can be taken here:
 the $GL_n(\CC)$-orbit of any $\phi\in \mend^{r,ss}(\CC^n)$  contains (and is determined by)
unique $S_n$-orbit of diagonal matrices.
 These orbits give rise to an $S_n$-Galois cover of $\mend^{r,ss}(\CC^n)$ and this is how cameral covers (Section \ref{cameral}) come into being. 
	    \subsection{The Spectral Sheaf}

We have  encoded the spectra of  $\Phi\in\mend (V\ctimes\cO_S)$
 in the  cover 
$\widetilde{S}\to S$,  but there is more to an endomorphism than just
its spectrum: one has to keep track of the invariant subspaces (or quotients).
We shall hence  ``decorate'' the spectral cover with a sheaf which encodes these data.

Let  $\textrm{pr}_S=\textrm{pr}_1: \tot\cO_S=S\times\CC\to S$ be the canonical projection.
Recall (\cite{H}, Exc.II.5.17.) that
the functor $\textrm{pr}_{S\ast}$ induces   an equivalence of
categories
\[
\xymatrix@1{\left\{ \textrm{Quasi-coherent } \cO_{S\times\CC}\textrm{-modules}\right\}\ar[r]^-{\textrm{pr}_{S\ast}} & \left\{\textrm{Quasi-coherent } \sym^\bullet\cO_S\textrm{-modules}\right\} }.
\]
It sends a quasi-coherent $\cO_{S\times\CC}$-module $\scF$   to the quasi-coherent $\cO_S$-module
 $\textrm{pr}_{\ast}\scF$, which  is
$\cO_S$-\emph{coherent} precisely when
$\textrm{supp }\scF\hookr \tot\cO_S\to S$ is a finite morphism.

The structure of a  $\sym^\bullet\cO_S$-module on 
an $\cO_S$-(quasi)coherent sheaf $E$ is determined by a section
$\fii\in \mend E$. Such a  pair $(E,\fii)$ could be called a
\emph{$\cO_S$-valued (quasi)coherent Higgs sheaf}, but we delay the definition until the next section.
 The    $\sym^\bullet\cO_S$-action   on   $E=\textrm{pr}_{\ast}\scF$  is 
generated by $\textrm{pr}_\ast(\lambda\cdot\ )$, where
 $\lambda\in \Gamma\left(\tot\cO_S, \textrm{pr}_S^\ast \cO_S\right)$ is the tautological section.
An inverse equivalence (in the coherent case) can be obtained  by sending the pair
$(E,\fii)$ to the  kernel of the restriction of
$  \left(\lambda\boldsymbol{1} - \textrm{pr}_S^\ast\fii\right)$ to ${{\textrm{supp }\textrm{cok} \left(\lambda\boldsymbol{1} - \textrm{pr}_S^\ast\fii\right)}}$.
The composition of the two inverse equivalences is  \emph{isomorphic} to the
 identity functor  and  this isomorphism
is controlled by the ramification data of the cover. 
We shall discuss this in  detail for the   cases of interest
in Section \ref{abel_lin}, Proposition \ref{4trm}, and refer to
\cite{don-gaits} or \cite{don-pan}, Appendix A for more general situations.

For the case at hand,
the  $\cO_S$-module $E=V\ctimes \cO_S$ becomes  a $\sym^\bullet\cO_S$-module 
via the action of $\Phi\in \mend (V\ctimes \cO_S)$.
In this situation we can give another, somewhat analytical, description of the spectral sheaf.
Denote by
  $\cE\in\Gamma(T_{\mend V})$ the Euler vector field on $\mend V$.
Then
 $p_1^\ast\cE -\lambda \boldsymbol{1}$ 
is a section
of 
the trivial bundle $\mend V\ctimes \cO_{\mend V\times\CC}  = p_1^\ast\send \left(\cO_{\mend V}\ctimes V\right)$.
It determines a short exact sequence
\[
 \xymatrix@1{ 0\ar[r]& \cO\ctimes V\ar[r]^-{\ \cE -\lambda \boldsymbol{1}} &\cO\ctimes V\ar[r] &\cQ\ar[r]& 0 }
\]
 on $\mend V\times \CC$, and $\cQ$ is a torsion sheaf, supported on $h^\ast\widetilde{\CC^n}$.
Upon restriction to $h^\ast\widetilde{\CC^n}$   we obtain a 4-term exact sequence
\[
 \xymatrix@1{ 0\ar[r]&\cK\ar[r]& \cO\ctimes V\ar[r] &\cO\ctimes V\ar[r] &\cQ\ar[r]& 0 },
\]
which  can then be  pulled back to $\widetilde{S}\subset S\times\CC$ by $\Phi$.

So far we have not imposed any restrictions on the image of $\Phi$, i.e., on the kind of endomorphisms
$\{\phi_s\}$ which arise in the family.
Suppose now that 
$\img\Phi \subset \mend^{r,ss}(V)$. Then the cover $\widetilde{S}\to S$ is \'etale, 
$\widetilde{S}$ is smooth,
 $\widetilde{L}=\Phi^\ast(\cK)$ is a line bundle and so  is $\Phi^\ast\cQ$.
The family $\{\phi_s\}$ is then determined by the pair $\left(p: \widetilde{S}\to S, \widetilde{L}\right)$ or by
 $\left(p: \widetilde{S}\to S, \Phi^\ast\cQ \right)$.
In particular,
 $p_\ast \widetilde{L}\simeq   \cO_S\ctimes V$.

If  the condition $\img \Phi \subset \mend^{r,ss}(V)$ is not satisfied, various complications ensue.
Indeed, the subsheaf $\widetilde{L}\subset \cO_{\widetilde{S}}\otimes V$
may not be a line bundle anymore:  its rank jumps if the endomorphisms $\phi_s$ have eigenspaces of  dimension greater than one. Moreover, 
$\widetilde{S}$ may be singular (or even reducible or non-reduced) and the cover
$\widetilde{S}\to S$  may be ramified.
Thus our construction is not very flexible so far. For one, since $\mend^{r,ss} V$  is the complement of
a divisor $ \scV(\fD)\subset \mend V$, the 
  eigenvalues  of $\{\phi_s\}$ can coalesce over divisors in $S$. A possible amendment is to 
require that
$\img\Phi\subset \mend^{reg}V$, the complement of a codimension three
subset in $\mend V$.
Regular endomorphisms have
 exactly one eigenline per eigenvalue, so while $\widetilde{S}\to S$ may be ramified,
 $\widetilde{L}$ is still of rank one.

A second deficiency of this simplified setup is that
  $\mend(V)$ is affine and thus the base $S$ cannot be proper. 
Aside from this, one may want   to replace $\mend V$ with another complex Lie algebra.
Thus, regularity aside,  there  are at least three  modifications that
one can attempt in order to obtain richer examples:
      \begin{enumerate}
       \item Replace  $S\times V=\tot V\ctimes\cO_S$ with a (suitable)  vector bundle $E\to S$. Replace 
 $\Phi\in\mend \left(V\ctimes\cO_S\right)$ with a global section $\fii\in H^0(S,\send E)$.
       \item Replace  $\tot\cO_S$ with a (suitable)
line or vector bundle $K$ (``coefficient bundle''). Combined with $(1)$, this means that
one considers sections 
$\fii\in H^0(S,\send E\otimes K)$. 
       \item Replace $GL(V)$ with a complex reductive group $G$. Respectively, replace $ E\to S$ with a principal $G$-bundle $P\to S$
and
replace $\send E$ with $\ad P$.
      \end{enumerate}
In these lectures we shall  (mostly) confine ourselves to  the case when $S$ is   a compact Riemann surface $X$.
The group $G$ will be either $GL_n(\CC)$ or a simple complex Lie group.
The coefficient bundle $K$ will be a line bundle, mostly the canonical line bundle
$K_X$ or its twist $K_X(D)$ by a divisor $D \geq 0$.
But before restricting to these situations let us make some final general remarks.

        Spectral covers with  coefficient bundle $K$ having $\rk K>1$ arise, for example,  in C.Simpson's non-abelian Hodge theory (\cite{simpson_ICM}, \cite{hbls}),
but see also  \cite{don-gaits}. To make sense of these, one needs to impose additional restrictions on $\fii$. Indeed, 
a  trivialisation of $K$ on an open $U\subset S$ identifies  $\fii\in H^0(U,\send E\otimes K)$   with
$\sum_i\fii_i$, where
$\fii_i\in \mend E_U$. Each of the $\fii_i$ gives rise to  a spectral cover  $\widetilde{U}_i\to U$, and we can take
the fibred product of all these covers. To guarantee that the result  is independent of trivialisation and glues to a global object, one 
imposes the integrability condition $\fii\wedge\fii=0$.
In these notes, Higgs bundles with coefficients in a vector bundle will appear only briefly in section \ref{sk}.

The non-abelian Hodge-theoretic viewpoint allows one to use Higgs bundles and spectral data for
 describing (certain) $D$-branes of type $B$ on cotangent bundles to compact K\"ahler manifolds, see 
\cite{pantev_wijnholt}, \cite{kop_dbranes}.

      \section{$K_X$-valued $G$-Higgs bundles on curves}\label{higgs_intro}

Let us  fix now  a compact Riemann surface $X$ of genus $g\geq 2$ and  denote by $K_X=\Omega^1_X$ its canonical bundle.
\begin{defn}
 A \emph{$K_X$-valued $G$-Higgs bundle on $X$} is a pair $(P,\theta)$, consisting of    a holomorphic principal $G$-bundle
$P$ and a section
$\theta\in H^0(X,\ad P\otimes K_X)$. If $V$ is a holomorphic vector bundle on $X$, a \emph{$V$-valued $G$-Higgs bundle on $X$} is
a pair $(P,\theta)$, consisting of a holomorphic principal $G$-bundle $P$ and a section
 $\theta\in H^0(X,\ad P\otimes V)$.
\end{defn}
In these notes we are going to limit ourselves to the case $G=GL_n(\CC)$ or $G$ a simple complex Lie group.
If $G\subset GL_n(\CC)$ is  a classical group, one can also work with the associated (Higgs) vector bundle. This is a  pair
$(E,\fii)$, consisting of a vector bundle $E=P\times_G \CC^n$, possibly equipped with some additional structure
(e.g., a quadratic form) and a twisted endomorphism $\fii$,
preserving this extra structure.
 For example, an    $SL_2(\CC)$ Higgs vector bundle is a
  pair
$(E,\fii)$, where:
    \begin{itemize}
     \item $E$ is a rank-2 vector bundle with $\det E\simeq \cO_X$
     \item $\fii\in H^0(X,\send_0 E\otimes K_X)$.
    \end{itemize}
We use the standard notation
  $\send_0 E\subset\send E$ for the sheaf
of trace-free endomorphisms of $E$.

If $P$ is a (holomorphic) principal $G$-bundle, a choice of Killing form gives an isomorphism $\ad P\simeq \ad P^\vee$
and consequently  $H^1(X, \ad P)^\vee\simeq H^0(X,\ad P\otimes K_X)$ by 
Serre duality.
The infinitesimal deformations of a pair $(P,\theta)$ are controlled, as shown in
 \cite{Biswas-Ramanan},  \cite{hitchin_teich}, \cite{nitsure}, by
$\HH^1(\scC_{(P,\theta)}^\bullet)$, where $\scC_{(P,\theta)}^\bullet$ is the complex
\begin{equation}\label{BR_complex}
\xymatrix@1{\ad P\ar[r]^-{\ad\theta} & \ad P\otimes K_X}.
\end{equation}
Being a shifted cone, $\scC_{(P,\theta)}^\bullet$   fits in the short   exact sequence of complexes
\begin{equation}\label{cone}
 \xymatrix@1{0\ar[r]&\ad P\otimes K_X[-1]\ar[r]&\scC_{(P,\theta)}^\bullet\ar[r]&\ad P\ar[r]&0}.
\end{equation}
Taking Euler characteristics and applying the Hirzebruch--Riemann--Roch theorem gives
\[
 -\chi\left(\scC^\bullet_{(P,\theta)}\right) = \chi(\ad P\otimes K_X) -\chi(\ad P) = \dim G \deg K_X = \dim G \cdot (2g-2),
\]
which is in fact the dimension of the (local) moduli space of $(P,\theta)$ if $G$ is simple.
Furthermore, the   long exact sequence of hypercohomology, corresponding to (\ref{cone}) gives
   \begin{equation}\label{cone_hyper}
  \xymatrix@1{(0)\ar[r] & \cok h^0(\ad\theta )\ar[r]& \HH^1(\scC_{(P,\theta)}^\bullet)\ar[r]^-{\pr}& \ker h^1(\ad\theta)\ar[r]& (0)}    ,
    \end{equation}
where $h^i(\ad\theta): H^i(\ad P)\to H^i(\ad P\otimes K_X)$ are the natural maps induced by $\theta$. Explicitly, 
$\cok h^0(\ad\theta)= \frac{H^0(\ad P\otimes K)}{\ad \theta (H^0(\ad P))}$.

The Grothendieck--Serre duality pairing for $\scC^\bullet_{(P,\theta)}$, i.e.,  Serre duality for hypercohomology,
gives, by the self-duality of $\scC^\bullet_{(P,\theta)}$, a pairing
$\HH^1(\scC^\bullet)^{\otimes 2}\to \CC$. As shown in 
 \cite[Theorem 4.3]{Biswas-Ramanan},  this pairing is
  a symplectic form $\omega_{can}$ on $\HH^1(\scC_{(P,\theta)}^\bullet)$.

In   the  case when $H^0(X, \ad P)=0$ 
we get that $\omega_{can}$ is the bilinear form determined by
      \[
       \xymatrix{
	    0\ar[r]& H^0(\ad P\otimes K_X)\ar[d]_-{-SD^\vee}^-{\simeq}\ar[r]& \HH^1(\scC_{(P,\theta)}^\bullet)\ar@{-->}[d]_-{\omega_{can}}\ar[r]& H^1(\ad P)\ar[r]\ar[d]_-{\simeq}^-{SD}& 0\\
	    0\ar[r]& H^1(\ad P)^\vee\ar[r]& \HH^1(\scC_{(P,\theta)}^\bullet)^\vee\ar[r]& H^0(\ad P\otimes K_X)^\vee\ar[r]& 0\\
	  }.
      \]

This is the case,  for instance, if  $(E,\fii)$ is   an $SL_2(\CC)$ Higgs pair with $\maut E= Z(SL_2(\CC))\simeq \ZZ/2$,
since for such a bundle $\mend_0 (E)=0$.

It is easy to see that a choice of hermitian metric on $\ad P$ gives a splitting of the  extension (\ref{cone}),
 whence $\omega_{can}$ is identified with
the canonical symplectic form on $\HH^1(\scC_{(P,\theta)}^\bullet)\oplus \HH^1(\scC_{(P,\theta)}^\bullet)^\vee$.
  Moreover, such a choice (of Hermitian metric) allows us to
use Hodge theory and work with Dolbeault representatives of hypercohomology. In terms of the decomposition
\[
  A_X^1(\ad P)\simeq A_X^{0}(\ad P\otimes K_X)\oplus A_X^{0,1}(\ad P)
\]
the symplectic form is
\begin{equation}\label{symplectic_dolbeault}
 \omega_{can}((\eta',\eta''),(\xi',\xi''))=\textrm{tr}\int_X \eta'\wedge \xi''-\xi'\wedge\eta'' =\textrm{tr}\int_X(\eta'+\eta'')\wedge (\xi'+\xi'').
\end{equation}

We turn now to discussing some properties of the coarse moduli spaces of Higgs bundles. Such moduli spaces do in fact exist, but for a mildly restricted
class of Higgs bundles -- the (semi)stable ones.
    \begin{defn}
     An $SL_n(\CC)$-Higgs bundle $(E,\fii)$ is \emph{(semi-)stable} if  any proper subbundle $F\subset E$, which is $\fii$-invariant in the sense that 
$\fii(F)\subset F\otimes K_X$, satisfies $\deg F<0$ (resp. $\deg F\leq 0$).
    \end{defn}
More generally,  a $GL_n(\CC)$-Higgs bundle $(E,\fii)$  is  semi-stable if  for any proper $\fii$-invariant subbundle $F\subset E$, the inequality
$\mu(F) \leq \mu(E)$ holds. Here $\mu=\deg/\rk$ denotes  the \emph{slope}  of a vector bundle (\cite{SSS}).
Notice that if a subbundle $F\subset E$ is  $\fii$-invariant,  then $\fii$ induces a Higgs field $\overline{\fii}$
on $E/F$.

Any ($GL_n(\CC)$) Higgs pair $(E,\fii)$ has  unique maximal semi-stable Higgs subbundle and hence
 admits a canonical increasing (Harder--Narasimhan) filtration $(E_\bullet,\fii_\bullet)$. It has the (defining) property  that
$(\textrm{gr}_k E_\bullet,\textrm{gr}_k \fii_\bullet)$ is the maximal semi-stable Higgs subbundle of $(E/E_{k-1},\overline{\fii})$.
Next,  every semi-stable Higgs
pair $(E,\fii)$ admits an increasing (Jordan--H\"older) filtration $(E_\bullet,\fii_\bullet)$,
   whose associated graded pieces are stable Higgs pairs with $\mu(\textrm{gr}_k E_\bullet)=\mu(E)=\mu(E_k)$ for all $k$.
Any two such filtrations $(E_\bullet,\fii_\bullet)$ and $(E'_\bullet,\fii'_\bullet)$ have the same length and
consequently the associated graded Higgs pairs $\left(\oplus_k \textrm{gr}_k E_\bullet ,\oplus_k\textrm{gr}_k \fii_\bullet\right)$ and
 $\left(\oplus_k \textrm{gr}_k E'_\bullet ,\oplus_k\textrm{gr}_k \fii'_\bullet\right)$
are isomorphic. If the pair $(E,\fii)$ is stable to begin with, its Jordan--H\"older filtration is trivial and the associated
graded object is isomorphic to  the original pair. 
Two semi-stable ($GL_n(\CC)$ or $SL_n(\CC))$ Higgs pairs $(E,\fii)$ and $(F,\psi)$ are said to be \emph{S-equivalent}, 
written $(E,\fii)\sim_S (F,\psi)$, if their associated graded pairs for the Jordan--H\"older filtration
are isomorphic. This gives rise to an equivalence relation which on stable Higgs pairs coincides with isomorphism.
The S-equivalence class of a semi-stable pair  contains unique polystable pair, i.e., a Higgs bundle, which is a direct
sum of stable Higgs bundles of the same slope.

The first results on moduli spaces of (rank two) Higgs bundles are due to Hitchin. In particular, he
proved the following theorem.

    \begin{thm}[\cite{hitchin_sd}] 
There exists a connected, quasi-projective  coarse moduli space
\[
 \Higgs_{SL_2(\CC),X}=\left\{(E,\fii) \textrm{ semi-stable } SL_2(\CC)\textrm{-Higgs bundle }\right\}/\sim_S
\]
of dimension $6g-6= \dim SL_2(\CC)\deg K_X$. It contains as a Zariski-dense open subvariety $T^\vee_{\textrm{Bun}^{reg}_{SL_2}}$, the
cotangent bundle to the smooth locus of the coarse moduli space of $SL_2$-bundles.
    \end{thm}

An analogous result holds for $SL_n(\CC)$-Higgs bundles (\cite{hitchin_sb}, \cite{moduli2}),
 the dimension of the moduli space now being $2(n^2-1)(g-1)$.

To define  semi-stability for a  $G$-Higgs bundle $(P,\theta)$ one can consider the associated Higgs vector bundle
$(\ad P,\ad\theta)$ and require the latter to be semi-stable. However, in this transition from $G$ to $G^{ad}=G/Z(G)\subset GL(\fg)$ some subtle information
is lost.
 To have  an adequate and intrinsic notion
of stability (not just semi-stability)  one needs to extend Ramanathan's approach (\cite{ramanathan}) to the case of Higgs bundles.
Notice that if a vector bundle $E$ admits a subbundle $F$, then its structure group is reduced from $GL_n$ to a maximal parabolic subgroup.
More generally, given a $G$-bundle $\pi: P\to X$, Ramanathan considers reductions $\sigma:X\to P/H$ of $P$ to a maximal parabolic subgroup $H\subset G$.
Let $\pi_H: P/H\to X$ be the reduced bundle and 
$T_{\pi_H}=\ker d\pi_{H}$  its relative tangent bundle.
 Definition 1.1, \emph{ibid.}  states  that $P$ is (semi-)stable, if for any such $\sigma$ the inequality
$\deg \sigma^\ast T_{\pi_H}>0$, respectively $\deg \sigma^\ast T_{\pi_H}\geq 0$, holds. 
Now $P\to P/H$ is a principal $H$-bundle, $P_H$, and any reduction $\sigma$ gives rise to a projection
$\ad P\to \ad P\left/ \sigma^\ast P_H\right.$. Given a 
Higgs bundle $(P,\theta)$ we say that a reduction $\sigma$ is $\theta$-invariant (or \emph{a Higgs reduction}), if $\theta$ is in
the kernel of this projection. Finally, we say that $(P,\theta)$
 is (semi-)stable if $\deg \sigma^\ast T_{\pi_H}>0$ (resp.  $\deg \sigma^\ast T_{\pi_H}\geq 0$) holds
for any Higgs reduction $\sigma$.

It turns out that $(P,\theta)$ is semi-stable if and only if $(\ad P,\ad\theta)$ is semi-stable, but a stable Higgs bundle may
have a strictly semi-stable adjoint Higgs bundle.
For comparisons between the different notions of (semi, poly) stability and discussion of $S$-equivalence, Harder--Narasimhan  and Jordan--H\"older filtrations
of $G$-Higgs bundles on curves and on higher-dimensional varieties,
we direct the reader to \cite{biswas_anchouche_eh}, \cite{arijit_parthasarathi_hn},
\cite{bea_jh} and \cite{ugo_bea_ss}. 

Coarse moduli spaces $\Higgs_{G,X}$ of semi-stable $G$-Higgs bundles exist, more generally,  for (affine) reductive groups $G$, as shown in
\cite{hitchin_sd}, \cite{hitchin_sb}, \cite{hbls} and \cite{moduli2}.
For simple $G$ the  spaces $\Higgs_{G,X}$  are singular normal quasi-projective varieties (\cite{hbls}). They have 
(\'{e}tale or analytic) local models $\HH^1(\scC_{(P,\theta)}^\bullet)/\maut(P,\theta)$, which are
orbifold singularities whenever $\maut(P,\theta)\varsupsetneq Z(G)$. The connected components of the moduli space
are labelled by $\pi_1(G)$, i.e., by the topological type of the $G$-bundle, underlying the Higgs pair,
see (\ref{conn_comp}).

The spaces $\Higgs_{G,X}$ carry rich geometry, as was observed  by Hitchin in \cite{hitchin_sd} for 
  $G=SL_2(\CC)$ and $G=PGL_2(\CC)$.

	\begin{thm}[\cite{hitchin_sd}]
      The space $\Higgs_{SL_2,X}$ carries a holomorphic symplectic structure $\omega_{can}$
which extends the canonical
(Liouville) symplectic structure on $T^\vee_{\textrm{Bun}^{reg}_{SL_2}}$. The map
\[
 h: \Higgs_{SL_2,X}\longrightarrow \cB_{SL_2}=H^0(X,K_X^2), \quad h\left([E,\fii]\right)=\det\fii
\]
is a proper, surjective morphism with Lagrangian fibres. 
Every  $b\in\cB$
determines   a natural double cover $\widetilde{X}_b\to X$. It is  non-singular precisely when
$b\in\cB$ has simple zeros. In that case, 
  the fibre 
$h^{-1}(b)$ is non-singular and
is a translate of the Prym variety $\Prym_{\widetilde{X}_b/X}\subset \pic \widetilde{X}_b$ of $\widetilde{X}_b\to X$.
	\end{thm}

The Hitchin map $h$ is not only proper -- it is in fact projective, and is the affinisation morphism of $\Higgs_{SL_2,X}$.
Non-singular Hitchin fibres are torsors over abelian varieties, which are (compactified) Prym varieties of (integral) double covers of $X$.
We shall return to this discussion, but in a more general setup, in Section \ref{abel_lin}.

For arbitrary (affine) reductive groups $G$
the moduli spaces $\Higgs_{G,X}$ are  holomorphic symplectic as well, see \cite{hitchin_sb}, \cite{Biswas-Ramanan}, \cite{moduli2}. 
Their symplectic structure can be expressed either via the duality  pairing for $\HH^1(\scC^\bullet_{(P,\theta)})$,
 or in Dolbeault terms,
as in (\ref{symplectic_dolbeault}). 
There is a natural linear form $\lambda_{(P,\theta)}\in \HH^1(\scC_{(P,\theta)}^\bullet)^\vee$, namely $\lambda(v)=\pr (v)(\theta)$,
see (\ref{cone_hyper}), which gives rise to a 1-form on $T^\vee\bun_G^{reg}$.
 In
\cite[Theorem 4.3]{Biswas-Ramanan} it is shown that
$\omega_{can}=-d\lambda$.  Thus the  open embedding
$T^\vee_{\textrm{Bun}^{reg}_{G}}\subset \Higgs^{reg}_{G,X}$ is a symplectomorphism.

As noted at the end of Section \ref{toy_version}, one can consider Higgs bundles with different ``coefficients''.
In some of the later sections we shall deal with $K_X(D)$-valued Higgs bundles, also known as \emph{meromorphic Higgs bundles}.
These are pairs $(P,\theta)$, consisting of a principal $G$-bundle $P$ and a Higgs field $\theta\in H^0(X,\ad P\otimes K_X(D))$,
where $D$ is a sufficiently positive effective divisor on $X$.
For such pairs  an appropriate deformation complex can be written, and the duality pairing gives rise to a holomorphic Poisson structure, 
see \cite{markman_thesis}, \cite{bottacin} or  the surveys \cite{donagi_markman}, \cite{markman_sw}, \cite{survey_mero_higgs}.
We shall denote the corresponding coarse moduli space by $\Higgs_{G,D,X}$ or $\Higgs_{G,D}$.
As discovered by Markman (\cite{markman_thesis}) and Bottacin (\cite{bottacin}), $\Higgs_{G,D}$  carries the structure of a Poisson 
completely integrable system, sometimes called \emph{the generalised Hitchin system}. 
We shall return to it in Sections \ref{cameral} and \ref{dm_cubic}.

The inclusion $T^\vee_{\textrm{Bun}^{sm}_{G}}\subset \Higgs_{G,X}$ is strict, and this will
be crucial  in Section \ref{flow}. To see this in the case $G=SL_2(\CC)$ we can follow \cite[\S\S 1, 10]{hitchin_sd}
and  fix a theta-characteristic (spin structure) $K_X^{1/2}$.
The bundle $E=K_X^{1/2}\oplus K_X^{-1/2}$ can be endowed with a nilpotent Higgs field $\fii = \begin{pmatrix}
                                                             0&0\\
							     1&0\\ 
                                                            \end{pmatrix}$.
The pair $(E,\fii)$ is a stable one,  
 while the bundle $E$ is unstable, being destabilised by $K_X^{1/2}\subset E$. 
An analogous example can be constructed for other simple groups by taking a principal
homomorphism $SL_2(\CC)\to G$,
 see \cite{hitchin_teich} and Section \ref{flow}.

      \section{Abelianisation for $GL_n$ and $SL_n$}\label{abel_lin}
We review in this section some basic properties of spectral curves and  abelianisation for the special and general linear group.
For further properties and insights we direct the reader to \cite{hitchin_sb}, \cite{bnr}, \cite{donagi_spectral_covers},
 \cite{kouv-pan} and \cite{hausel_pauly}.

	  \subsection{The General Linear Group}
We start by introducing the \emph{$GL_n$-Hitchin base}. For us this will be 
 the vector space of global sections of the rank $n$ vector bundle $\bU=\bigoplus_{k=1}^n K_X^n$, i.e.,
\[
 \cB_{GL_n}=H^0(X,K_X)\oplus\ldots \oplus H^0(X,K_X^n).
\]
Whenever there is no chance of confusion we  write  $\cB$ for $\cB_{GL_n}$.
By the Riemann--Roch theorem the dimension of the base is
\[
 1 + \sum_{k=1}^n\left(k \deg K_X + (1-g)\right)= 1 + \frac{n(n+1)}{2}\deg K_X - n(g-1)=n^2(g-1)+ 1,
\]
or, more intrinsically,
\[
  \dim \cB_{GL_n}= \frac{1}{2}\dim GL_n(\CC)\deg K_X + \dim Z(GL_n).
\]
  Notice the ``miraculous  
 coincidence of dimensions'' (\cite{hitchin_sb})  $\dim\cB_{GL_n}=\dim \textrm{Bun}_{GL_n}$.

Let us denote by $Y$ the non-compact surface $\tot K_X$ and by $\pi$ the
bundle projection $Y=\tot K_X\to X$.
 Any
  $b=(b_1,\ldots, b_n)\in \cB_{GL_n}$ gives rise to a  section 
\begin{equation}\label{eqn_spectral}
 \lambda^{\otimes n}+ \pi^\ast b_1 \lambda^{\otimes (n-1)}+\ldots + \pi^\ast b_n
\end{equation}
of $\pi^\ast K_X^n$,
whose vanishing locus is, by definition, 																																																																																
\emph{the spectral curve $\widetilde{\cC}_b\subset Y$ associated to $b\in\cB_{GL_n}$}.
If we  let $b\in\cB_{GL_n}$ vary, then (\ref{eqn_spectral}) becomes a section of $p_Y^\ast \pi^\ast K_X^n$
on $\cB\times Y$, whose vanishing locus is the 
\emph{universal spectral curve}
$\widetilde{\cC}\to \cB_{GL_n}\times X$:
\[
 \xymatrix{\widetilde{\cC}_b\ar@{^{(}->}@<-0.5ex>[r]\ar[d]_-{\pi_b}&\widetilde{\cC}\ar[r]\ar[d]_-{\pi}& Y= \tot K_X\ar[d]\\
	  \{b\}\times X\ar@{^{(}->}@<-0.5ex>[r]& \cB_{GL_n}\times X\ar[r]^-{\textrm{ev}}&\tot K_X^n\\
 }.
\]

In this way we have ``stacked together'' the  different spectral curves:
$\widetilde{\cC}=\bigcup_{b\in\cB}\widetilde{\cC}_b$. The individual spectral curves may be
singular, or reducible, or non-reduced, but the total space of the family $\widetilde{\cC}$
is in fact smooth, see \cite[Corollary 1.2]{kouv-pan}.

We proceed  now with the analogues of the various constructions from Section \ref{toy_version}.
First, since $Y=  \gspec (\sym^\bullet K_X^{-1})$, we have an isomorphism of quasi-coherent $\cO_X$-modules
\[
 \pi_\ast\cO_Y =  \sym^\bullet K_X^{-1}= \cO_X\oplus K_X^{-1}\oplus K_X^{-2}\oplus\ldots .
\]
Denoting by $\cI_b$ the ideal sheaf
\[
 \cI_b= \textrm{Im}\left(b_n,\ldots, b_1, 1, 0,0,\ldots\right)^t: K_X^{-n}\longrightarrow \sym^\bullet K_X^{-1},
\]
we can describe the spectral curve as the global spec
\[
 \widetilde{\cC}_b= \gspec \left(\sym^\bullet K_X^{-1}/\cI_b \right)
\]
and correspondingly
\begin{equation}\label{spectral_push}
 \pi_{b \ast} \cO_{\widetilde{\cC}_b}= \sym^\bullet K_X^{-1}/\cI_b \simeq_{\cO_X} \cO_X\oplus K_X^{-1}\oplus \ldots \oplus K_X^{-n+1}.
\end{equation}

To get an explicit local description of $\widetilde{\cC}_b$ over an affine open $U\subset X$ we choose a nowhere vanishing section (generator)
$u\in K_X^{-1}(U)$, so that $\left. \tot K_X\right|_U =\spec \cO_X(U)[u]$.
Let
$\overline{b_k}=b_k(u^k)=\langle b_k, u^k \rangle\in \cO_X(U)$ denote the outcome of the  evaluation pairing between $b_k$ and $u^k$.
 Then  the cover $\left. \widetilde{\cC}_b\right|_U\to U$ is determined by the ideal
 $\cI_b(U)\subset \cO_X(U)[u]$, generated by
 \[
u^n + \overline{b_1}u^{n-1}+\ldots + \overline{b_n}=\langle b_n+\ldots +b_1+1, u^n\rangle\in\cO_X(U)[u].
\]
Taking determinants in (\ref{spectral_push}) we obtain
 \[
\det \pi_{b \ast} \cO_{\widetilde{\cC}_b}= K_X^{-\frac{n(n-1)}{2}},\ \deg \pi_{b \ast} \cO_{\widetilde{\cC}_b}=- n(n-1)(g-1).
 \]
Next, since  $K_X^n\otimes \pi_\ast\cO_Y = \bigoplus_{k\geq 0} K_X^{n-k}\simeq \bU\oplus \pi_\ast\cO_Y$, the projection formula gives rise to an isomorphism
\begin{equation}\label{sections_Y}
 H^0(Y,\pi^\ast K_X^n)= H^0(X,\pi_\ast\pi^\ast K_X^n)\simeq H^0(K_X^n\otimes \pi_\ast\cO_Y) = H^0(\cO_X\oplus \bU) = \CC\oplus \cB,
\end{equation}
\[
 b_0 \lambda^{\otimes n}+ \pi^\ast b_1 \lambda^{\otimes (n-1)}+\ldots + \pi^\ast b_n \longmapsto (b_0,b_1,\ldots , b_n).
\]
Let us  compactify $\pi:Y\to X$ to the ruled surface $\pi: \overline{Y}=\PP(\cO\oplus K_X)=\underline{\textrm{Proj}}\sym^\bullet(K_X^{-1}\oplus\cO_X)\to X$, 
and denote by 
 $\cO_{\overline{Y}}(1)$  the relative hyperplane bundle, satisfying $\pi_\ast \cO_{\overline{Y}}(1)\simeq K_X^{-1}\oplus \cO_X$.
Let us also denote by
$\mu\in H^0(\overline{Y},\cO_{\overline{Y}}(1))$ and $\lambda\in H^0(\overline{Y},\pi^\ast K\otimes\cO_{\overline{Y}}(1))$
  the infinity and zero sections.
Then the  isomorphism (\ref{sections_Y}) extends to
\[
 H^0(\overline{Y}, \pi^\ast K_X^n\otimes \cO_{\overline{Y}}(n))\simeq H^0(X,\cO_X) \oplus \cB_{GL_n},
\]
\[
  b_0\lambda^n+ \pi^\ast b_1\lambda^{n-1}\mu+\ldots + \pi^\ast b_n \mu^n \longmapsto (b_0,\ldots, b_n).
\]

 One defines an \emph{$n$-sheeted spectral curve} as an element of
the linear system $\left|\pi^\ast K_X^n\otimes \cO_{\overline{Y}}(n)\right|$ which is contained in $Y\subset\overline{Y}$.
We see that this is an affine open which can be identified with  $\cB_{GL_n}$, and so any spectral curve is the
spectral curve, associated to some $b\in\cB_{GL_n}$.
The  linear system  $\PP\left(\CC\oplus\cB_{GL_n}\right)$ is base-point free, as shown in \cite[\S5]{hitchin_sb}.
 Indeed, any base point must occur along the zero section $X\subset \tot K_X$ (since the linear system contains $\lambda^n$) and thus has to be
 a base point of $K_X^n$. But  $\left| K_X^n\right|$ is base-point free for $n\geq 2$. In the exceptional case $n=1$ (which we usually ignore), 
$\cB_{\CC^\times}=H^0(X,K_X)$ and $|K_X|$
  is  base-point free if and only if $X$ is not hyperelliptic. By Bertini's theorem, the generic spectral curve
$\widetilde{\cC}_b$ is smooth.

The genus $\widetilde{g}$ of a spectral cover $\widetilde{\cC}_b$ can be computed by the
Grothendieck--Riemann--Roch theorem: the equality
\[
 1-\widetilde{g}= \chi(\cO_{\widetilde{\cC}_b})=\chi(\pi_{b \ast}\cO_{\widetilde{\cC}_b})= - n(n-1)(g-1) - n(g-1)
\]
implies
$\widetilde{g}= n^2(g-1)+1 = \dim \cB_{GL_n}$. 
Suppose  that $\widetilde{\cC}_b$ is smooth, and choose a line bundle $[L]\in \pic^d \widetilde{\cC}_b$.
 Then $E=\pi_\ast L$ is a rank $n$ vector bundle on $X$ and $\chi(L)=\chi(\pi_{b\ast}L)$
implies
$\deg E= d - (n^2-n)(g-1)$. The same is true if $\widetilde{\cC}_b$ is integral and $L$ is rank one, torsion-free sheaf.

Thus the pair $(\widetilde{\cC}_b, L)$ determines a $GL_n(\CC)$ Higgs (vector) bundle $(E,\fii)= (\pi_{b\ast L},\pi_{b\ast}(\lambda\otimes) )$ on $X$.
Conversely, to any  pair $(E,\fii)$ we  assign a point $b=(b_1,\ldots, b_n)\in\cB_{GL_n}$, to be denoted 
$h(E,\fii)$, by setting  $b_i=(-1)^i\textrm{tr}(\Lambda^i\fii)$.
 If the corresponding cover $\widetilde{\cC}_b$ is smooth,
we have a $\fii$-eigenline bundle over $\widetilde{\cC}_b\backslash \{\pi^{-1}(B)\}$,
   $B:=\textrm{Bra}(\pi_b)$.
This line bundle extends to a line bundle $L$ over  all of $\widetilde{\cC}_b$, but with a twist along the ramification $\cO(R)=\pi_b^\ast K_X^{-n}$, as shown in
  (\cite{hitchin_sb} and \cite{bnr}).
Here is the precise statement, in a mildly generalised setup.

\begin{Prop}[\cite{bnr}, \cite{hitchin_sb},\cite{hurtubise_is}]\label{4trm}
 Let $K$ be a line bundle on $X$  and $(E,\fii)$ a rank $n$,  $K$-valued Higgs bundle on $X$.
Suppose that the spectral cover $\pi_b:\widetilde{\cC}_b\to X$ corresponding to $b=h(E,\fii)$ is non-singular.
Then there exists a line bundle $L$ on $\widetilde{\cC}_b$ which  fits in the   exact sequence
\begin{equation}
  \xymatrix@1{0\ar[r]& L(-R)\ar[r]&\pi_b^\ast E\ar[r]^-{\ \pi_b^\ast\fii-\lambda \boldsymbol{1}}&\pi_b^\ast E\otimes \pi_b^\ast K\ar[r]& \pi_b^\ast K\otimes L\ar[r]& 0}
\end{equation}
and satisfies $\pi_{b\ast }L\simeq E$.
\end{Prop}
\emph{Proof:}
We outline the proof of this very useful statement, essentially following \cite[\S 4.3]{hurtubise_is}.
See also \cite[Remark 3.7]{bnr}  and \cite[\S 5]{hitchin_sb}.

 We set, slightly abusively, $Y=\tot K$,
even though $K$ need not be the canonical bundle.
There is a
short exact sequence
\begin{equation}\label{3trm_1}
\xymatrix@1{0\ar[r] & \pi^\ast(E\otimes K^{-1}) \ar[r]^-{\ \pi^\ast\fii-\lambda \boldsymbol{1}} &\pi^\ast E\ar[r]& \cQ\ar[r] & 0}, 
\end{equation}
where $\cQ$ is a torsion sheaf, supported on $\widetilde{\cC}_b\subset Y$. We compactify  $Y$ to
the ruled surface $\pi:\overline{Y}=\PP(K\oplus\cO_X)\to X$ and write $\lambda$ and $\mu$ for the zero and
infinity sections, as before. 
Then $\mu\otimes \pi^\ast\fii-\lambda\otimes \boldsymbol{1}$ is a global section  of
$\pi^\ast\left( \send E\otimes K\right)\otimes \cO_{\overline{Y}}(1)$, and, since $\cO_{\overline{Y}}(1)$ is trivial on $Y\subset \overline{Y}$,
the sequence (\ref{3trm_1}) is the restriction to $Y$ of the exact sequence
      \begin{equation}\label{3trm_2}
  \xymatrix@1{0\ar[r] & \pi^\ast(E\otimes K^\vee)\otimes \cO_{\overline{Y}}(-1) \ar[r]^-{\ \mu \fii-\lambda \boldsymbol{1}}  &\pi^\ast E\ar[r]& \cQ\ar[r] & 0}
      \end{equation}
on $\overline{Y}$.
We can determine $\pi_\ast\cQ$ by pushing (\ref{3trm_2}) down  to $X$.  Indeed, the
long exact sequence of $R^\bullet\pi_\ast$, combined with the projection formula gives
\[
 \xymatrix@1{ E\otimes K^\vee\otimes\pi_\ast \cO_{\overline{Y}}(-1) \ar[r]& E\otimes\pi_\ast\cO_{\overline{Y}}\ar[r]& \pi_\ast\cQ\ar[r]& 
E\otimes K^\vee\otimes R^1\pi_\ast \cO_{\overline{Y}}(-1) }
\]
But   $\pi_\ast\cO_{\overline{Y}}\simeq \cO_X$  since $\pi$ has connected  fibres,
  and $R^i \pi_\ast  \cO_{\overline{Y}}(-1) =0$,
since $H^i(\PP^1,\cO_{\PP^1}(-1))=0$, so
$E\simeq \pi_\ast\cQ$.
Restricting the sequence (\ref{3trm_2}) to $\widetilde{\cC}_b\subset Y\subset \overline{Y}$ we obtain a 4-term exact sequence
      \begin{equation}\label{4trm1}
  \xymatrix@1{0\ar[r]& \cK\ar[r] & \pi_b^\ast(E\otimes K^\vee)  \ar[r]  &\pi_b^\ast E\ar[r]& L\ar[r] & 0},
      \end{equation}
and $L=\left. \cQ\right|_{\widetilde{\cC}_b}$ satisfies $\pi_{b\ast }L=\pi_\ast\cQ\simeq E$.
Finally,  splitting (\ref{4trm1}) into two short sequences and taking determinants, we obtain $\cK=L\otimes \pi_b^\ast K^{-n}$.
\qed

We formulate now the
 ``spectral correspondence'' for $GL_n(\CC)$.

      \begin{Prop}[\cite{bnr},Proposition 3.6]\label{bnr_correspondence}
Let   $b\in \cB_{GL_n}$ and  $\widetilde{\cC}_b\subset \tot K$
  an integral spectral curve. There is a natural bijection between
the sets
	\begin{itemize}
	 \item Isomorphism classes of torsion-free, rank one  $\cO_{\widetilde{\cC}_b}$-modules
	 \item Isomorphism classes of $GL_n(\CC)$ Higgs bundles $(E,\fii)$ on $X$,  satisfying $h(E,\fii)=b$.
	\end{itemize}
The bijection sends the class of a sheaf $L$ on $\widetilde{\cC}_b$ to the class of the Higgs bundle
$(E,\phi)= \left(\pi_{b\ast} L, \pi_{b\ast}(\lambda\otimes)\right)$ on $X$.
      \end{Prop}

\emph{Proof:}
Since the morphism $\pi_b: \widetilde{\cC}_b\to X$ is  affine, by \cite[Exc.II.5.17]{H}  the functor $\pi_{b\ast}$
induces an equivalence between the category of quasi-coherent $\cO_{\widetilde{\cC}_b}$-modules
and the category of quasi-coherent $\pi_{b\ast}\cO_{\widetilde{\cC}_b}$ modules, i.e.,
 quasi-coherent $\cO_X$-modules, having the structure of a $\pi_{b\ast}\cO_{\widetilde{\cC}_b}$-module.

We argue now that $\pi_{b\ast}$  preserves the subcategories of torsion-free sheaves.
Let $\cR_X$ and $\cR_{\widetilde{\cC}_b}$ be sheaves of rational functions on $X$ and $\widetilde{\cC}_b$, respectively.
By integrality, both of these  sheaves  are constant  in the Zariski topology and moreover 
$\cR_{\widetilde{\cC}_b}\simeq \pi_b^\ast\cR_{X}$.
A sheaf $\cF$ is torsion-free if and only if $\textrm{tor}\cF = \ker\left(\cF\to \cF\otimes_\cO\cR\right)$
vanishes. So pushing forward  an injection $L\hookr L\otimes_{\cO_{\widetilde{\cC}_b}}\cR_{\widetilde{\cC}_b}$ we
obtain an injection $\pi_{b\ast}L\hookr \pi_{b\ast }L\otimes \cR_X$ by the projection formula. Notice here
that while $\cR_X$ is not locally free, it is quasi-coherent (\cite[Ex.5.2.5]{H}) and the projection formula holds for affine morphisms
and pairs of quasi-coherent sheaves. Conversely, given an injection $E\hookr E\otimes \cR_X$ and an isomorphism
$E\simeq \pi_{b\ast }L$, we obtain that $\pi_{b\ast}\textrm{tor}L=0$ and hence $\textrm{tor}L=0$.

The rank of a sheaf is determined as the rank at the generic point (whose local ring is the field of rational functions).
Since $\deg\pi_b=n$, $\pi_{b\ast}$ identifies rank one torsion-free sheaves on $\widetilde{\cC}_b$ with rank $n$ torsion-free sheaves 
on $X$. But $X$  is non-singular, so these are locally free, i.e.,  rank $n$ vector bundles, having  a
 $\pi_{b\ast}\cO_{\widetilde{\cC}_b}$-module structure.

We now argue that the 
 structure of a   $\pi_{b\ast}\cO_{\widetilde{\cC}_b}\simeq\sym^\bullet K^{-1}/\cI_b$-module
on a rank $n$ vector bundle $E$ is equivalent to the data of a Higgs field $\fii$ with characteristic polynomial $b$.
Indeed, such a module structure on $E$
is determined by an algebra homomorphism 
$\sym^\bullet K^{-1}/\cI_b\to \send E$. By (\ref{spectral_push}) this is equivalent to the data of  an $\cO_X$-module homomorphism
$\fii: K^{-1}\to \send E$, i.e., a Higgs field $\fii\in H^0(X, \send E\otimes K)$ which  satisfies a degree-$n$ polynomial equation, 
i.e., sends (\ref{eqn_spectral}) to zero. Since the rank of $E$ is $n$ and $\widetilde{\cC}_b$ is reduced and irreducible,
  the characteristic polynomial of $\fii$ coincides with its minimal polynomial and with
 $b\in\cB_{GL_n}$. Conversely, given a pair $(E,\fii)$ with characteristic polynomial $b=h(E,\fii)$, we have
that $\fii$ satisfies its own characteristic equation by the Cayley--Hamilton theorem, and hence determines an algebra homomorphism
$\sym^\bullet K^{-1}/\cI_b\to \send E$. Clearly two pairs $(E,\fii)$ and $(E',\fii')$ with $b=h(E,\fii)=h(E',\fii')$
are isomorphic
precisely when $E$ and $E'$ are isomorphic as $\sym^\bullet K^{-1}/\cI_b$-modules. 

Passing to isomorphism classes we obtain the required bijection.

\qed

We make now some further comments and remarks concerning the last Proposition.

 To relate explicitly the module structure on
the spectral sheaf with the data of a Higgs bundle, notice that
if $L$ is a sheaf of abelian groups on $\widetilde{\cC}_b$ which admits an $\cO_{\widetilde{\cC}_b}$-module structure,
then such a structure is determined by an algebra homomorphism $\cO_{\widetilde{\cC}_b}\to \send L$. Pushing forward this homomorphism
we obtain the algebra homomorphism $\pi_{b\ast}\cO_{\widetilde{\cC}_b}\to \send (\pi_{b\ast} L)$, which in turn is determined by a
Higgs field $\fii$ on $E=\pi_{b\ast}L$.

If $K=K_X$, the locus of integral spectral curves is a non-empty Zariski open of  codimension at least $g-1$ in $\cB_{GL_n}$, see \cite[\S 1]{kouv-pan}.
For arbitrary $K$ the open of integral curves is nonempty if there exists
  $b_n\in H^0(X, K^n)$  which is not of the form $b_n=g^m$ for $g\in H^0(X,K^{n/m})$,
$n/m\in \NN$,
see \cite[Remark 3.1]{bnr}.

Let $\cB^o\subset \cB$ be the locus of integral spectral curves, and $\widetilde{\cC}^o\to\cB^o$ the restriction of the
universal spectral cover to that locus. By Grothendieck's theorem (\cite{Grothendieck_descente_6}, \cite[\S 9.4]{FGA_explained}) there exists a
relative Picard scheme $\pic_{\widetilde{\cC}^o/\cB^o}\to \cB^o$. Its fibres away from the locus $\scB\subset \cB^o\subset\cB$ of non-singular spectral curves
have non-proper connected components. In particular (\cite[\S 1.7]{kouv-pan}), for a generic point $b$ of the discriminant divisor $\cB\backslash\scB$
the spectral curve $\widetilde{\cC}_b$ has unique ordinary double point as a singularity and the components of $\pic_{\widetilde{\cC}_b}$ are 
$\CC^\times$-bundles over abelian varieties.
By integrality of fibres, the results from \cite{altman_kleiman_cj} and \cite{dsouza_cj} imply the existence of compactified
Picard schemes $\overline{\pic}_{\widetilde{\cC}^o/\cB^o}\to \cB^o$, whose fibres are fine moduli spaces for rank-one torsion-free sheaves on the respective
spectral curves.
In fact, since all spectral curves $\widetilde{\cC}_b\subset \tot K$ are planar (and hence their singularities have embedding dimension no bigger than two), 
the boundary locus in $\overline{\pic}_{\widetilde{\cC}_b}$ is ``not too big'', i.e., $\pic_{\widetilde{\cC}_b}\subset \overline{\pic}_{\widetilde{\cC}_b}$, $b\in\cB^o$ is a dense open, 
see \cite{altman_kleiman_cj} and \cite{dsouza_cj}.

The Beauville--Narasimhan--Ramanan Proposition, together with the properness of the Hitchin map $h^0:\Higgs_X^o\to \cB^o$ (\cite{hitchin_sd}) implies that
we can identify (non-canonically) the Hitchin fibre $h^{-1}(b)$ and $\overline{\pic}_{\widetilde{\cC}_b}$ for $b\in \cB^o$. Globally
$\Higgs^o$ is a $\overline{\pic}$-torsor, with multiplication $\overline{\pic}_{\widetilde{\cC}^o/\cB^o}\times_{\cB^o}\Higgs_X^o\to \Higgs_X^o$
induced by tensoring the spectral sheaf by a rank-one torsion-free sheaf. As checked previously, this action preserves the connected components, 
i.e. $\overline{\pic}^{d}_{\widetilde{\cC}^o/\cB^o}$ acts simply-transitively on $\Higgs^{o, d-(n^2-n)(g-1)}_{GL_n,X}$.

We note that the properness of the Hitchin map $h:\Higgs_X\to\cB$ is proved for the moduli space of \emph{semi-stable} Higgs bundles, but in fact the Higgs bundles which 
have \emph{integral} spectral curves are \emph{stable}. Indeed, by \cite[ Proposition 1.1]{kouv-pan} an irreducible component of a spectral curve is again a spectral curve
(corresponding to a Higgs bundle of smaller rank). Hence if $(E,\fii)$ has a proper $\fii$-invariant subbundle $V$, the spectral curve of
$(V,\left. \fii\right|_V)$ will be an irreducible component of the spectral curve of $(E,\fii)$, thus violating the integrality of the latter.
Thus unstable and strictly semi-stable Higgs bundles are to be found in Hitchin fibres $h^{-1}(b)$ for $b\in \cB\backslash\cB^o$.

Finally, as shown in \cite{schaub_spectral}, the conclusion of the Proposition remains true even if $\widetilde{\cC}_b$ is not assumed to be integral,
but is an arbitrary spectral curve.

%
%
%
%
%
%
	\subsection{The Special Linear Group}
The Beauville-- Narasimhan-- Ramanan result (Proposition \ref{bnr_correspondence}) can easily be upgraded to a statement about
$SL_n(\CC)$-Higgs bundles.
For that we consider  pairs
$(E,\fii)$, satisfying $\det E=\cO_X$, $\textrm{tr}\fii=0$ and having fixed characteristic polynomial $h(E,\fii)=b$
with an integral spectral curve $\widetilde{\cC}_b$. 
Here
$b$ must be a point in the
  $SL_n(\CC)$-Hitchin base, which is defined as the codimension-$g$ subspace
\[
 \cB_{SL_n}\simeq H^0(K_X^2)\oplus\ldots \oplus H^0(K_X^n)\subset \cB_{GL_n}.
\]
To preserve the correspondence, upon passing from $\cB_{GL_n}$ to $\cB_{SL_n}$ we must accordingly restrict the
class of $\cO$-modules $L$ on the spectral curve: we require these  to satisfy the condition $\det\pi_{b\ast} L\simeq \cO_X$.

The condition $\det\pi_{b\ast} L\simeq \cO_X$ can be understood in terms of the norm homomorphism associated  to the
covering $\pi_b: \widetilde{\cC}_b\to X$. We recall the definition here and refer to \cite[6.5.5]{EGAII}  and \cite[21.5]{EGAIV} 
for the full details.  Identify  $\pi_{b\ast} \cO_{\widetilde{\cC}_b}$ with an $\cO_X$-submodule of 
$\send(\pi_{b\ast}\cO_{\widetilde{\cC}_b})$
by mapping a local section to the endomorphism given by multiplication with that section. 
We obtain then
a homomorphism (of multiplicative monoids)
$\det: \pi_{b\ast}\cO_{\widetilde{\cC}_b}\subset \send(\pi_{b\ast}\cO_{\widetilde{\cC}_b})\to \cO_X$.
Next, to an invertible $\pi_{b\ast}\cO_{\widetilde{\cC}_b}$-module $\cL$ we can associate an invertible $\cO_X$-module $N_{\widetilde{\cC}_b/X}\cL$
by mapping transition functions (with respect to a cover) to their determinants. Finally, we define a group homomorphism
\[
 Nm_{\widetilde{\cC}_b/X}: \pic \widetilde{\cC}_b\to\pic X
\]
via  $[L]\mapsto [N(\pi_{b\ast}L)]$. The fact that this is a group homomorphism follows from the properties of determinants.
If the spectral cover is non-singular, $Nm$ can be identified with the usual pushforward of divisors,
given by
$\cO_{\widetilde{\cC}_b}(\sum n_i p_i)\mapsto \cO_X(\sum n_i \pi_b(p_i))$, see \cite[21.5]{EGAIV} .

We refer for more details to EGA and to \cite{hausel_pauly}, where applications to non-reduced spectral curves are discussed.
If $E$ is a torsion-free, rank-$r$ $\cO_{\widetilde{\cC}_b}$-module and $L$ an invertible $\cO_{\widetilde{\cC}_b}$-module,
by Proposition 3.8 of \cite{hausel_pauly} one has $\det \pi_\ast(E\otimes L) =\det \pi_\ast E\otimes Nm(L)^{\otimes r}$,
which implies the well-known formula $\det\pi_\ast L= Nm(L)\otimes \det(\pi_{b\ast}\cO_{\widetilde{\cC}_b})$.
This formula  can also be taken as a definition of $Nm$.

Having computed $\det (\pi_{b\ast}\cO_{\widetilde{\cC}_b})$ by  Riemann--Roch  we obtain
\[
 \det\pi_{b \ast} L=  Nm_{\widetilde{\cC}_b/X}(L)\otimes K_X^{-\frac{n(n-1)}{2}}.
\]
Hence the  fibres of $h:\Higgs_{SL_n}\to \cB_{SL_n}$ over
$b\in \scB_{SL_n}\subset \cB_{SL_n}$ (i.e.,  corresponding to non-singular spectral curves)
are identified with the solutions of the inhomogeneous linear equation 
$Nm(L)=K_X^{\frac{n(n-1)}{2}}$
in $\pic \widetilde{\cC}_b $:
\[
 \xymatrix{Nm^{-1}\left(K_X^{\frac{n(n-1)}{2}}\right)\ar@<-0.5ex>@{^{(}->}[r]\ar[dr] &\pic \widetilde{\cC}_b \ar[d]^-{Nm}\\
						     &        \pic X                            \\  	
  }.
\]
This  implies that if $h^{-1}(b)\neq \varnothing$, it
is a torsor for $Nm_{\widetilde{\cC}_b/X}^{-1}(\cO_X)=: \textrm{Prym}(\widetilde{\cC}_b/X)$.
The non-emptiness of the fibre follows from the surjectivity of $Nm$, which can be seen as follows. If $M$ is a line bundle on $X$ and
$\pi_b^\ast M\simeq \cO_{\widetilde{\cC}_b}$, we must have
 $M\otimes \pi_{b\ast}\cO_{\widetilde{\cC}_b}\simeq \pi_{b\ast}\cO_{\widetilde{\cC}_b}$, and the latter bundle has a 
non-zero section, being isomorphic to $\oplus_{i=0}^{n-1}  K_X^{-i}$.
But then $\deg M=0$ and it  admits a section (since $\deg K_X>0$), so $M\simeq \cO_X$. In other words, $\pi_b^\ast: \pic X\to \pic \widetilde{\cC}_b$
is injective, and hence $Nm$, being its transpose, is surjective. See also \cite[3.10]{bnr}.

If the spectral curve $\widetilde{\cC}_b$ is integral, 
a similar  description of $h^{-1}(b)$ can be given in terms of compactified Jacobians, see  \cite[ \S 1.7.]{kouv-pan}.
      \subsection{Other groups}\label{other_grps}
A natural question  arising at this point is the question of identifying some kind of ``spectral data''
that can be used to describe $G$-Higgs bundles for an arbitrary (e.g., reductive) group $G$.
Hitchin in  \cite{hitchin_sb}, \cite{hitchin_G2} treated  the case of classical groups 
by reducing the study of principal $G$-bundles to the study of holomorphic vector bundles with extra structure (bilinear form)
and related the Hitchin fibre to a Prym of Jacobian variety.
Such an approach can also be adapted to work for the exceptional group $G_2$, see \cite{hitchin_G2}.
For a complete treatment, then, one has to study the spectral correspondence for pairs consisting of a group $G$ and
a representation $G\to GL_n$, see \cite{donagi_spectral_covers}. Another feature is that
  even the generic  spectral curve may happen to be reducible, e.g., if  $G=SO(2n+1)$ 
(in which case one  uses  a particular irreducible component for setting up the spectral correspondence).
There is a   uniform approach (\cite{donagi_spectral_covers}), using more Lie theory, where one
replaces spectral curves with \emph{cameral curves}. These are (ramified) Galois covers with covering group
the Weyl group $W$ of $G$.
These curves come with an embedding in the vector bundle $\tot \ft\ctimes K_X$ and  spectral curves
can be obtained as  appropriate quotients thereof.
 This is the approach taken in  \cite{faltings}, \cite{donagi_spectral_covers},  \cite{hurtubise_kk}, \cite{scognamillo_elem}, \cite{don-gaits}.
The last reference also contains a discussion of ``abstract'' cameral covers, where the r\^ole of the coefficient bundle (here $K_X$) is separated from
that of the intrinsic spectral data.
 We are going to discuss cameral covers in Section \ref{cameral}.

For the group  $SL_2(\CC)$ the spectral and cameral covers coincide. In this case $\ft\simeq\CC$, $W\simeq \ZZ/2\ZZ$ and
$\cB_{SL_2}=H^0(X,K^2_X)$.
The double cover, associated with $b\in \cB$ is
\[
 \widetilde{X}_b= \scV\left( \lambda^{\otimes 2}- \pi^\ast b\right) \subset \tot K_X.
\]
It  can be   reducible (if $b=b_0^2$,  $b_0\in H^0(K_X)$), or
non-reduced (if $b=0$). It is easy to check that $\widetilde{X}_b$ is smooth precisely when $b\in H^0(K_X^2)$ has only simple zeros, 
and we assume this to be the case from now on. 
Then  $\pi_b$ is branched at  the $(4g-4)$  zeros
of $b$ and by the theorem of Riemann--Hurwitz its genus  is $\widetilde{g}=4g-3$. 
The spectral cover has
an involution $\sigma\in\maut(\widetilde{X}_b)$, induced by $\lambda\mapsto -\lambda$, i.e., 
 multiplication by $(-1)\in\CC$ in the vector bundle $\tot K_X$. The automorphism $\sigma$ generates an action of $W$  on $\widetilde{X}_b$,
which corresponds to interchanging the eigenvalues of the Higgs field.

As we have argued, the Hitchin fibre $h^{-1}(b)$ can be described as
\[
h^{-1}(b)\simeq \left\{\left. L\right| \det \pi_{b\ast} L=\cO_X \right\}= Nm_{\widetilde{X}_b/X}^{-1}(K_X) \subset \pic^{2g-2} \widetilde{X}_b.
 \]
 By the definition of the norm map,  $Nm^{-1}(\cO_X)$ consists of  line bundles $\cO_X(D)$ for which $\sigma^\ast D\equiv -D$.
Hence $h^{-1}(b)$ is a torsor over the Prym variety
\[
\textrm{Prym}_{\widetilde{X}_b/X}= \left\{\left. L\in\pic \widetilde{X}_b\right| \sigma^\ast L\simeq L^{-1}\right\}.
\]
We can think of the Prym variety in yet another way. The cocharacter lattice $\Lambda\simeq \ZZ$ of $SL_2$ carries
a natural $W$-action, namely the sign representation of $\ZZ/2\ZZ$, and so does
the cover $\widetilde{X}_b$. The latter action  induces
an action of $W$ on $\pic_{\widetilde{X}_b}$ by pullback of divisors.
The Prym variety is the set of \emph{invariant} elements for the combined $W$-action on  $\pic\widetilde{X}_b\otimes_\ZZ \Lambda$.
Up to isogeny, this turns out to be the correct description for arbitrary structure groups, as we shall briefly discuss in (\ref{gen_abel}).

      \section{Cameral covers}\label{cameral}
	    \subsection{Adjoint quotients}
Let us fix a complex simple Lie group $G$, with Lie algebra $\fg$, and consider  the coordinate ring $\CC[\fg]=\sym\ \fg^\vee$, on which
$G$ acts via the adjoint representation.
By a theorem of Chevalley, the algebra of invariants  $\CC[\fg]^G\subset \CC[\fg]$ is a free commutative (i.e., polynomial) algebra on $l$ generators:
$\CC[\fg]^G\simeq \CC[I_1,\ldots, I_l]$, where $I_j$ are homogeneous polynomials of degree $d_j=m_j+1$, and $l$ is the \emph{rank} of $\fg$, i.e., 
$l=\dim \ft$, where $\ft\subset \fg$ is (any) Cartan subalgebra. 
The generators of
the ring of invariants are not canonical in any way, but the set $\{d_j\}$ is  independent of the chosen $\{I_j\}$.

Recall that  $x\in \fg$ is called \emph{semisimple} (respectively, \emph{nilpotent})
if $\ad x\in \mend(\fg)$ is semisimple (respectively, nilpotent). 
It is \emph{regular}, if its  centraliser is of the smallest dimension possible, i.e.  $\dim \ker\ad x=  l$.
We use $\fg^{reg}$, $\fg^{ss}$, $\fg^{nlp}$ to denote the subsets of $\fg$, consisting of regular, semisimple and nilpotent
elements, respectively.
By Jordan decomposition \cite{jacobson_jordan}, every $x\in \fg$ has a unique representation as a sum $x=x^{ss}+x^{nlp}$, where
$x^{nlp}\in \fg^{nlp}$, $x^{ss}\in \fg^{ss}$ and $[x^{ss},x^{nlp}]=0$.

We consider now the different guises of the  ``adjoint quotient morphism''. The GIT quotient of $\fg$ under $G$ is, by definition, 
$\fg\sslash G=\spec \CC[\fg]^G$, and the inclusion $\CC[\fg]^G\subset \CC[\fg]$ corresponds to a morphism of
affine varieties $\chi: \fg\to \fg\sslash G$.
 A choice of $\{I_j\}$
identifies $\fg\sslash G$ with  $\CC^l$, and   $\chi(x)=  (I_1(x),\ldots I_l(x))$ in terms of this identification.
Choosing a basis  in $\ft$ (e.g. by fixing simple coroots) we identify the latter  with a morphism
$\chi: \CC^l\to\CC^l$, see Section (\ref{kostant}) and (\ref{G2}) for concrete examples.

We fix now  Cartan and Borel subgroups $T\subset B\subset G$, and denote by $W\subset GL(\ft)$ the Weyl group.
By another  result of Chevalley, the inclusion
$\ft\hookr \fg$  induces an algebra isomorphism $\CC[\fg]^G\simeq \CC[\ft]^W$,
\[
      \xymatrix{
		    \CC[\fg]^G\ar[d]_-{\simeq}\ar@{^{(}->}@<-0.5ex>[r]     &\CC[\fg]\ar@{->>}[d]^-{\textrm{res}}\\
		    \CC[\ft]^W  \ar@{^{(}->}@<-0.5ex>[r]                   &\CC[\ft]\\
},
\]
and hence $\ft/W\simeq \fg\sslash G$.
We can thus    interpret the adjoint quotient as the morphism of affine varieties
$\fg\to \ft/W$, corresponding to the algebra homomorphism $\CC[\ft]^W\simeq \CC[\fg]^G\subset \CC[\fg]$.
More concretely, any element of $\fg$ is $G$-conjugate to an element of $\ft$ that is determined uniquely up to $W$-conjugation.

Finally, we  recall that there is a bijection $\fg\sslash G\simeq \fg^{ss}/G$, corresponding to
  $I_j(x)=I_j(x^{ss})$.
All in all, given $x\in\fg$, we can describe $\chi(x)$ in any of the following ways:
\[
 \xymatrix{  \fg\sslash G        &\simeq & \fg^{ss}/G  & \simeq &\ft/W                              &\simeq &\CC^l\\
	     \overline{G\cdot x}&       &G\cdot x^{ss}&        &\left(G\cdot x^{ss}\right)\cap \ft &       &\underline{I}(x) \\
}.
\]

 A special r\^ole in what follows is played by the \emph{regular} elements of $\fg$.
For the moment we only say that
by \cite[Theorem 2]{kostant} there is  
a bijection $\fg^{reg}/G\simeq \fg\sslash G$. Moreover,  Theorem 7, \emph{ibid.} describes a section $\ft/W\to \fg^{reg}$, which will be discussed
in (\ref{kostant}), to be used also in (\ref{flow}).

	    \subsection{$K_X$-valued cameral covers}
The affine variety $\fg\sslash G$, while not a vector space in a canonical way,  carries a canonical $\CC^\times$-action. 
This action is
 induced by
the homothety action on $\fg$ by forcing $\chi$ to be $\CC^\times$-equivariant, i.e., by setting $t\cdot \chi(v):=\chi(tv)$.
 After a choice of $\{I_j\}$ this becomes the action $t\cdot (z_1,\ldots, z_l)= (t^{m_1+1}z_1,\ldots , t^{m_l+1}z_l)$.
Notice that there is a unique fixed point, the closure of the orbit of the origin, which makes $\fg\sslash G$ into a pointed space.
We can now consider  any principal $\CC^\times$-bundle and  form the associated bundle with fibre $\fg\sslash G$.
We apply this to $\tot K_X^\times$, 
 the canonical bundle $K_X$ with the zero section removed. 
Naturally, we can use the isomorphism (of varieties) $\fg\sslash G\simeq \ft/W$, and consider the $\CC^\times$-action on $\ft/W$ as
induced by homotheties in $\ft$. Either way, we are interested in the (fibre)bundle
\[
 \bdU = \left(\ft\ctimes K_X\right)/W\simeq K_X^{m_1+1}\oplus\ldots \oplus K_X^{m_l+1}.
\]
Its  space of global sections is the $G$-Hitchin base:
\[
 \cB_{G}= H^0(X,\bdU)= H^0(X,\ft\ctimes K_X/W)\simeq H^0(K_X^{m_1+1})\oplus \ldots \oplus H^0(K_X^{m_l+1}).
\]
We shall sometimes denote the base by $\cB_\fg$, since it  actually depends only on $\fg$, or equivalently, on $G^{ad}$, via the $\{d_j\}$.

The ramified $W$-cover $\ft\ctimes K_X\to \ft\ctimes K_X/W$ can be pulled back to $X$ along any section
$b:X\to \ft\ctimes K_X$ and hence
gives rise to a $\cB_\fg$-family of $W$-covers of $X$:
\[
 \xymatrix{\widetilde{X}_b\ar[d]^-{\pi_b}\ar@<-0.5ex>@{^{(}->}[r] &\widetilde{\cX}\ar[r]\ar[d]^-{\pi}&  \tot \ft\ctimes K_X\ar[d]\\
	   \{b\}\times X\ar@<-0.5ex>@{^{(}->}[r]  & \cB_\fg\times X\ar[r]^-{\textrm{ev}}&\tot (\ft\ctimes K_X)/W\\
 }.
\]
Here $\widetilde{X}_b$ is the cameral cover of $X$, corresponding to $b\in\cB_\fg$, while
 $\widetilde{\cX}$ is the \emph{universal cameral cover}.
Notice that both $\widetilde{X}_b$ and $\widetilde{\cX}$ come with an embedding in $\tot\ft\ctimes K_X$, and hence,
with a canonical $W$-action. Properly speaking, cameral covers of this kind are called 
\emph{$K_X$-valued} cameral covers.
Instead of $K_X$ we could have used an arbitrary line bundle (having sufficiently many sections) to host the cameral curve.
 In Sections \ref{dm_cubic}
and \ref{G2} we shall choose a sufficiently positive divisor $D$ on $X$ and consider $L:=K_X(D)$-valued cameral covers.
 This notion can be contrasted with the notion of an  ``abstract cameral cover'' of $X$, which is defined as
a $W$-cover that is \'{e}tale-locally the pullback of the $W$-cover $\ft\to \ft/W$. We recall that a $W$-cover is a finite flat morphism
$\pi:\widetilde{X}\to X$, such that $\pi_{\ast}\cO_{\widetilde{X}}$ is locally isomorphic to $\cO_X\otimes\CC[W]$. We refer to
\cite[\S 2.4]{don-gaits}  for more details.

Away from the ramification locus the covers $\widetilde{X}_b\to X$ and $\widetilde{\cX}\to \cB\times X$ are
$W$-Galois covers.
By repeated use of Bertini's theorem one can show that 
for generic $b\in\cB_\fg$ the cover $\widetilde{X}_b$ is non-singular, see \cite[\S 1]{scognamillo_elem}.

There is a natural  map
\[
 h: \Higgs_{G}\to \cB_{\fg}
\]
from the coarse moduli space of semi-stable $G$-Higgs bundles to the Hitchin base. It is a
 surjective proper  morphism with Lagrangian fibres. It is known that the components of $\Higgs$
 are indexed by the topological type of the principal bundle (see \cite[Lemma 4.2]{don-pan}
and \cite{oscar_andre_compts} for a different proof), i.e.
\begin{equation}\label{conn_comp}
 \Higgs_{G,X}=\coprod_{d\in\pi_1(G)}\Higgs_{G,X}^d 
\end{equation}
is a decomposition into connected components 
and the restrictions $h_d:\Higgs_{G,X}^d\to \cB_\fg$ are proper morphisms.

To construct $h$, notice that
for any $G$-bundle $P$, the quotient $\chi:\fg\to \fg\sslash G$ induces a map of associated bundles
$Ad(\chi): P\times_{Ad}\fg\to P\times_{Ad}(\fg\sslash G)\simeq P\times (\fg\sslash G)$. After twisting with $K_X$
and taking global sections, this gives a morphism of affine varieties $H^0(\ad P\otimes K_X)\to \cB_\fg$.
If $P$ is regularly stable, both the source and the target   are affine spaces of the same dimension
(equal also to $\dim \bun_{G,X}$) and, by Serre duality,
$H^0(\ad P\otimes K_X)=T^\vee_{\bun_G, [P]}$.
If $P$ is not stable, then 
$h^0(\ad P\otimes K_X)\geq \dim\cB_\fg$.
In that case, there is a Zariski-open (possibly empty) subset of $H^0(\ad P\otimes K_X)$,
corresponding to stable Higgs  structures on $P$.

The above construction can be  relativised. Indeed, given a complex manifold $S$ and a 
(holomorphic) $G$-bundle $\scP\to S\times X$, we obtain, exactly as a above, a polynomial
map $ H^0(S\times X, \ad\scP\otimes p_X^\ast K_X)\to H^0(S\times X, p_X^\ast \bdU)$. Since the
target vector space equals $H^0(S,\cO_S)\widehat{\otimes} \cB_\fg$,
  given an $S$-family of
Higgs bundles $(\scP,\Theta)$, $\Theta\in H^0(S\times X, \ad\scP\otimes p_X^\ast K_X)$,
one obtains a morphism $S\to\cB_\fg$. If $(\scP,\Theta)$  is a family of semi-stable Higgs bundles, this
morphism factors through the classifying map $S\to \Higgs_{G,X}$.

In Section (\ref{kostant})
we are going to consider  the Kostant section of $\chi$, the way it gives rise to a section 
of $h_0$, and the behaviour of the section under the Hamiltonian flow of linear functions on
$\cB_\fg$. The construction of the Hitchin section involves constructing a morphism
$\cB_\fg\to H^0(X,\ad P\otimes K_X)$ for an appropriately chosen unstable bundle $P$.

	    \subsection{General Abelianisation}\label{gen_abel}
It is beyond the scope of the current lectures to present a detailed account of abelianisation.
This is a beautiful and surprisingly intricate subject, beginning with \cite{hitchin_sd} and \cite{hitchin_sb}
where classical groups are treated via spectral curves, and continuing with works of many mathematicians,
such as \cite{donagi_decomposition}, \cite{faltings}, \cite{scognamillo_elem} and \cite{don-gaits}.
Here we just state the basic form of the result for the case of simple group $G$, which is all we shall need in the sequel.

Consider a point  $b\in \cB_G$, corresponding to a  non-singular cameral cover $\widetilde{X}_b$.
We identify, as usual, the lattice $\Lambda =\cchr_G\subset \ft_\RR$ with $\mhom(\CC^\times, T)$, and, 
consequently, $T\simeq \Lambda \otimes_\ZZ \CC^\times$.
The set of isomorphism classes of $W$-invariant $T$-bundles on $\widetilde{X}_b$ can then be identified as
$\left(\pic \widetilde{X}_b\otimes_\ZZ \Lambda\right)^W$, where $W$ acts both on $\pic \widetilde{X}_b$
(by pullback of divisors via the action on the cameral cover) and on $\Lambda$ by reflections.
The connected component containing the trivial $T$-bundle is in fact an abelian variety, a Prym-like variety.
The above set of $W$-invariant elements of $\pic \widetilde{X}_b\otimes_\ZZ \Lambda $
can also  be described as $H^1(X,\overline{\cT})$, where the sheaf of abelian groups $\overline{\cT}$ 
is defined as $\pi_{b\ast}^W\left(\Lambda\otimes_\ZZ\cO^\times_{\widetilde{X}_b}\right)$.

One can then construct (\cite{scognamillo_elem}) a morphism $h^{-1}(b)\to H^1(X,\overline{\cT})$, or,
if given also a point in $h_c^{-1}(b)$, a morphism from $h^{-1}_c(b)\to H^1(X,\overline{\cT})^0$.
The source and the target have the same dimension and the map is injective (for $G$-simple) if
$G\neq SO(2n+1)$, in which case the map has a finite kernel.
The complication is related to the presence of non-primitive coroots. One then  introduces
an appropriate subsheaf $\cT\subset \overline{\cT}$, consisting of sections, taking value $+1$ 
at all ramification points, and the generic Hitchin fibre is a torsor over the connected component
$H^1(X,\cT)^0$. For a description of the torsor and the spectral data for principal Higgs bundles
see \cite{don-gaits} and \cite[Appendix A]{don-pan}.
We shall stick to the notation of these references   and write $\Prym_{\widetilde{X}_b/X}$
for $H^1(X,\cT)$.
The relative Prym fibration (over $\scB_G\subset \cB_G$) will be denoted by $\Prym_{\widetilde{\cX}/\scB}$.
If $G$ is not of type ${\sf B}$ this is precisely 
$\left(\pic_{\widetilde{\cX}/\scB}\otimes_\ZZ \Lambda\right)^W$.

For topologically trivial Higgs bundles one can make a consistent choice of a point in the Hitchin fibre
(see Section \ref{flow}) and hence obtain a global identification 
$\left.\Higgs_{G,X}^0\right|_{\scB}\simeq \Prym^0_{\widetilde{\cX}/\scB}$. For other topological types this is possible
only locally on $\scB$.

      \section{Principal Subalgebras and Kostant's Section}\label{kostant}
We begin this section by  reviewing  some Lie-theoretic results, mostly due to B.Kostant.
	\subsection{Gradings}
Let $T\subset G$ be a Cartan subgroup and
\[
 \fg = \ft\oplus \bigoplus _{\alpha\in \cR}\fg_\alpha
\]
the corresponding root space decomposition. Here $\cR\subset \ft^\vee$ denotes the root system of $\fg$
 and $\fg_\alpha\subset \fg$ is the $\ft$-eigenspace with eigenvalue $\alpha$. That is,
$x\in \fg_\alpha$ if and only if 
\[
 [h,x]=\alpha(h)x
\]
for all $h\in\ft$.
A choice of Borel subgroup $B\supset T$ 
 determines a set of  simple (positive) roots $\Pi\subset \cR^+$ and hence   a grading of $\fg$ by ``height''. This is the unique
 Lie algebra grading
$\fg=\bigoplus_{m=-M}^M \fg_m$ satisfying
      \[
     \fg_0=\ft,\  \fg_{1}=\bigoplus_{\alpha\in \Pi}\fg_\alpha,\quad \fg_{- 1}=\bigoplus_{-\alpha\in \Pi}\fg_\alpha .
      \]
The natural number $M$  is the  \emph{Coxeter number} of $\fg$.

    \begin{Ex}\label{sl3_grading}
	Consider $G=SL_3(\CC)$ with Borel (respectively Cartan) subgroups, consisting of upper-triangular (respectively,  diagonal) elements of
$G$. Then the subspaces $\fg_m\subset \fs\fl_3(\CC)$, $-2\leq m\leq 2$ are 
\[
 \fg_m =  \left\{\left. (A_{ij})\in \fs\fl_3\right| A_{ij}=0 \textrm{ unless } j=i+m  \right\},
\]
 i.e., 
\[
\fg_{-2}=
	  \begin{pmatrix} 
		     0 &0 &0\\
		     0& 0 & 0\\
		     *& 0& 0\\
	  \end{pmatrix}, \quad
\fg_{-1}=
	  \begin{pmatrix} 
		     0 &0 &0\\
		     *& 0 & 0\\
		     0& *& 0\\
	  \end{pmatrix}, \quad
\fg_{0}=
	  \begin{pmatrix} 
		     * &0 &0\\
		     0& * & 0\\
		     0& 0& *\\
	  \end{pmatrix} \cap \ft,
\]
\[
\fg_{1}=
	  \begin{pmatrix} 
		     0 &* &0\\
		     0& 0 & *\\
		     0& 0& 0\\
	  \end{pmatrix}, \quad
\fg_{2}=
	  \begin{pmatrix} 
		     0 &0 &*\\
		     0& 0 & 0\\
		     0& 0& 0\\
	  \end{pmatrix}. 
\]
      \end{Ex}
Next,the height grading of $\fg$ induces   a Lie-algebra grading on $\mend(\fg)$
``by shift'', i.e., $\deg \phi=m$ if $\phi(\fg_k)\subset \fg_{m+k}$, for all $k$.
The adjoint representation  is compatible with these gradings, i.e.,
$\ad: \fg\to \mend (\fg)$
 is a homomorphism
of Lie algebras with grading, so $\ad(\fg_m)\subset \mend_m(\fg)$.
	\subsection{Principal subalgebras}
We recall now some standard properties of principal three-dimensional subalgebras.
Convenient references include Kostant's original papers (\cite{kostant_tds}, \cite{kostant}) and some later expositions, such as
   \cite{bei-drin},  \cite{hitchin_teich}, \cite{cg} and  \cite[Ch.VIII, \S 11]{bourbaki_lie_gr_alg}.

By an ``$\slt$-triple'' we mean a (non-zero) triple of elements 
$\{\cx,\cih,\cy\}$ of $\fg$,
which satisfy the relations
$[\cx, \cy]=\cih$, $[\cih, \cx]=2\cx$, $[\cih, \cy]=-2\cy$.
Thus we are considering
 subalgebras $\fa\subset \fg$,  isomorphic to $\slt$, together with a choice of 
``canonical'' generating set.
We note that the
 defining relations vary throughout the  references. The elements $\cx$ and $\cy$ are called the
nil-positive and nil-negative elements of the triple, while $\cih$ is its neutral element.

By considering the adjoint representation of $\fa=\CC\langle\cx, \cih, \cy\rangle$  on $\fg$ we see
that $\cx$ and $\cy$ are nilpotent and that
the subspace $\ker\ad_\cy\cap \textrm{Im }\ad_\cy$ of $\fg$ is  a 
 Lie subalgebra. It is in fact nilpotent, and $G_{\cy}:=\exp\left(\ker\ad_\cy\cap \textrm{Im }\ad_\cy\right)$
will stand for
the corresponding unipotent group.

Suppose we are given two triples,  $\{\cx,\cih,\cy\}$ and $\{\cx',\cih',\cy'\}$, spanning subalgebras
$\fa$ and $\fa'$, respectively. If  $\cy=\cy'$ and $\cih=\cih'$, then also
$\cx=\cx'$. 
If only $\cy=\cy'$ is known to hold, then $\cih-\cih'\in \ker\ad_\cy\cap\textrm{Im }\ad_\cy$.
In this case
the two triples are conjugate by some 
$g\in G_{\cy}$, i.e.,
 $\{\cx',\cih',\cy'\}= \{g\cdot \cx, g\cdot \cih, g\cdot\cy=\cy\}$.
The assignment $g\mapsto \{g\cdot\cx, g\cdot\cih, \cy\}$ establishes a one-to-one correspondence
between $G_\cy$ and the set of triples, containing $\cy$ as nil-negative element. Finally, the orbit
$G_\cy\cdot\cih$ is identified with the affine space $\cih+ \ker\ad_\cy\cap \textrm{Im }\ad_\cy$.
 For proofs, see 
 \cite[Theorem 3.6]{kostant_tds} or \cite[Ch.VIII, \S 11, Lemma 4]{bourbaki_lie_gr_alg}.
In  general,  the two triples are conjugate by some $g\in G^{ad}$ if and only
if $\fa'=g\cdot \fa$  if and only if $\cy'=g\cdot \cy$, see  Proposition 1, \emph{ibid.}.
Consequently, 
the assignment $\{\cy,\cih,\cx\}\mapsto \cy$ induces
 an injective map from the set of conjugacy 
classes of $\slt$-triples (or subalgebras)  to $\left(\fg^{nlp}\backslash\{0\}\right)/G^{ad}$, the set of conjugacy classes 
of non-zero nilpotents.
This map is in fact a bijection, since
every non-zero nilpotent element of $\fg$ can be completed to  an
$\slt$-triple, 
by the Jacobson--Morozov theorem (\cite[Theorem 3.4]{kostant_tds}).
To establish the correspondence (on the level of conjugacy classes) one could, in fact, use any (non-zero) 
nilpotent of the $\slt$-subalgebra, see Corollary 3.7, \emph{ibid.}.

The adjoint representation of 
an $\slt$-subalgebra  $\fa\subset \fg$ decomposes  $\fg$ into a sum of irreducible representations,
whose number is \emph{at least} $l=\rk(\fg)$ (Theorem 5.2, \emph{ibid.}).
The subalgebra $\fa$ is called \emph{principal} if $\fa\backslash \{0\}\subset \fg^{reg}$. 
The  minimal number of summands ($l$) is achieved precisely when
$\fa$ is principal, and in that case they are all odd-dimensional
(Corollary 5.2, \emph{ibid.}).
The dimensions of the summands are determined by the exponents of $\fg$ and so one has
\begin{equation}\label{decomp}
\fg = \bigoplus_{i=1}^{l} W_{m_i},
\end{equation}
where $W_{m_i}\simeq \sym^{2m_i}(\CC^2)$.
On each $W_{m_i}$ the eigenvalues of $\ad {\cih}$ are even integers $2m$, where $-m_i\leq m \leq m_i$.  
The centraliser  $\fz_\cx=\ker\ad \cx=\fg_\cx$ is spanned by the primitive (highest weight) vectors in the $W_{m_i}$'s.
The principal subalgebra $\fa$ appears in this decomposition as $W_1$.

      \begin{Ex}\label{sl3_decomp}  Consider $\fg=\fs\fl_3(\CC)$ and the $\fs\fl_2(\CC)$-triple 
\[
 \left\{\cx=
 \left(\begin{array}{rrr}
                0 & 2 & 0\\
                0 & 0 & 2\\
                0 & 0 & 0\\
        \end{array}\right), \cih=
 \left( \begin{array}{rrr}
                2 & 0 & 0\\
                0  & 0 & 0\\
                0  & 0 & -2\\
                \end{array}\right),
\cy= \left(\begin{array}{rrr}
            0 & 0 & 0\\
            1 & 0 & 0\\
            0 & 1 & 0\\
            \end{array}\right)\right\}.
\]
 The induced   decomposition is
\[
 \fs\fl_3(\CC)\simeq W_1\oplus W_2,
\]
where
$W_1 =\fa = \textrm{span}\{\cx,\cih,\cy\}$
and
$W_2$ is the span of the matrices
\[ 
   \left( \begin{array}{rrr}
                                        0&0&1\\
                                        0&0&0\\
                                        0&0&0\\
                                        \end{array}\right),
    \left( \begin{array}{rrr}
          0&-1&0\\
          0& 0& 1\\
          0& 0& 0\\
          \end{array}\right),
   \left( \begin{array}{rrr}
          1 &0&0\\
          0 & -2& 0\\
          0& 0& 1\\
          \end{array}\right),
 \left( \begin{array}{rrr}
                        0 & 0 & 0\\
                        1 & 0 & 0\\
                        0 & -1 & 0\\
                \end{array}\right),
 \left( \begin{array}{rrr}
                        0 & 0 & 0\\
                        0 & 0 & 0\\
                        1 & 0 & 0\\
                \end{array}\right) 
\]
      \end{Ex}

Principal $\slt$-subalgebras play a special r\^{o}le in what follows.
Under the correspondence between (conjugacy classes of) $\slt$-subalgebras and nilpotents,
 principal subalgebras are identified with (conjugacy classes of)
 regular nilpotent elements.
      \begin{Ex} Let $\fg=\slt$. Then 
\[
 \xymatrix@1{\slt^{nlp} =
 \chi^{-1}(0)\subset \slt\ar[r]^-{\chi= \det} &\CC}
\]
is a singular affine quadric (cone) in $\CC^3$.
 The complement of the  tip of the cone is
the set of regular nilpotents  $\slt^{reg,nlp}= \slt^{nlp}\backslash\{0\}$.
Any two regular nilpotents are conjugate.
      \end{Ex}
It turns out that the above example is indicative for the situation in general:
 ``most'' of the $\slt$-subalgebras of a simple Lie algebra
$\fg$ are principal, and any two principal subalgebras  are conjugate. Moreover, there is a ``preferred'' choice of a principal subalgebra,
associated with any choice of simple positive roots. We outline the argument in the paragraphs below.

Let $B\supset T$ be a Borel subgroup and $\Pi\subset\cR^+$ the corresponding choice of simple positive roots.
Recall that the (closed) Weyl chamber $D\subset \ft_\RR$ is the set of all $\fii\in\ft_\RR$ which satisfy
 the inequality $\alpha(\fii)\geq 0$ for all $\alpha\in \Pi$.
If $\psi\in \ft$ is semi-simple and $\ad\psi$ has real eigenvalues, then $\psi$ is ($G^{ad}$-)conjugate to a unique element in $D$.
Hence the question of identifying the conjugacy classes of (principal) $\slt$-subalgebras can be split into two:
identifying the $\slt$-triples containing a fixed $\cih\in \ft$ and identifying those
$\cih\in D\subset\ft_\RR$, which are contained in an $\slt$-triple.

The choice of semi-simple  $\cih\in\ft$
  gives rise to a Lie-algebra grading
$\fg=\bigoplus_m \fg_m$,
where $\fg_m$ is the $m$-th eigenspace of $\ad\cih$.
In particular,
$\fg_0=\zeta(\cih)$ is the centraliser of $\cih$ in $\fg$ and $\exp \fg_0= Z_{G^{ad}}(\cih)$ is its centraliser
 in $G^{ad}$.

Consider the morphism $\fg_{-1}\to \mhom_\CC(\fg_0,\fg_{-1})$, $f\mapsto \left. \ad f\right|_{\fg_0}$ and
let $\widehat{\fg_{-1}}\subset \fg_{-1}$ be the preimage of the set of surjective linear maps.
This is a connected, dense and (Zariski) open set (\emph{ibid.}, Lemma 4.2B).
By Theorem 4.2, \emph{ibid.}, any two triples, containing $\cih$ as a neutral element are
conjugate. Moreover, the set of such triples is a $Z_{G^{ad}}(\cih)$-torsor and the choice of one  such a triple,
say $\{\cx,\cih,\cy\}$, identifies this torsor with  $\widehat{\fg_{-1}}= Z_{G^{ad}}(\cih)\cdot \cy$.

      \begin{Ex}
Let $\fg=\fs\fl_3(\CC)$ and $\cih=\textrm{diag} (2,0,-2)$. The induced grading  of $\fg$ coincides with the one given
in Example \ref{sl3_grading} and $\widehat{\fg_{-1}}=\CC^\times e_{21}\times \CC^\times e_{32}$.
      \end{Ex}

We now turn to the question of identifying the semi-simple elements in $D\subset\ft_\RR$ which can be completed
to an $\slt$-triple. As discovered by Dynkin (\cite{dynkin_subalgebras}, Theorem 8.2) this condition is quite restrictive.

 Recall that we can associate with $G$ the lattices
\begin{equation}\label{root_char}
 \rts\subset \chr_G\subset \wts \subset \ft^\vee_{\RR}
\end{equation}
\begin{equation}\label{coroot_cochar}
 \crts\subset \cchr_G\subset \cwts\subset \ft_\RR.
\end{equation}
Here $\rts = \bigoplus_{\alpha_i\in \Pi} \ZZ\alpha_i$ and $\wts\subset \ft^\vee_\RR$ is defined as the set of all
$\beta$ which satisfy   $( \alpha,\beta) = 2\frac{\langle\alpha,\beta\rangle}{\langle \alpha,\alpha\rangle}\in \ZZ$
for all $\alpha\in \rts$,
 $\langle\ ,\ \rangle$ being the Killing form.
Explicitly, $\wts = \bigoplus \ZZ\omega_i$ and the fundamental weights
 $\{\omega_i\}$ are related to the simple roots by
 $\alpha_i = \sum_j N_{ji}\omega_j$, where $N=(N_{ij})$ is the 
Cartan matrix.
 The lattice $\crts= \wts^\vee$ is generated by the 
simple coroots $\check{\alpha}_i = 2\frac{\langle\alpha_i,\ \rangle}{\langle\alpha_i,\alpha_i\rangle}$, 
and $\cwts =\rts^\vee$ is generated by the fundamental coweights $\{\eps_i\}$. 
The pairs $\left( \{\alpha_i\}, \{\eps_i\}\right)$ and
$\left(\{\omega_i\}, \{\check{\alpha}_i\}\right)$ are pairs of dual bases.
The character and cocharacter lattices can be identified with (the differentials of) elements of
$\mhom(T,\CC^\times)$ and $\mhom(\CC^\times,T)$. Similar identifications, involving  $T^{sc}$
and $T^{ad}$, can be given for the other four lattices.
The three top lattices determine the group $G$ up to isomorphism. 
The bottom row determines the (isomorphism class of the) Langlands dual group $^LG$.

Now suppose $\cih=\sum_{i=1}^l a_i \epsilon_i$, $a_i\in \ZZ$. If $\cih\in D$, then
necessarily $a_i\in\NN$. If there exists $\cy\in \ft$, such that $[\cih,\cy]=-2\cy$, then
$\cy$ has only components involving negative root vectors. If $e_{\alpha_i}$ is a root vector, 
corresponding to $\alpha_i\in \Pi$, then $\ad_\cih([\cy,e_{\alpha_i}])=(a_i-2)[\cy,e_{\alpha_i}]$.
Thus $a_i>2$ would violate the simplicity of $\alpha_i$.
Hence $a_i\in \{0,1,2\}$ and there are \emph{at most} $3^l-1$ non-trivial choices of $\cih\in D$.

The choices of $\cih\in D$ that actually \emph{do occur} are significantly fewer than $3^l-1$ and
 were determined by Dynkin for each of the simple Lie algebras,
see \cite[ \S 9]{dynkin_subalgebras}, Tables 13 and 15.
Moreover, these $\cih\in D$
determine the conjugacy classes of $\slt$-subalgebras 
(\cite{dynkin_subalgebras},  \cite[Lemma 5.1]{kostant_tds}).

There is one choice of $\cih\in D$ which occurs always, i.e., is allowed for all simple $\fg$, namely
   $a_1=\ldots=a_l=2$. The semi-simple element arising in this way is known
as ``twice the dual Weyl vector'':
\[
 \cih = 2\check\rho = 2\sum_{i=1}^l\eps_i = \sum_{i=1}^l r_i\check\alpha_i\subset \cchr_G.
\]
The integers $r_i$ (``Kac labels'') are found by inverting the Cartan matrix, see  \cite[\S 13]{humphreys}.

We now describe an $\slt$-triple, containing $\cih$.
First we need to fix root vectors in the $\fg_\alpha$, $\alpha\in \pm \Pi$.
In fact, only half of these need to be fixed and we choose
$\{f_{\alpha_i}\in \fg_{-\alpha_i}^\times,\ \alpha_i\in \Pi\}$.
 These uniquely determine $e_{\alpha_i}\in \fg_{\alpha_i}^\times$ by the requirement that
$\{ f_{\alpha_i}, e_{\alpha_i}, \check\alpha_i = [e_{\alpha_i}, f_{\alpha_i}]\}$ be an $\slt$-triple.
Here, given a root $\alpha\in\cR$ we write $\fg_\alpha^\times$ for the set  $\fg_\alpha\backslash\{0\}\simeq\CC^\times$.
    \begin{defn}
Let $G$ be a simple complex Lie group and $T\subset G$ a Cartan subgroup. By 
\emph{\'epinglage data} (compatible with $T$) we shall mean a choice of Borel subgroup $B$ (such that
     $T\subset B\subset G$), together with a choice of negative root vectors
$\{f_{\alpha_i}\in \fg_{-\alpha_i}^\times,\ \alpha_i\in \Pi\}$.
    \end{defn}
Note that usually the definition of \'epinglage involves a choice of \emph{positive} root vectors, but the above convention is better adapted
to our intended applications in the next section.  

One then checks immediately that
\begin{equation}\label{epinglage_triple}
 \left\{\cy=\sum_{i=1}^l f_{\alpha_i},\ \cih = 2\check\rho, \  \cx=\sum_{i=1}^l r_i e_{\alpha_i} \right\}.
\end{equation}
is an $\slt$-triple, which we shall call \emph{the $\slt$-triple associated to a choice of \'epinglage data}.
The $\slt$-triple in Example \ref{sl3_decomp} is an example of such a triple.

We now apply the previously stated results of \cite[\S 4]{kostant_tds} to this triple.
First, the centraliser $Z_{G^{ad}}(\cih)= T^{ad}= T/Z(G)\simeq \left(\CC^\times\right)^l$.
Moreover,
\[
\widehat{\fg_{-1}}=  Z_{G^{ad}}(\cih)\cdot \cy= \bigoplus_{\alpha\in\Pi} \fg_{-\alpha}^\times \subset \fg_{-1}
\]
and
\[
 B_{-}\cdot\cy:= \exp\left(\bigoplus_{m\leq -1} \fg_m\right)\cdot \cy = \widehat{\fg_{-1}}\times \bigoplus_{m= -2}^{M} \fg_m .
\]

By Theorem 5.3 \emph{ibid.}, the latter set is precisely the set of regular (principal) nilpotents in the
(negative) nilpotent Lie algebra $\bigoplus_{m\leq -1}\fg_m$, that is, 
\[
    \left(\bigoplus_{\alpha\in \cR^-}\fg_{\alpha}  \right)\cap \fg^{reg} =  \bigoplus_{\alpha\in\Pi} \fg_{-\alpha}^\times \times \bigoplus_{m= -2}^{M} \fg_m .
\]
Explicitly, these are the elements of $\bigoplus_{m\leq -1}\fg_m$, which have non-zero components along \emph{each} negative simple
root space.
      \begin{Ex} Let $\fg=\fs\fl_3(\CC)$. Then
\[
 \left( \fg_{-2}\oplus \fg_{-1}\right)\cap \fg^{reg} = \left\{\left. \begin{pmatrix}
									0&0&0\\
									a&0&0\\
									c& b&0\\
								      \end{pmatrix} \right| ab\neq 0
 \right\}.
\]
      \end{Ex}
Finally, by Corollary 5.4, \emph{ibid.}, \emph{any} nilpotent element of $\fg$ is $G^{ad}$-conjugate
to a nilpotent element in $\bigoplus_{m\leq -1}\fg_m$.
Hence, by
Corollary 5.5, 
regular nilpotents $\fg^{reg,nlp}$  form a single $G^{ad}$-orbit,  which is a connected, dense open subset
$\fg^{reg,nlp}\subset \fg^{nlp}=\chi^{-1}(0)$. 

In view of the correspondence between (conjugacy classes of) $\slt$-subalgebras and
nilpotents elements, the above implies that
 all principal subalgebras  are conjugate under $G^{ad}$, and the choice of
\'epinglage data singles out a representative of the conjugacy class.

Let us  mention in conclusion there was no real need to fix
  the \'epinglage data.
Instead, one could just fix
$B\supset T$ and consider $\slt$-subalgebras (or triples) which are ``compatible'' with
this choice and contain $\cih=2\check{\rho}$ as a neutral element.
Then the set of principal compatible triples  is just a
  $(\CC^\times)^l$-torsor, $\bigoplus_{\alpha\in\Pi} \fg_{-\alpha}^\times$.  
Fixing the \'epanglage trivialises the torsor and identifies the compatible principal triples with
$\left\{\sum c_i f_{\alpha_i}, 2\check{\rho}, \sum c_i^{-1}r_i e_{\alpha_i}\right\}$, $c_i\in\CC^\times$.

	    \subsection{The Kostant section}
As in Section \ref{cameral}, let us denote by $\chi:\fg\to \ft/W\simeq \fg\sslash G$ the adjoint quotient
morphism, and let
 $\fa=\{\cy,\cih,\cx \}_\CC$ be a principal 3-dimensional subalgebra of $\fg$.
 In Theorem 7 from \cite{kostant} it is proved  that
the restriction of $\chi$ to the affine space $\{\cy\}+\ker\ad\cx$
is an isomorphism (of algebraic varieties) $\left. \chi\right|_{\cy+ \fz(\cx)}$.
We thus have a morphism $\bk:\ft/W\to \fg$, \emph{the Kostant section}, such that
$\chi\circ \bk=id$.

The choice of $\fa$ determines  a (vector space) bigrading of $\fg$: 
    \[
     \fg = \bigoplus_{k,i} \fg_{k,i}, \ \fg_{k,i}=\fg_k\cap W_i,
    \]
and in this case the centralisers of $\cx$ and $\cy$ are  the sums of highest, respectively lowest weight spaces:
    \[
     \fz(\cx)=\ker\ad \cx = \bigoplus_{i=1}^{l} \fg_{m_i,i},\ \fz(\cy)=\ker\ad \cy = \bigoplus_{i=1}^{l} \fg_{-m_i,i}.
    \]
As is clear from the preceding discussion, here $\dim \fg_{\pm m_i,i}=1$.
Thus the image of the Kostant section, also known as \emph{the Kostant slice}, is
\[
 \textrm{Im } \bk = \{\cy\}+ \fz(\cx) =\{\cy\}+ \bigoplus_{i=1}^l \fg_{m_i,i} \subset \fg^{reg}.
\]

    \begin{Ex}\label{sl3_kostant}
	Consider $\fg=\fs\fl_3(\CC)$ with the principal subalgebra from Example \ref{sl3_decomp}.
 Let $\chi: \fs\fl_3(\CC)\to \CC^2\simeq \ft/W$ be defined as
$\chi(A)=(a_1(A),a_2(A))$, where $\det(A-\lambda 1)=-\lambda^3 + \lambda a_1(A) + a_2(A)$.
%
Then $\bk: \CC^2\to \fs\fl_3(\CC)^{reg}$ is defined by
	\[
	 \bk(a,b)=    \begin{pmatrix}
			    0 &a/2 & b\\
			    1  & 0  & a/2\\
			    0& 1& 0 \\
		        \end{pmatrix}.
	\]
      \end{Ex}

Kostant's   Theorem 7 actually says a bit more. 
If one insists on considering $\chi$ as a morphism $\fg\to\CC^l$, one must choose
generators $I_j$ of $\CC[\fg]^G$.
While there is freedom in the choice of generators, 
  their behaviour on the slice is constrained.
  The second part of Theorem 7 states  that there is always a choice of highest weight vectors $v_i\in W_i\cap\fg_{m_i}$, 
such that
$I_j(\cy +\sum_{i=1}^l a_i v_i)= a_j + p_j(a_1,\ldots, a_{j-1})$, for some polynomials $p_j$ without constant term,
and similarly for the inverse map.
In particular, if we fix \emph{some} collection of $v_i$, then  $\{I_j\}$ can always be chosen in a way that $I_j(\cy +\sum_{i=1}^l a_i v_i)= a_j$.
This was already hinted at in the above example, where the factors of $1/2$ are carefully chosen.

	    \subsection{Principal homomorphisms and $\CC^\times$-actions}

Let  $\{\cy,\cih, \cx\}$ be a distinguished principal $\slt$-triple, compatible with a choice of
$T\subset B\subset G$. By mapping $\left\{\begin{pmatrix}
                                           0&0\\ 1&0\\
                                          \end{pmatrix},
					  \begin{pmatrix}
					   1&0\\0&-1\\
					  \end{pmatrix}, 
					   \begin{pmatrix}
					    0&1\\0&0\\
					   \end{pmatrix}
\right\}$ to $\{\cy,\cih, \cx\}$ we obtain a homomorphism $\slt\to \fg$.
 Let $\varrho: SL_2(\CC)\to G^{sc}$ be the corresponding ``principal homomorphism'' and
$\overline{\varrho}$ the composition $SL_2(\CC)\to G^{sc}\to G^{ad}\subset GL(\fg)$.
Consequently, the $1$-parameter subgroup
  $\CC^\times\subset  SL_2$, $t\mapsto \textrm{diag}(t,t^{-1})$
 gives rise to a $\CC^\times$-action on $\fg$, having weight $(2m)$ on the subspace $\fg_{m}\subset \fg$.
Then the Kostant slice  $\bdk(\ft/W)$ is preserved under the ``shifted'' action $t^2\overline{\varrho(t^{-1})}$:
\[
 t^2 \overline{\varrho(t^{-1})} \cdot \left( \{\cy\}+\fz(\cx) \right)= \{\cy\} +  t^2 \overline{\varrho(t^{-1})}\cdot \fz(\cx).
\]

	\section{The Hitchin section and its flow}\label{flow}
	\subsection{The section}
 In \cite{hitchin_teich} Hitchin ``promoted'' Kostant's section to a section of the morphism
$h_0:\Higgs_{G,X}^0\to \cB_\fg$ (see Section \ref{cameral}).
Hitchin provided not just an  existence statement, but an actual construction, 
depending on the choice of a principal subalgebra and a theta-characteristic on $X$.
 We shall review briefly this construction below.
Apart from  \cite{hitchin_teich}, my understanding of this topic has been largely enhanced by the expositions in 
\cite{ngo} and \cite{don-pan}.

  Recall first that a \emph{theta-characteristic},
or equivalently, a \emph{spin structure}, on $X$ is a pair $(\zeta,i)$, where $\zeta$ is a line bundle
and $i\in \mhom(\zeta^{\otimes 2}, K_X)$ is an isomorphism. 

Due to the divisibility of $\pic X$ spin structures exist for any $X$.
The degree of such a line bundle $\zeta$ must be $g-1$ and 
there are $2^{2g}$ choices of its isomorphism class. These classes are identified with  the fibre $sq^{-1}(K_X)$
of the squaring morphism $sq: \pic X\to \pic X$,
$sq([L])=  [L]^{\otimes 2}$.
For each chosen $\zeta$ there is a $\CC^\times$-worth of  choices of $i$.
For the most part we shall  simply  write 
$K_X^{1/2}$ instead of $(\zeta,i)$.
Notice that $i$ induces a canonical morphism $1: K_X^{1/2}\to K_X^{-1/2}\otimes K_X$.
	    \begin{thm}(\cite{hitchin_teich})
A choice of principal $\slt$-triple $\{\cx,\cih, \cy\}$ in $\fg$ and a theta-characteristic $K_X^{1/2}$ determines a Lagrangian section
$\fv$ of the restriction $h_0$ of the Hitchin map to the neutral connected component of the (coarse) moduli space of semi-stable
$K_X$-valued $G$-Higgs bundles on $X$:
\[
 \xymatrix@1{  \Higgs^0_{G,X}\ar[r]^-{h_0}&\cB_\fg\ar@/^1.5pc/[l]^-{\fv} }.
\]
	    \end{thm}
For some groups  (such as $G=SL_{2n+1}(\CC)$) one can construct a section without choosing $K_X^{1/2}$.
We recall also that $\Higgs_{G,X} = \coprod_{d\in\pi_1(G)}\Higgs^{d}_{G,X} $ and
the non-neutral connected   components are $\Higgs_{G,X}^0$-torsors which may not admit   global sections
of $h_d$.

We start with the basic special case $G=SL_2(\CC)$, which was considered in detail already in \cite{hitchin_sd}.
The section gives rise to a collection  of rank-2 Higgs (vector) bundles indexed by
 $\cB_{\slt}=H^0(X,K^2_X)$, which is
    \begin{equation}\label{section_sl2}
\left\{\left(V_b,\fii_b\right)\right\}_{b\in\cB}=  \left\{ \left(K_X^{1/2}\oplus K_X^{-1/2},\begin{pmatrix}
                                                     0& b\\
						     1& 0\\
                                                    \end{pmatrix}
 \right)    \right\}_{b\in\cB}.
    \end{equation}
For $b=0$ this is a ``Toda Higgs bundle'', whose  Hermite--Yang--Mills metric 
is induced by the unique metric on $X$, descending from the constant negative curvature metric
on the unit disk. Moreover, deformations of this bundle can be tied with deformations of the underlying curve $X$.
This is largely the subject of \cite[\S 10]{hitchin_sd}  and 
is central  in the subject of ``higher Teichm\"uller theory''.

    A family of Higgs bundles parametrised by $\cB_\fg$ (e.g., the section $\fv$) 
is not just a collection
 $\left\{\left(V_b,\fii_b\right)\right\}_{b\in\cB}$, as in (\ref{section_sl2}):
one must also specify how the members of the family   fit together.  I.e., we must exhibit  a $p_X^\ast K_X$-valued $G$-Higgs bundle
on $\cB_{\fg}\times X$, which restricts to $(V_b,\fii_b)$ on $\{b\}\times X$.
 We outline the construction
of such a  $p_X^\ast K_X$-valued Higgs bundle below.

Any  choice of $T\subset B\subset  G$ determines 
 a natural $GL(\fg)$-bundle (rank $\dim G$-vector bundle) on $X$, namely
\begin{equation}\label{Lie_bundle}
 \bdE= \bigoplus_{m=-M}^M \fg_m\ctimes K_X^m.
\end{equation}

 The natural linear maps
\[
 \xymatrix@1{\fg_{-1\ }\ar@{^{(}->}@<-0.5ex>[r]^-{\ad} &\mend_{-1}(\fg)\simeq \mhom_{-1}(\fg,\fg)\simeq \mhom_{-1}(\bdE,\bdE\otimes K_X) }
\]
allow us to identify any   element of $\fg_{-1}\subset \fg$ with  a Higgs field on  $\bdE$.

We  fix next a theta-characteristic $K_X^{1/2}=(\zeta,i)$ and  
a principal
triple $\{\cy,\cih,\cx\}$,   compatible with $T\subset B$.
Let $P_0$ be the frame bundle  $\misom_{\cO_X} \left(K_X^{1/2}\oplus K_X^{-1/2},\cO_X^{\oplus 2}\right)$
and 
$\bdP=P_0\times_{\varrho} G^{sc}$ the $G^{sc}$-bundle  associated to it via $\varrho$.
We then have a vector bundle isomorphism  $\bdE\simeq \ad \bdP = P_0\times_{\overline{\varrho}}\fg$.

  Consider  (as in Section \ref{cameral})  the bundle
$\bdU=\ft\ctimes K_X/W$ and
let $u:\tot \bdU \to X$ be the   bundle projection.
The choice of $\{\cy,\cih,\cx\}$  fixes a
Kostant section $\bk: \ft/W\hookr \fg$ and hence
an injection of vector (cone) bundles
  $\bdk:\tot\bdU\hookr \tot \bdE\otimes K_X$.
In turn, the latter determines a global  section $\bph_{\bdk}$ of $u^\ast\left(\bdE  \otimes K_X\right)$,
and thus   a $u^\ast K_X$-valued $G$-Higgs bundle  $\left(u^\ast\bdP, \bph_{\bdk} \right)$ on $\tot \bdU$.
Pulling  it to $\cB_\fg\times X$ by the evaluation morphism gives rise to a family $(p_X^\ast\bdP, \textrm{ev}^\ast\bph_{\bdk})$
of $K_X$-valued $G$-Higgs bundles, parametrised
by $\cB_\fg$:
\[
 \xymatrix{ p_X^\ast\left(\bdE\otimes K_X \right)\ar[d]\ar[r]  &  u^\ast\left(\bdE\otimes K_X\right)\ar[d]\ar[r] &\bdE\otimes K_X\ar[d]  \\
	    \cB_\fg\times X\ar[r]^-{\textrm{ev}}&\tot \bdU\ar@/^/[u]^-{\bph_{\bdk}} \ar[r]^-{u} & X\\ } .                                                                                                                                           	
\]
Hitchin in \cite{hitchin_teich} showed that this 
is a family of semistable Higgs bundles. By  the coarse moduli space property, there is a classifying map
$\cB_\fg\to \Higgs^0_{G,X}$, whose composition with $h_0$ is the identity.

Here is a more concrete description. Recall that in Section \ref{cameral} we defined, for a principal bundle $P$, a morphism
$H^0(X,\ad P\otimes K_X)\to\cB$, induced by $\chi$. 
Above we have constructed a particular bundle, $\bdP$, with $\bdE=\ad\bdP$ (see \ref{Lie_bundle}), together with an  embedding
of the base $\cB_\fg$ 
\[
 \bdk(\cB_\fg)= \{\cy\} + \bigoplus_{i=1}^l \fg_{m_i,i}\otimes H^0(X,K_X^{m_i+1})\subset H^0(X,\bdE\otimes K_X), 
\]
in a way that $h_0:\bdk(\cB_\fg)\to \cB_\fg$ is a linear isomorphism. 
The family $\fv(\cB_\fg)$  is obtained by
restricting to $\bdk(\cB_\fg)$ the tautological family  of Higgs structures on $\bdP$, parametrised by $H^0(\bdE\otimes K_X)$.

After choosing \'epinglage data, i.e., isomorphisms
$\fg_{-m_i,i}\simeq \CC$, we obtain bases of $\fg_{m_i,i}\simeq \CC$ by applying appropriate powers of $\ad\cih$.
 By \cite[Theorem 7]{kostant}, for any such choice
one can  choose  $G$-invariant polynomials $\{I_j\}$ so that
 $h_0: \bdk(\cB_\fg)\simeq \cB_\fg$ is identified with translation by $-\cy$, followed by the induced linear map
\[
 \bigoplus_{i=1}^l \fg_{m_i,i}\otimes H^0(X,K_X^{m_i+1}) \simeq \bigoplus_{i=1}^{l} H^0(X,K_X^{m_i+1}).
\]
	    \begin{Ex}\label{sl3_hitchin}
Consider $G=SL_3(\CC)$ with the standard Borel and Cartan subgroups. The exponents are $m_1=1$, $m_2=2=M$
and
\[
\bdE = \bigoplus_{m=-2}^2 \fg_{m}\ctimes K_X^{m},
\]
where the grading $\fg=\oplus \fg_m$ is described in
 Example \ref{sl3_grading}. Consequently,
\[
H^0(\bdE\otimes K_X)= \bigoplus_{m=-1}^2 \fg_{m}\ctimes H^0(K_X^{m+1})\]
is isomorphic to $\CC^2\oplus H^0(K_X)^{\oplus 2}\oplus H^0(K_X^2)^{\oplus 2}\oplus H^0(K_X^3)$.

 Choosing the distinguished principal $\slt$-triple as in Example \ref{sl3_decomp}, we obtain for the highest weight spaces
\[
 \fg_{m_1,1} = \CC\cdot \left(\begin{array}{rrr}
                0 & 1 & 0\\
                0 & 0 & 1\\
                0 & 0 & 0\\
        \end{array}\right),\ 
\fg_{m_2,2}= 
\CC\cdot \left(\begin{array}{rrr}
                0 & 0 & 1\\
                0 & 0 & 0\\
                0 & 0 & 0\\
        \end{array}\right).
\]
Choosing $\{I_1,I_2\}$ as in Example \ref{sl3_kostant}	 we have an  $\cB_{\fs\fl_3(\CC)}\simeq H^0(K_X^2)\oplus H^0(K_X^3)$
and $\bdk:\cB_{\fs\fl_3(\CC)}\to H^0(\bdE\otimes K_X)$ is identified with
	\[
	 (a,b)\longmapsto    \begin{pmatrix}
			    0 &a/2 & b\\
			    1  & 0  & a/2\\
			    0& 1& 0 \\
		        \end{pmatrix}.
	\]
	    \end{Ex}

	\subsection{The flow}\label{flow_flow}
Given a choice of principal subalgebra and a theta-characteristic, we can describe the Hitchin section
very concretely. A natural question, then, is to try to explicate the evolution of  the section  under the hamiltonian flow
of linear functions on the base.

We shall outline here the setup and the main ingredients of the construction. More details and some applications  can be found
in \cite{thesis} and the preprint \cite{symplectic_kur}. A more extensive treatment of this topic shall be given at another occasion.

In this subsection we assume that a set of generators $\{I_i\}$ of $\CC[\fg]^G$ has been fixed, so $\cB= H^0(X, \bU)$
will be given a  vector space structure. We focus our attention on the non-singular locus $\Higgs_{G,X}^{0,reg}$ of the neutral connected
component $\Higgs_{G,X}^0$. We shall write  simply  $\Higgs^{0,reg}$ whenever there is no danger of confusion.
We emphasise that $\Higgs^{0,reg}$ contains both non-singular Hitchin fibres and smooth loci of singular fibres, and hence
 $h_0^{-1}(\scB)=\Higgs^0\backslash h_0^{-1}(\Delta)\subsetneq \Higgs^{0,reg}$.
We shall also write $h_0$  or $h$ instead of $\left. h_0\right|_{\Higgs^{0,reg}}$ for better readability.

Let us discuss some consequences of having to deal  with the non-singular  locus of the moduli space. First, if one works with
the entire moduli space of semi-stable Higgs bundles, 
$h_0:\Higgs^0_G\to \cB_\fg$ is proper and even though it has some singular fibres, there are
isomorphisms $R^k h_{0\ast}\cO_{\Higgs}\simeq \Omega_\cB^k$. For $k=0,1$ and $G=SL_2(\CC)$ this is discussed in \cite[\S 6]{hitchin_sb}
and \cite[\S 5]{hitchin_flat_GQ}. For $SL_n$ and higher direct images this is \cite[Theorem 15]{arinkin_compactified} and the general
result is due to  \cite{teleman_frenkel}.
 What is important for us is that $h_{0\ast}\cO_{\Higgs^0}=\cO_{\cB}$ and 
\[
 T^\vee_\cB \simeq R^1 h_{0\ast}\cO_{\Higgs^0}.
\]
Next, since $\cB_\fg$ is a vector space, the Leray spectral sequence implies that for any coherent sheaf $\cF$
on $\Higgs_G^0$ we have  an isomorphism $H^0(\cB,R^k h_{0\ast} \cF)\simeq H^k(\Higgs^0_G,\cF)$ and all higher cohomologies
of $R^k h_{0\ast} \cF$ vanish. For example, $H^0(\cB,T^\vee_\cB)\simeq H^1(\Higgs^0_G,\cO)$.

If we restrict ourselves to $\Higgs^{0,reg}_G$ now, the Hitchin map is not proper anymore. In particular, the direct image sheaves,
such as $h_{0\ast}\cO_{\Higgs^{reg}}$ and $h_\ast h^\ast T_\cB=T_\cB\otimes h_{\ast}\cO_{\Higgs^{reg}}$,
need not be coherent anymore. However, since the moduli space is normal (\cite{moduli2}), by Leray and Hartogs' theorems for
higher cohomologies (\cite{higher_hartogs}) we have an isomorphism
\[
 H^0(\cB,R^k h_{0\ast}\cF)\simeq  H^k(\Higgs^{0,reg}_G,\cF)\simeq H^k(\Higgs^{0}_G,\cF)
\]
for any coherent sheaf $\cF$ on $\Higgs^0_G$. In particular, while the canonical map
$\cO_\cB\to h_{0\ast}\cO_{\Higgs^{reg}}$ is not an isomorphism, it induces an isomorphism on global sections:
\[
 H^0(\cB,\cO_\cB)= H^0(\Higgs^{0,reg}_G,\cO)= H^0(\Higgs^{0}_G,\cO).
\]

Since  $\Higgs^{0,reg}_G$ is holomorphic symplectic and
$h_0: \Higgs_G^{0,reg}\to\cB_\fg$ is a (non-proper) holomorphic submersion with Lagrangian fibres, we have an isomorphism
$h_0^\ast T_\cB^\vee\simeq T_{\Higgs^{reg}/\cB}$ induced by $\omega^{-1}\circ dh_0^\vee$.
We notice that we are working with the full Hitchin base, whose (co)tangent bundle is canonically trivial: $T^\vee_\cB=\cB^\vee\ctimes \cO_\cB$, where 
$\cB_\fg^\vee=\mhom(\cB_\fg,\CC)$ is  isomorphic to $\bigoplus_i H^1(X,K_X^{-m_i})$ by Serre duality.

In the usual fashion, we obtain a locally transitive infinitesimal action of $h_0^\ast T^\vee_\cB$ along the fibres of
$\Higgs^{0,reg}$.  Over $\scB$ (and in fact, over the larger open of non-singular Hitchin fibres) this action integrates
to an action of the fibres $T^\vee_{\scB,b}=\cB^\vee$ (considered as abelian groups), and the Hitchin section $\fv$ determines a holomorphic flow (``exponential'') map
\[
 \tot T^\vee_{\scB}\times_{\scB} \left(\Higgs^{0,reg}\backslash h_0^{-1}(\Delta)\right) \longrightarrow  \Higgs^{0,reg}\backslash h_0^{-1}(\Delta).
\]
The infinitesimal action of $T^\vee_\cB$ is defined  at all smooth points of $\Higgs^0$ and we denote by  $\cU\subset \tot T^\vee_\cB$
 the largest  open 
on which the flow map can be defined. We note that $\cU$ contains $\cB\subset \tot T^\vee_\cB$ (the zero-section) as well as
$\tot T^\vee_\scB$.

To be completely explicit,  $\tot T^\vee_\cB=\cB\times \cB^\vee$, 
\begin{equation}\label{flow1}
 \tot T^\vee_\cB\times_{\cB}\Higgs^{0,reg}= \tot h_0^\ast T^\vee_{\cB}= \Higgs^{0,reg}\times \cB^\vee
\end{equation}
and  the flow map
\begin{equation}\label{flow2}
 \tot T^\vee_\cB\supset \cU\longrightarrow \Higgs_{G,X}^0,\ 
\end{equation}
is given by $(b,\alpha)\longmapsto \exp X_\alpha \cdot \fv(b)$,
where $X_\alpha$ is the Hamiltonian vector field, corresponding to $\alpha\in \cB^\vee$. 
An explicit expression for this map is provided in the third item of the following theorem.

      \begin{thm}
Let $K_X^{1/2}$ be a theta-characteristic on $X$ and
$\{\cy, \cih,\cx\}$  the distinguished principal $\slt$-triple associated with a choice of \'epinglage data.
Let $\fv:\cB_\fg\to \Higgs_{G,X}^0$ and  $(\bdP,\theta)$ be the
  Hitchin section and the uniformising Higgs bundle, associated with these data.
In particular $[(\bdP,\theta)]=\fv(0)$,
$\ad\bdP=\bdE$ is the bundle (\ref{Lie_bundle}) and $\theta=\cy\in H^0(\bdE\otimes K_X)$. Let
$H^0(\bdE\otimes K_X)^{st}$ be the non-empty open set of stable Higgs structures on $\bdP$. Let
\[
 \widetilde{\fii}: H^0(\bdE\otimes K_X)^{st}\longrightarrow \textrm{Gr}\left(\dim \Bun_G, H^1(\bdE)\right)
\]
be the map $\eta\longmapsto \left[\ker h^1(\ad_\eta)\right]$ and let $\cF=\widetilde{\fii}^\ast \cS\subset H^1(\bdE)\ctimes \cO$
be the pullback of the tautological vector bundle $\cS$ on the Grassmannian. Then:
	    \begin{enumerate}
	     \item The total space $\tot \cF\subset\tot T^\vee_{H^0(\bdE\otimes K_X)^{st}}$ carries a tautological analytic family of Higgs bundles
$\left(\scP\to \tot\cF,\Theta\right)$ which
is a deformation of $(\bdP,\theta)$. 
If $\cU\subset \tot\cF$ is the (non-empty)  open set of semi-stable Higgs bundles and
 $\Phi: \cU\to \Higgs^0_{G,X}$ the classifying map, then $\Phi$
 maps fibres of $\cF$ to Hitchin fibres and the zero-section of $\cF$ to $\fv(\cB)$.
	     \item The restriction of $(\scP,\Theta)$ to   $\tot \left(\left. \cF\right|_{\bdk(\cB_\fg)}\right)$ is  a locally universal family of deformations of
$(\bdP,\theta)$. 
 Moreover,  $\Phi^\ast\omega_{\Higgs}=\omega_{can}$, where $\omega_{can}$ is the restriction of the canonical symplectic form on
$\tot T^\vee_{H^0(\bdE\otimes K_X)^{st}}$.
	      \item The choice of principal subalgebra determines a trivialisation
\[
 \tot \left(\left. \cF\right|_{\bdk(\cB_\fg)}\right)\simeq \tot T^\vee_{\cB}\simeq \bdk(\cB)\times \left(\bigoplus_{i=1}^l \fg_{-m_i,i}\otimes H^1(K_X^{-m_i})\right)
\]
\[
 \left(h,v\right) \longmapsto \left(h, \sum_{k=0}^{M}(-1)^k\left(\cP\circ \ad_\cih\right)^k(v)\right),
\]
where $\cP$ is the splitting of $\ad\cy$ induced by $\ad\cx$ via the inclusion
$\mend_1(\fg)\subset \mhom_1\left(H^1(\bdE\otimes K_X), H^1(\bdE)\right)$.
	    \end{enumerate}

      \end{thm}

      \emph{Sketch of proof:}
We outline the idea of the proof here. While a more detailed discussion will be given at another occasion, the main ingredients are
to be found in  \cite[Chapter  7]{thesis} and  \cite[ \S 6]{symplectic_kur}.

We construct the family $(\scE,\Theta)$ by using a small amount of Kodaira--Spencer--Kuranishi theory
and differential graded Lie algebras (dgla). The deformations of the pair $(\bdP,\theta)$ are controlled by the
Biswas--Ramanan  complex (\ref{BR_complex}). The controlling dgla (in the sense of \cite{goldman-millson} and \cite{maurer})
is given by the global sections of the Dolbeault resolution of  (\ref{BR_complex}), see also \cite{simpson_hodge}
and \cite{moduli2}.
 Explicitly, the dgla in question is the vector space $\bigoplus_{p,q} A^p(X,\ad\bdP\otimes\Omega_X^q)$, 
graded by total degree ($p+q$), with bracket induced by combining wedge product and commutators, and having differential $\dbar_{\bdP}+\ad \theta$.
The solutions of the Maurer--Cartan equation are then pairs of elements $(h,v)\in A_X^0(\bdE\otimes K_X)\oplus A_X^{0,1}(\bdE)$,
satisfying $\dbar h+ [\theta+h,v]=0$. 
If moreover  $h\in H^0(\bdE\otimes K_X)\subset A_X^0(\bdE\otimes K_X)$, this reduces to
$v\in \ker h^1(\ad (\theta+h))$ with notation as in (\ref{cone_hyper}).

For $h=0$ one has $\ker h^1(\ad\theta)= \left(\bigoplus_{i=1}^l \fg_{-m_i,i}\otimes H^1(K_X^{-m_i})\right)$,
which is identified with $\cB_\fg^\vee= T^\vee_{\cB_\fg,0}$.
As $h$ varies over $\bdk(\cB_\fg)\subset H^0(\bdE\otimes K_X)$ one has a varying family of
``centralisers''  $\ker h^1(\ad (\theta+h))\subset H^1(\bdE)$.
It turns out that the choice of principal subalgebra determines a trivialisation of this family, that is, an
isomorphism 
 $\tot \left(\left. \cF\right|_{\bdk(\cB_\fg)}\right)\simeq \cB_\fg\times \cB^\vee_\fg$, where the latter is identified with 
an affine subspace of
$H^0(\bdE\otimes K_X)\times H^1(\bdE)$, namely
$\bdk(\cB)\times \left(\bigoplus_{i=1}^l \fg_{-m_i,i}\otimes H^1(K_X^{-m_i})\right)$.
Indeed, we have an isomorphism $\fz_\cx=\ker\ad_\cx\simeq \textrm{coker}(\ad_\cy)$, and hence
$\cP\in\mhom_1 (\textrm{Im }\ad\cy,\fg)$, a splitting of $\theta=\ad\cy$.
 It can be identified with a degree-1 homomorphism from $H^1(\bdE\otimes K_X)$ to $H^1(\bdE)$,  denoted with the same letter. 
One then checks directly that the formula in $(3)$ provides a trivialisation of the bundle of centralisers and that this trivialisation is
symplectic.
The formula in $(3)$ can be obtained as a  ``symplectic version'' of the formal power series solution of the Maurer--Cartan equation
(see e.g. \cite{kodaira_nirenberg_spencer}) in which Green's operator is replaced by the splitting $\cP$.

The local universality follows from Hodge-theory, bijectivity of the Kodaira--Spencer map at $(0)$ and the fact that $(\bdP,\theta)$ is regularly stable.

\qed

We make now several brief comments about applications and related results.
The construction from the Theorem provides Darboux coordinates in a neighbourhood of the
Hitchin section (and not in a neighbourhood of a smooth fibre, as is usual). The explicit description of
the flow map allows one to give some approximation of the image of the brane of opers under the non-abelian
Hodge correspondence, see \cite{thesis}. For more on the relation between opers, non-abelian Hodge
theory and physics see \cite{KW} and the recent preprint \cite{mulase_neitzke_etal}.
As suggested in \cite{KW}, the hamiltonian flow along the Hitchin fibres can be considered as an analytic analogue 
of the so-called ``Hecke operators'' (see \cite{don-pan} for the definition).
Our formula $(3)$ is compatible (but not identical!) with a similar  formula of C.Teleman for 
 $G=GL_n$, see  \cite[\S 7.3]{teleman_langlands}.

	  \section{Special K\"ahler Geometry }\label{sk}
	\subsection{}
In this section we review briefly a differential-geometric structure called \emph{special K\"ahler geometry},
which was first discovered by physicists (\cite{bcov}, \cite{seiberg_witten_1}) in the context of $N=2$ supersymmetry in four dimensions.
This structure exists on the base of any algebraic completely integrable Hamiltonian system (away from the discriminant locus). Conversely, 
such  data give rise to an algebraic integrable system.
The case of interest for us is the Hitchin base $\scB_\fg\subset\cB_\fg$ (or the slightly larger locus of non-singular Hitchin fibres).
Part of these data, the \emph{Donagi--Markman cubic},  is purely holomorphic and can be identified as  the infinitesimal period map of the integrable system.
For the Hitchin system the cubic has been computed by Balduzzi and Pantev.
Together with Ugo Bruzzo we have extended the Balduzzi--Pantev calculation to the case of 
 the \emph{generalised} Hitchin system. This is the topic of Section \ref{dm_cubic}.

I was introduced to this subject during Tony Pantev's lectures in 2003. One of the most detailed and elegant intrinsic introductions to
this material is \cite{danfreed}. Other illuminating references are \cite{bartocci_mencattini_sk}, \cite{ugo_lag_tf} and \cite{markman_sw}. The relation to
$tt^\ast$-geometry is discussed in \cite{hertling} and \cite{hertling_hoev_posthum}.
	\subsection{Intrinsic definition}
We begin with the intrinsic definition of special K\"ahler geometry, leaving the extrinsic (coordinate) definition for the next subsection.
      \begin{defn}[\cite{bcov},\cite{danfreed}]\label{sk_def}
Let $(B,\omega)= (M,I,\omega)$ be a K\"ahler manifold with symplectic form $\omega$,  almost complex structure $I$ and  underlying real manifold $M$.
A \emph{special K\"ahler (SK)}  structure on $(B,\omega)$ is a connection $\nabla$ on $T_M= T_{B,\RR}$, 
which is:
	\begin{enumerate}
	 \item Flat
	 \item Symplectic 
	 \item Torsion-free 
	 \item Special 
	\end{enumerate}
or, with formulae,
	\begin{enumerate}
	 \item  $(d^\nabla)^2=0$
	 \item	$\nabla \omega=0$
	 \item  $d^\nabla(\boldsymbol{1})=0$
	 \item $d^\nabla (I)=0$.
	\end{enumerate}
      \end{defn}

There is a natural notion of morphism of special K\"ahler manifolds: a morphism of K\"ahler manifolds, preserving the connections.

For a connection $\nabla$ on $T_M$, the operator
 $d^\nabla: \cA_\RR^p(T_M)\to \cA_\RR^{p+1}(T_M)$ is an extension of  $\nabla: T_M\to \cA^1(T_M)$ via the exterior differential, i.e.
$d^\nabla(\alpha\otimes v) = d\alpha\otimes v + (-1)^p\alpha\wedge \nabla v$ for a $p$-form $\alpha$ and a vector field $v$.
 Next, we recall that $\nabla$ induces a connection $\nabla^\vee$
on $T^\vee_M=\cA^1_\RR$ via $(\nabla^\vee_X \Phi) \cdot v = X(\Phi\cdot v) - \Phi\cdot (\nabla_X v)$,  and hence a connection 
$\widetilde{\nabla}= \nabla^\vee\otimes 1+ 1\otimes \nabla$ on $T^\vee_M\otimes T_M$.
Usually this connection is denoted simply by  $\nabla$, as in condition
 $(2)$, but for now we   keep the (somewhat pedantic) notation $\widetilde{\nabla}$.
In this section we also denote by $\cdot$ the canonical pairing between $1$-forms and vector fields.

Since we have a sheaf isomorphism
$\send (T_M)\simeq \cA^1_\RR(T_M)$, we can compare $d^\nabla$ to (the anti-symmetrisation of) $\widetilde{\nabla}$. It is not hard to check
that the obvious diagram
      \[
    \xymatrix{\send (T_M)\ar[r]^-{\widetilde{\nabla}} \ar[d]_-{\simeq} & \cA^1_{\RR}(\send T_M)\simeq (\cA_\RR^1)^{\otimes 2}\otimes T_M\ar[d]^-{Alt}\\
	     \cA^1_\RR(T_M)\ar[r]^-{d^\nabla} & \cA^2_\RR(T_M)  }
      \]
 \emph{does not} commute, and the failure is  the torsion of $\nabla$:
 $Alt\circ \widetilde{\nabla}= d^\nabla - d^\nabla (\boldsymbol{1})$. Explicitly, this means
that 
\[
(\widetilde{\nabla}_X \Phi) \cdot Y - (\widetilde{\nabla}_Y \Phi) \cdot X =(d^\nabla \Phi)(X,Y)- \left(\nabla_X Y- \nabla_Y X-[X,Y]\right).
\]
Next, the complexification $I_\CC$ of $I$ decomposes $T_{B,\CC}$ into $\pm i$ eigenbundles and the 
$(1,0)$  projector
\[
  T_{B,\CC}\simeq T_B^{1,0}\oplus T_B^{0,1}\longrightarrow T_B^{1,0}
\]
is  precisely $\pi^{1,0}=\frac{1}{2}(\boldsymbol{1}-iI_\CC)$. 
 Taking real and imaginary parts, we can rephrase
  conditions $(3)$ and $(4)$ as a $d^\nabla$-horizontality condition for $\pi^{1,0}$:
\[
 d^\nabla(\pi^{1,0}) =0 \Longleftrightarrow   
		\left|
	    \begin{array}{l}
	     d^\nabla(\boldsymbol{1})=0\\
	     d^\nabla(I)=0\\
	    \end{array} 
		\right. 
\Longleftrightarrow
		\left|
	    \begin{array}{l}
	     d^\nabla(\boldsymbol{1})=0\\
	     \nabla(I)\in \Gamma(\sym^2 T^\vee_M\otimes T_M)\\
	    \end{array} 
		\right. .
\]
For simplicity we do not distinguish notationally between $\nabla$ (a connection on $T_{M}=T_{B,\RR}$) and its
complexification $\nabla^\CC$ on $T_{B,\CC}$.

We emphasise that
conditions $(3)$, $(4)$ do not imply the vanishing of  $\nabla (I)$ but just put a symmetry restriction on it.
In particular,  $\nabla$ need not be the Levi-Civita connection for the K\"ahler metric.
Of course, being flat K\"ahler implies being special K\"ahler, but not conversely. 

We also emphasise that we have not imposed any compactness restrictions on $M$. In fact, as shown by Lu
(\cite{cpct_sk}),
 the only  \emph{compact}
special K\"ahler manifolds are the compact flat  K\"ahler ones.
	  \subsection{Coordinate description}
Conditions (1), (2) and (3) imply that $M$ admits   flat local Darboux coordinates. In other words, near every point 
$p\in M$ there is an open $U\subset M$, with local coordinates
 $\{x_i, y_i\}$, such that 
$\omega_U =\sum dx_i\wedge dy_i$ and  $\nabla(dx_i)=0$, $\nabla (dy_i)=0$. Indeed, since
$\nabla$ is flat symplectic, $\omega_U$  can be written in terms of flat local  sections of $T^\vee_M$. But
for a torsion-free connection a horizontal 1-form is  closed (and hence locally exact).
Any two choices of such coordinates differ by an affine-linear transformation, whose linear part is symplectic.
I.e., if $(B,\omega,\nabla)$ satisfies $(1)$, $(2)$, $(3)$ then $M$ is equipped with a flat symplectic structure
and in particular, it is endowed with an affine structure. A choice of basepoint $o\in M$ and a trivialisation
 $T_{M,o}\simeq \RR^{2\dim B}$ determine a monodromy representation $\pi_1(B,o)\to \textrm{Sp}(2\dim B,\RR)$.

Since $(d^\nabla)^2=0$, the  condition $d^\nabla\pi^{1,0}=0$ implies that  locally
$\pi^{1,0}_U=\nabla\zeta$, for some complex vector field $\zeta\in \Gamma(U,T_{B,\CC})$.
Such  a $\zeta$ is unique up to a $\nabla$-flat vector field.
 Then
\[
 \zeta=\frac{1}{2}\left(\sum_i z_i\frac{\partial}{\partial x_i}- w_i\frac{\partial}{\partial y_i}\right),
\]
for uniquely determined functions $z_i$, $w_i\in C^\infty(U)$. By the flatness of the Darboux frame
\[
 \nabla \zeta = \frac{1}{2}\left(\sum_i dz_i\otimes \frac{\partial}{\partial x_i}- dw_i\otimes \frac{\partial}{\partial y_i}\right).
\]
Finally, since $\pi^{1,0}\in A^{1,0}(T^{1,0})$, we have  $\dbar z_i=0=\dbar w_i$, so $z_i,w_i\in \cO_B(U)$. Since
$\textrm{Re} \pi^{1,0}=\frac{1}{2}\boldsymbol{1}$, we have
\[
 \textrm{Re }(dz_i)=dx_i, \ \textrm{Re }(dw_i)= -dy_i.
\]
Consequently,   $z_i$ are complex coordinates on $U\subset M$ (and so are $w_i$)
and we get
\begin{equation}\label{projector}
 \pi^{1,0}=\nabla \zeta = \sum_i dz_i\otimes \frac{\partial}{\partial z_i}= \frac{1}{2}\left(\sum_i dz_i\otimes \frac{\partial}{\partial x_i}- dw_i\otimes \frac{\partial}{\partial y_i}\right).
\end{equation}

The coordinates $\{z_i\}$  are called \emph{special coordinates}, adapted to the flat Darboux coordinates $\{x_i, y_i\}$ 
and $\{\omega_i\}$ are the \emph{dual (conjugate) special coordinates}, see \cite{danfreed}. 

Since we have two sets of local coordinates, we can consider the matrix of functions $\tau = (\tau_{ij})\in \textrm{Mat}_{\dim B}(\cO_B(U))$ relating the
respective  coframes:
\begin{equation}\label{def_tau}
 dw_i =\sum_j \tau_{ji}dz_j.
\end{equation}
  Consequently, by (\ref{projector}) we have
\begin{equation}\label{special_coords}
 \frac{\partial}{\partial z_i}= \frac{1}{2}\left(\frac{\partial}{\partial x_i}- \sum_j \tau_{ij}\frac{\partial}{\partial y_j} \right).
\end{equation}
The K\"ahler condition on $(B,\omega)$ imposes significant restrictions on $\tau$.
As    (the complexification of) $\omega$ is of type $(1,1)$,
substituting (\ref{special_coords})  in $\omega \left(\frac{\partial}{\partial z_i}, \frac{\partial}{\partial z_j} \right)=0$
we obtain  that $\tau= \tau^t$.
More conceptually, the symmetry of $\tau $ follows from  the equality
\begin{equation}\label{20_closed}
  0=  \omega^{2,0}= d \left(\sum_i w_i dz_i \right)= \sum_{i,j}\tau_{ij} dz_i\wedge dz_j.
\end{equation}
But then $\sum_i w_i dz_i = d\cF$ for some holomorphic function  $\cF\in \cO_B(U)$, possibly after  shrinking $U$.
  Consequently,  $\tau$ must be a Hessian of $\cF$:
\begin{equation}\label{tau_hessian}
 \tau_{ji}= \frac{\partial w_i}{\partial z_j}= \frac{\partial^2\cF}{\partial z_i\partial z_j},\ \tau_{ij}=\tau_{ji}.
\end{equation}
 Consequently, the K\"ahler form is
\begin{equation}\label{sk_form}
 \omega = \sum_i dx_i\wedge dy_i = \frac{i}{2}\sum_{j,k}Im(\tau)_{jk}dz_j\wedge d\overline{z_k}
\end{equation}
and comes from a K\"ahler potential $\frac{1}{2}\textrm{Im}\left(\sum_k w_k\overline{z}_k\right)$.
Finally, the condition $\omega>0$ implies $\textrm{Im}(\tau)>0$. 
In this way we obtain a holomorphic map $\tau: U\to \HH^{\dim B}$, where
\[
\HH^{\dim B} =\left\{ Z\in \textrm{Mat}_{\dim B}\left| Z^t=Z, \textrm{Im}Z>0\right.\right\}
\]
is Siegel's upper half space. This map is of a very special form: it arises as a Hessian of a holomorphic function.

The special K\"ahler metric determines  $\cF$  only up to affine-linear
terms. Conversely, any choice of such an $\cF$ determines   locally the special K\"ahler structure.
From equation (\ref{20_closed}), i.e., $\sum_{i} w_i dz_i =d\cF$  we see that after modifying by $\cF$ an affine-linear term we have
\begin{equation}\label{dual_coords}
 w_i=\frac{\partial \cF}{\partial z_i}, i=1\ldots \dim B.
\end{equation}
The holomorphicity of $\tau$ (and $\cF$) also gives a convenient description of the special K\"ahler connection. Namely,
equation (\ref{special_coords}) and the flatness of $\{x_i, y_i\}$ imply that
\begin{equation}\label{sk_conn_cubic}
 \nabla\left(\frac{\partial}{\partial z_i}\right) = -\frac{1}{2} \sum_j \partial \tau_{ij}\otimes\frac{\partial}{\partial y_j} =
-\frac{1}{2} \sum_{i,j,k}\frac{\partial^3\cF}{\partial z_k \partial z_j\partial z_i} dz_k\otimes \frac{\partial}{\partial y_j}.
\end{equation}
In particular, $\nabla^{0,1}=\dbar_{T_B}$.

In the physics literature the function $\cF$ is called ``holomorphic pre-potential''.
      \subsection{Related Geometries}
In this subsection
we review  some equivalent ways to repackage the special K\"ahler geometry and review its  relations to other mathematical structures.
      \subsubsection{Weight-1 $\RR$VHS}
Recall  that a weight-one real variation of Hodge structures ($\RR$VHS) on $B$ is given by a quadruple
\[
 \left(\cF^1\subset \cF^0, \cF_\RR\subset \cF^0, \nabla^{GM}, Q\right) 
\]
consisting of a length-one flag of holomorphic bundles $\cF^\bullet$, a real subbundle $\cF_\RR$ of $\cF^0$,
a holomorphic flat connection $\nabla^{GM}$ (Gauss--Manin connection) and a polarisation $Q$. These data have to satisfy certain compatibility
conditions (\cite[Ch.III]{voisin1}).

In order to avoid confusion, let us recall our  notation for the different  tangent bundles that we use. First, 
$T_{B,\CC}\simeq T^{1,0}\oplus T^{0,1}$ is the complexified tangent bundle of $B=(M,I)$, $T_B\subset T^{1,0}= T_B\otimes_{\cO_B} \scC_B^\infty$
is  the holomorphic (respectively,  $(1,0)$-) tangent bundle and
$T_M= T_{B,\RR}\subset T_{B,\CC}$ is the real tangent bundle.

As discovered in \cite{hertling} (see also \cite{bartocci_mencattini_sk}),
the data $(M,I,\omega, \nabla)$ is equivalent to the data  of a certain weight one, polarised $\RR$VHS.
Let 
$\nabla = \nabla^{1,0}+\nabla^{0,1}$ be the type decomposition of the (complexification of the) special K\"ahler connection.
Then $\nabla^2=0$ implies that $\nabla^{0,1}$ is a holomorphic structure on the vector bundle
$T_{B,\CC}$ and one takes as $\cF^0$ precisely that bundle, $(T_{B,\CC},\nabla^{0,1})=\ker\nabla^{0,1}$.
The flatness of $\nabla$ also implies that $\nabla^{1,0}$ is a flat holomorphic connection on $\cF^0$.
Finally, the positive-definitness  of $\omega$ implies that it can be used as a polarisation.
Altogether, 
the $\RR$VHS associated to the special K\"ahler  data is the quadruple
\[
\left(T_B\subset (T_{B,\CC},\nabla^{0,1}), T_M, \nabla^{1,0}, \omega   \right).
\]
Notice that the  polarisation $Q$ determines an isomorphism $\cF^0/\cF^1=\overline{\cF^1}\simeq_Q \cF^{1\vee}$
and so $\cF^0$ fits in  an extension
\begin{equation}\label{sk_ext}
 \xymatrix@1{0\ar[r]& T_B\ar[r]&\cF^0\ar[r]& T^\vee_B\ar[r]&0}.
\end{equation}
      \subsubsection{Integral special K\"ahler geometry}
Let us consider a weight-1 $\RR$VHS which is induced by a weight-1 $\ZZ$VHS,
 i.e., $\cF_\RR= \cF_\ZZ\otimes \scC_B^\infty$, for
some locally constant sheaf of lattices $\cF_\ZZ\subset T_M$. 
Then  the dual lattice $\widehat{\cF}_\ZZ\subset T^\vee_B$
 gives rise to a holomorphic family of complex tori 
\[
 h: \cH=\cF^1/\widehat{\cF}_\ZZ= T^\vee_B/\widehat{\cF}_\ZZ \longrightarrow B,
\]
which  is in fact a family of (polarised) abelian varieties.

If $\widehat{\cF}_\ZZ\subset T^\vee_B$ is
Lagrangian  the canonical symplectic form $\omega_{can}$ on $\tot T^\vee_B$  descends to $\cH$ and the fibres
$h^{-1}(b)\subset \cH$ are Lagrangian.
In general, $\widehat{\cF}_\ZZ$ is \emph{not} Lagrangian and none of the above holds.

      \begin{defn}
 An \emph{integral} special K\"ahler structure is a special K\"ahler structure $(B=(M,I),\omega,\nabla)$ together with a
$\nabla$-parallel  sheaf of lattices $\cF_\ZZ\subset T_B$, such that $\widehat{\cF}_\ZZ\subset T^\vee_B$ is a complex-Lagrangian submanifold.
      \end{defn}

An important result, due to Donagi and Witten (\cite{donagi_witten}, \cite{donagi_swis}, \cite[Theorem 3.4]{danfreed}, \cite[Theorem 2.1]{markman_sw})
establishes an equivalence between the data of an integral special K\"ahler structure on $B$ and the data of a (polarised) algebraic completely integrable hamiltonian system
(ACIHS) $h:\cH\to B$ whose fibres are polarised abelian varieties. Above we have indicated how to construct an ACIHS from the integral special K\"ahler
data, so here we outline the inverse construction.

Let $\omega_\cH$ be the holomorphic symplectic form on $\cH$.
The assumption that the fibres of $h$ are abelian varieties implies that $h$ admits a section, $\fv:B\to \cH$. We may assume that $\fv$ is Lagrangian
(any section  becomes Lagrangian after replacing $\omega_\cH$ with $\omega_\cH - (\fv\circ h)^\ast\omega_\cH$).
Since $h$ is proper, $\omega_\cH$ induces, as in Section (\ref{flow_flow}), an isomorphism 
\[
T^\vee_B\simeq h_\ast T_{\cH/B}= h_\ast \Omega^{1\vee}_{\cH/B}
\]
induced by $\omega_\cH^{-1}\circ dh^\vee$.
Consider then $R^1\pi_\ast \ZZ$, the local system of first integral cohomologies of the fibres of $h$. Its dual,
the local system of homologies, admits a natural embedding
$\mhom(R^1\pi_\ast \ZZ,\ZZ)\hookr h_\ast \Omega^{1\vee}_{\cH/B}$ given by ``relative integration''. Pointwise
this is the integration homomorphism $H_1(\cH_b,\ZZ)\hookr H^0(\cH_b,\Omega^1)^\vee$. We then denote by $\widehat{\cF}_\ZZ$
the corresponding lattice in the cotangent bundle:
\[
 \widehat{\cF}_\ZZ = \textrm{Im}\left(\mhom(R^1\pi_\ast \ZZ,\ZZ)\hookr h_\ast \Omega^{1\vee}_{\cH/B} \simeq T^\vee_B\right).
\]
Its dual lattice $\cF_\ZZ\subset T_B$ is isomorphic to $R^1 h_\ast\ZZ$. We now have an isomorphism of varieties
$T_B^\vee/\widehat{\cF}_\ZZ\simeq \cH$ induced by $\fv$, as in Section (\ref{flow_flow}). This is in fact
a (local) symplectomorphism (see \cite[\S 44]{guillemin_sternberg} ) and hence $\widehat{\cF}_\ZZ$ is Lagrangian,
  being the preimage in $T^\vee_B$ of the Lagrangian section $\fv$.
 Let $n:=\dim B$, $U\subset B$ a (contractible) open and $\{\gamma_1,\ldots,\gamma_{2n}\}$ a symplectic basis for the
polarisation (say, $\rho$)
over $U$. This means that for $1\leq i,j \leq n$ we have $\rho(\gamma_i,\gamma_j)=0$, $\rho(\gamma_i,\gamma_{n+j})=\delta_i \delta_{ij}$, 
where $\delta_1$, \ldots, $\delta_n\in \NN$ are the ``divisors of polarisation'' of the abelian varieties $\cH_b$, $b\in B$.
Since $\widehat{\cF}_\ZZ\subset T^\vee_B$ is lagrangian, the holomorphic 1-forms on $U\subset B$ (corresponding to) $\{\gamma_i\}$ are closed and hence (locally) exact.
Thus there exist holomorphic functions $\{u_i, u_i^D\}$ on $U$, such that
\[
 \{\gamma_1,\ldots, \gamma_n, \gamma_{n+1},\ldots, \gamma_{2n}\}= \{du_1, \ldots, du_n, du^D_1,\ldots, du^D_n\}.
\]
The corresponding family of period matrices is $\Pi=\left(\Delta_\delta\left| \tau\right.\right)$, where
$\tau: U\to \HH^n$ describes locally the period map of $\cH\to B$ and
  $\Delta_\delta=\textrm{diag}(\delta_1,\ldots \delta_n)$.
 Correspondingly (\cite[Ch.2, \S 6]{GH}) we have 
the relation $du_i^D=\sum_{j}\tau_{ij}\delta_j^{-1}du_j $.
But then the 1-form $\sum_i \delta_i^{-1} u_i^D du_i$ is closed and hence (after possibly shrinking $U$) equal to $d\cF$, for some $\cF\in \cO_B(U)$.
Choosing such an $\cF$ and comparing with equations (\ref{def_tau}), (\ref{20_closed}), we get a pair of dual special coordinates on $U$, namely
$z_i:= u_i/\delta_i$, $w_i:= u_i^D$.

Of course, the above gives only a local description of the special K\"ahler data.
Globally, the  condition that an ACIHS has a  section is, of course, very restrictive. On the other hand with a given ACIHS $h:\cH\to B$
we can associate natural smooth families of abelian varieties:
the relative Albanese fibration $\textrm{Alb}_{\cH/B}$ and the relative Jacobian fibration $\textrm{Jac}_{\cH/B}$.

Let us  emphasise once more  that the ACIHS appearing here have smooth and connected fibres. In more general situations $\cH$ is only Poisson, so
one must consider individual symplectic leaves. On such a leaf  the preimage of the discriminant locus must be removed before applying 
the Donagi--Witten correspondence.

      \subsubsection{The Donagi--Markman cubic}
In the preceding discussion, we have assigned to an ACIHS $h:\cH\to B$ (with smooth connected fibres),  a contractible  open set $U\subset B$
and a trivialisation of $\mhom(R^1 h_{\ast}\ZZ,\ZZ)$ the data of 
 a one-form $dF\in \Gamma(T^\vee_U)$ and $\textrm{Hess}\cF\in \Gamma(\sym^2 T^\vee_U)$,
where the pre-potential $\cF\in \cO(B)$ is defined only up to an affine-linear transformation.
It turns out (\cite{donagi_markman_cubic}) that the third derivative of $\cF$, i.e.,
\[
 c_U= d (\textrm{Hess}\cF) 
= \sum_{i,j,k}\frac{\partial^3 \cF}{\partial z_i\partial z_j\partial z_k}
dz_i\otimes dz_j\otimes dz_k.
\]
is actually the restriction (to $U$) of a global section $c\in H^0(B, \sym^3 \Omega_B^1)$, usually called ``the Donagi--Markman cubic''.
Denoting by  $\Phi: B\to \HH^{\dim B}/\Gamma_\delta$  the classifying map for $h:\cH\to B$, 
can identify the cubic as $c=d\Phi \in H^0(\sym^3\Omega^1_B)\subset H^0(\Omega_B^1\otimes \sym^2\Omega_B^1)$.

The cubic can be arrived at via another route.
Namely, given a family $h: \cH\to B$ of abelian varieties
(or complex tori), satisfying $\dim \cH=2\dim B$, one may ask whether  there exists a holomorphic symplectic  structure 
$\omega_\cH$ on $\cH$, 
for which the fibres of $h$ are  Lagrangian. 
This can be translated to the familiar local picture. Namely,  we have an analytic open $\cU\subset \CC^n$, a holomorphic map $\tau:\cU\to \HH^n$ and a set of polarisation
divisors $(\delta_1,\ldots, \delta_n)$. 
These determine a  group $\Gamma$ of holomorphic automorphisms of 
$T^\vee\cU=\cU\times\CC^n$, generated by
\[
 (z,v)\longmapsto (z, v+ \Pi(z)(e_j)),\ j=1\ldots 2n,
\]
where $\Pi = \left(\left. \Delta_\delta\right|\tau \right):\cU\to \mhom(\CC^{2n},\CC^n)$.
Consider the $2n$ sections $s_j\in\Gamma(\cU, T^\vee\cU)$, defined by the columns of the period matrix $\Pi$, that is,
$s_j(z)= (z, \Pi(z)(e_j))$.
The symplectic structure on $T^\vee\cU$ descends to $T^\vee\cU/\Gamma\simeq_U \cH_U$ and the torus fibres are lagrangian precisely when the
images of the sections $s_j$ are Lagrangian subvarieties of $\cU\times \CC^n$. This happens  precisely when
$d\tau(e_j)=0$, i.e., $\partial_i \tau_{kj}= \partial_k\tau_{ij}$, for all $j$. This, in turn, is equivalent to
the existence of a holomorphic function $\cF: \cU'\to \CC$ on a  possibly smaller subset $\cU'\subset \cU$, such that
$\tau = \textrm{Hess}(\cF)$.
In terms of $h:\cH\to B$ this means that the infinitesimal period map is a section of $\sym^3 \Omega^1_B$.

Since the special K\"ahler data can be repackaged in terms of VHS, one can look for a more direct description of $c$ in Hodge-theoretic
terms. Such a description involves infinitesimal variation of Hodge structures and Higgs bundles arising from system of Hodge bundles,
see \cite{simpson_uniformisation}.
More concretely, to any polarised,  weight-1 $\RR$VHS $(\cF^\bullet, \cF_\RR, \nabla^{GM}, Q)$ on $B$ we can associate 
a $\Omega^1_B$-valued Higgs (vector) bundle $(E,\theta)$. 
This is the vector bundle
$E= \textrm{gr}\cF^\bullet = \cF^1\oplus \cF^{1\vee}$, equipped with the Higgs field
\[
 \theta = \textrm{gr}\nabla^{GM}:\cF^1\to \cF^{1\vee}\otimes \Omega^1_B.
\]
In our case (\ref{sk_ext}) this Higgs  pair is $E=T_B\oplus \Omega^1_B$,  $\theta= \begin{pmatrix}
                                                                  0 & 0\\
								  c & 0\\
                                                                 \end{pmatrix}$,
for $c\in H^0(B,\Omega_B^{1\otimes 3})$.

The relation between systems  of Hodge bundles and Yang--Mills theory was the starting point
of Simpson's study of non-abelian Hodge theory, see e.g., \cite{simpson_uniformisation}. 
In the context of special K\"ahler  geometry this relation is based on the observation (\cite{bcov}, \cite{danfreed})  that
\[
 \nabla = \nabla^{LC}+\theta +\theta^\ast.
\]
That is, the Hitchin--Simpson connection for the above Higgs pair is precisely the special K\"ahler connection.
Moreover,  
Hermite--Yang--Mills equation is the $tt^\ast$-equation (\cite[ 1.32]{danfreed}, \cite{hertling}).
Finally, the Donagi--Markman cubic $c$ is given  explicitly in terms of the special K\"ahler data
(\cite{danfreed}) as
\[
 c=-4\omega(\pi^{1,0},\nabla\pi^{1,0})\in \Gamma (T^\vee_{B,\CC}\otimes \sym^2 T^\vee_{B,\CC}).
\]
Indeed, in  special coordinates we have, from (\ref{sk_conn_cubic}), that
\[
 -\omega\left(\sum_i dz_i\otimes \frac{\partial}{\partial z_i},\sum_j \nabla \left(dz_j\otimes \frac{\partial}{\partial z_j}\right)  \right) =
\frac{1}{4}\sum_{i,j,k}\frac{\partial^3 \cF}{\partial z_i\partial z_j\partial z_k}
dz_i\otimes dz_j\otimes dz_k.
\]
      \subsection{Relations to physics}
Any thorough discussion of the appearances of  special K\"ahler geometry and physics is beyond the scope of the
current lectures. 
 The structures that we have discussed  have entered the mainstream
physics literature around 1984 
from two directions simultaneously: supersymmetry and supergravity.
Probably  the most influential examples have been
 \cite{bcov} and \cite{seiberg_witten_1}, \cite{seiberg_witten_2}.

Let $G_c\subset G$ 
 be a compact real form of the simple complex group $G$.
Seiberg and Witten considered (pure)
 $N=2$ supersymmetric $G_c$ Yang--Mills theory  in four dimensions.
The vacuum of this theory is infinitely degenerate, with $\ft/W$ being the moduli space of vacua.
Seiberg and Witten discovered the presence of special K\"ahler geometry on the complement of the discriminant locus
in $\ft/W$.
The special coordinates $\{z_i\}$ (respectively $\{w_i\}$) describe the electric (respectively, magnetic)
charges of the theory.
As shown in \cite{seiberg_witten_1}, the low-energy effective Lagrangian of the theory can
be expressed in terms of a single function of the electric charge, the prepotential $\cF=\cF(\{z_i\})$.
Supersymmetry implies that all functions involved, in particular, $\cF$ and $\{z_i\}$, are holomorphic.
The  matrix $\tau$ plays the r\^ole of complexified gauge coupling of the theory.
The $S$-duality transformation acts on $\{z_i,w_i\}$ by a finite index subgroup of
$\textrm{Sp}(2l,\ZZ)$.
If adjoint matter is added to theory, the duality acts by affine-linear symplectic transformations,
whose translational part is determined by the added masses.

Seiberg and Witten identified the prepotential and electric charges of the theory for $G=SL_2(\CC)$ by
studying the global properties of the algebraic integrable system that arises in this way.
Later Donagi and Witten (\cite{donagi_witten})  studied the ACIHS which can arise in this way.
In particular, they  proposed that for $G=SL_n(\CC)$
the Seiberg--Witten integrable system can be realised as
 a generalised Hitchin system
over an elliptic curve $X$.
Donagi and Witten considered also  different limits of the theory, such as  keeping the mass fixed and letting the elliptic curve degenerate
or keeping the elliptic curve and taking the limit of zero mass. In the former case one obtains pure $N=2$ theory and in the latter
an $N=4$ theory. 
For more details we direct the reader to
the survey \cite{donagi_swis}, which describes the work in \cite{donagi_witten} with an emphasis 
on the mathematical development.

We also remark that the relation between twisted $N=4$ super Yang--Mills theory in four dimensions and the
(generalised) Hitchin system is at the base of \cite{KW}.

	\section{The Donagi--Markman cubic for the Hitchin system}\label{dm_cubic}
In the study of Hamiltonian systems there is a great difference between proving complete integrability and
linearising (realising) the flow.
In many cases  the period map cannot be determined explicitly, so its derivative, the cubic $c$, is the
next best thing.
It is only natural, then, to try to
compute $c$ for the Hitchin integrable system and its various generalisations.
For the Hitchin system per se this was done by T.Pantev (for $SL_n$, unpublished, but see \cite{ddp} for $SL_2$),
and by D.Balduzzi (\cite{balduzzi}) for arbitrary reductive $G$. For \emph{meromorphic} Higgs bundles this was done by
U.Bruzzo and the author (\cite{ugo_peter_cubic}).
We shall recall the statement of the main theorem and sketch the key steps of the proof.
For more details one can refer to \cite{ugo_peter_cubic} or \cite{survey_mero_higgs}.

We fix an effective divisor $D$ on the Riemann surfaces $X$, and assume that 
$K_X(D)^2$ is very ample. We also set $L=K_X(D)$. 

The coarse moduli space $\Higgs_{G,D}$ contains  a connected component $\Higgs_{G,D}^c$ for each topological type $c\in \pi_1(G)$
(although it is not known whether these are the only connected components of the moduli space).
We write $\cB_{G,D}$ for the generalised Hitchin base $H^0(X,\ft\ctimes L/W)\simeq \oplus H^0(X,L^{m_i+1})$
and consider the Hitchin map(s) $h_c: \Higgs_{G,D}^c\to \cB_{G,D}$.
Markman (\cite{markman_thesis}, \cite{markman_sw}) and Bottacin (\cite{bottacin}) have shown that, whenever
$\Higgs_{G,D}^c$ is not empty, $h_c$ endows it with the structure of a holomorphic ACIHS in the Poisson sense.
Our goal then is to compute $c$ for the (generic) symplectic leaves of $\left. \Higgs^c_{G,D}\right|_{\cB\backslash\Delta}\to \cB_{G,D}\backslash \Delta$,
where $\Delta\in\cB_{G,D}$ is the discriminant locus.
Consider the vector subspace  $\cB_0=\oplus_i H^0(L^{d_i}(-D))\subset \cB_{G,D}$. 
Markman (\cite{markman_thesis}, \cite{markman_sw}) established a bijection between
$\cB_{G,D}/\cB_0$ and the set of closures of symplectic leaves.
The bijection  assigns to a $\cB_0$-coset the closure of the unique symplectic leaf  of maximal rank, contained in the fibre of $\overline{h_c}$,
where $\overline{h_c}$ is the composition of $h_c$ with the quotient $\cB_{G,D}\to \cB_{G,D}/\cB_0$.

Let $\scB\subset \cB_{G,D}$ be the set of  cameral covers which are generic (i.e., nonsingular and with simple ramification). This set is open under the assumptions on 
$D$. We choose a point $o\in\scB$, and consider $\bB=\left(\{o\}+\cB_0\right)\cap\scB$. 
We also denote by $\pi_o:\widetilde{X}_o\to X$ the cameral cover, corresponding
to $o$.
Now  $\left. S\right|_\bB= h_c^{-1}(\bB)\to \bB$ is an integrable system in the
symplectic sense and our theorem is a statement about its infinitesimal period map.
It turns out that $c$, a section of $\sym^3 T^\vee_\bB$ can be computed in terms of cameral data, and
more precisely, as  a quadratic residue of a ``logarithmic derivative'' of the discriminant of $\fg$.
We recall that the discriminant $\fD\in \sym^{|\cR|}(\ft^\vee)$ gives rise to a section
of the line bundle $p_X^\ast L^{|\cR|}$ over $\cB_{G,D}\times X$, denoted by the same letter.

      \begin{thm}[\cite{ugo_peter_cubic}, Theorem A]\label{thm_cubic}
       There is a natural isomorphism 
\[
  T_{\bB,o}\simeq H^0\left(\widetilde{X}_o, \ft\ctimes K_{\widetilde{X}_o}\right)^W, 
\]
denoted by $Y_\xi\mapsto \xi$.   Under this isomorphism, $c_o: T_{\bB,o}\to \sym^2\left(T_{\bB,o}\right)^\vee$
is identified with
\[
 c_o : H^0(\widetilde{X}_o,\ft\ctimes K_{\widetilde{X}_o})^W\longrightarrow \sym^2 \left(H^0(\widetilde{X}_o,\ft\ctimes K_{\widetilde{X}_o})^{W \vee} \right),
\]
\[
 c_o(\xi)(\eta,\zeta) = \frac{1}{2}\sum_{p\in\textrm{Ram }\pi_o}\textrm{Res}^2_p \left( \pi_o^\ast\left(\left.\frac{\cL_{Y_\xi}\fD}{\fD}\right|_{\{o\}\times X}\right) \eta\cup\zeta\right).
\]
      \end{thm}
\emph{Sketch of Proof:}
Let 
$N$ denote the normal bundle of 
$\widetilde{X}_o\subset \tot \ft\ctimes L$ and let $r$
 denote the bundle projection $\tot\ft\ctimes L\to X$. 
The total space of $\ft\ctimes L$ carries
a canonical $\ft$-valued 2-form  $\omega_\ft\in H^0(\tot \ft\otimes L, \Omega^2(r^\ast D))$, generalising
the Liouville symplectic form on $\tot K_X$.
By restriction  we get a map
$\omega_{\ft}:N\to \tot \ft\otimes K(r^\ast D)$. It  induces an isomorphism
$H^0(\widetilde{X}_o, N(-r^\ast D))^W\simeq_{\omega_{\ft}} H^0(\widetilde{X}_o, \ft\otimes K_{\widetilde{X}_o})^W$,
since by \cite{kjiri}, the generalised Hitchin system satisfies the ``rank-2 condition'' of Hurtubise and Markman
(\cite{hurtubise_markman_rk2}). On the other hand,   $H^0(\widetilde{X}_o, N(-r^\ast D))^W$ 
is canonically isomorphic to $T_{\bB,o}$.

Over sufficiently small (analytic) opens $\cU\subset \bB$ the fibration $\Prym^o_{\cX/\cU}\to \cU$ admits sections and
can be identified with $\left. S\right|_{\cU}\to \cU$.
One then shows by a fairly standard argument  that the symplectic structure on $\left. S\right|_{\cU}$ 
can be identified with the canonical symplectic structure on $\Prym^o_{\cX/\cU}$.
 We have that $T_{P_o}\simeq H^1(\widetilde{X}_o, \ft\otimes \cO)^W\otimes \cO_{P_o}$, and the canonical symplectic structure
is built by splitting the tangent space to  $\Prym_{\cX/\cU}$ into self-dual spaces.
 Since the complex structure of
the Prym is induced by the complex structure of the cameral curve $\widetilde{X}_o$ and 
$H^1(T_{P_o})=H^1(\widetilde{X}_o,\ft\otimes \cO)^{W\otimes 2}$.
In this way we have reduced the question of computing the infinitesimal period map of 
$\left. S\right|_{\cU}\to \cU$ to the question of computing it for $\cX\to \cU$.
However, by a theorem of Griffiths (\cite{griffiths_periods}) the  infinitesimal period map can be obtained from the Kodaira--Spencer map
$\kappa : T_{\bB,o}\to H^1(\widetilde{X}_o, T_{\widetilde{X}_o})$ via
\[
 c_o(\xi)(\eta,\zeta) = \frac{1}{2\pi i}\int_{\widetilde{X}_o}\kappa(Y_\xi)\cup \eta\cup \zeta.
\]
Finally, since $\cX\subset \tot\ft\ctimes L$ is a complete intersection, cut out by the invariant polynomials $I_k$,
$\kappa(Y_\xi)$ can be computed on an appropriate open cover, using the genericity assumption.

\qed

	    \section{The generalised $G_2$ Hitchin system and Langlands duality}\label{G2}
In this section we shall discuss some very basic properties of the generalised Hitchin system for the 
exceptional group $G_2$.
The (usual) $G_2$ Hitchin system has been extensively studied in \cite{katz_pan_G2} and \cite{hitchin_G2},
while some aspects   of the generalised (ramified) setup have been discussed, for  $X=\PP^1$, in
\cite{adlervm_cis}.
For general properties of the generalised Hitchin systems we refer the reader to the references given
in the surveys \cite{donagi_markman}, \cite{markman_sw} and \cite{survey_mero_higgs}.

The exposition below fits within the context of an ongoing joint project with U.Bruzzo. It
 is largely motivated by trying to grasp the circle of ideas discussed in
 \cite{don-pan},  \cite{hitchin_G2} and \cite{argyres_kapustin_seiberg}, and their implications
for the ramified Hitchin system and its special K\"ahler geometry.
Below we state a result concerning the invariance of the Donagi--Markman cubic under the Langlands involution of the Hitchin base.
It is a direct extension to the ramified case of a result of Hitchin in \cite{hitchin_G2} and its complete proof will be discussed
in a forthcoming work.

We start with  an important Lie-theoretic observation (\cite[Remark 3.1]{don-pan}).
If $\fg$ is a simple Lie algebra of type ${\sf B}$ or ${\sf C}$, any choice of 
 Killing form $\langle, \rangle$ determines an isomorphism $\ft\simeq ^L\ft$ by composing  the 
isomorphism $\ft^\vee\simeq_{\langle, \rangle} \ft$ with the isomorphism $\ft^\vee= ^L\ft$, independent of any choices.
Given  two isomorphic simple Lie algebras, $\fg_1$ and $\fg_2$, there is a canonical isomorphism 
$W_1=W_2=:W$
between their  Weyl groups.
Moreover, there exists a $W$-equivariant isomorphism $\mu: \ft_1\simeq \ft_2$, unique up to the $W$-action,
such that $\mu^\vee(\cR_2)=\cR_1$.
 Such an isomorphism can be constructed, for instance, by choosing
simple coroots for both root systems, and using these to determine a linear map $\ft_1\simeq \ft_2$.
So  if $\fg$ is not of type $\sf{B}$ or $\sf{C}$, the Lie algebras $\fg$ and $^L \fg$ are (abstractly) isomorphic and we can apply to them
the above consideration, thus obtaining a preferred
  Killing form $\langle, \rangle$ :  the one for which  the composition
\[
 \ft\longrightarrow \ft^\vee=^L\ft\simeq _{\langle,\rangle}\ft
\]
 sends short coroots to long coroots. It turns out that this automorphism is in $W$ if $\fg$ is of type
$\sf{ADE}$, and not in $W$ if $\fg$ is of type $\sf{FG}$. In the latter case, however, the \emph{square} of that
automorphism \emph{is} in $W$.  In this way we obtain an involution ${\sf l}: \ft/W\to \ft/W$ and consequently, 
an involution of the Hitchin base (when $G$ is not of type $\sf{BC}$). We shall call $\sf{l}$ ``the Langlands 
involution on the (Hitchin) base''.
In \cite{argyres_kapustin_seiberg} ${\sf l}$ was interpreted as an $S$-duality transformation, acting non-trivially
on the moduli space of $N=4$ $G_2$ super Yang--Mills theories.
 
We recall the following result.
      \begin{thm}[\cite{don-pan}, Theorem A]
There is an isomorphism $\sf{l}:\cB_G\to \cB_{^LG}$, unique up to overall scalar,   which identifies the discriminant
locus $\Delta\subset \cB_G$ with $^L\Delta\subset \cB_{^LG}$ and which lifts to an automorphism $\widetilde{\sf{l}}$
of $\widetilde{\cX}$ such that
\[
 \xymatrix{\widetilde{\cX}\ar[r]^-{\widetilde{\sf{l}}}\ar[d] &^L\widetilde{\cX}\ar[d]\\
	      \cB_\fg\ar[r]^-{\sf{l}}& ^L\cB_{^L\fg}}.
\]
For each $b\in\scB=B\backslash\Delta$ there is an isomorphism of polarised abelian varieties 
$h^{-1}_0(b)\simeq \pic^0(^Lh_0^{-1}({\sf l}(b)))$, induced by a global duality of
$\Higgs^0_G$ and $\Higgs^0_{^LG}$ away from the discriminant.
      \end{thm}
An  interesting and complicated question is to understand  to what extent does the above statement extend to the case of the generalised (``meromorphic'')
Hitchin system. We refer to \cite{gukov_witten} for some insight about this case. As a preliminary check, one may ask whether the special K\"ahler metric
or the Donagi--Markman cubic on the (generalised) base $\cB_{G,D}$ is preserved by the (analogue of the) involution ${\sf l}$.
For $D=0$ and $G=G_2$ the cubic turns out to be invariant, as shown by Hitchin in \cite{hitchin_G2}.
Hitchin's approach carries over to the ramified situation more-or-less directly since the Balduzzi--Pantev formula
for the cubic  carries over to the ramified case, as discussed in (\ref{dm_cubic}).

We proceed  by recalling some basic results about the exceptional Lie algebra $\fg_2$.

The Lie algebra $\fg_2$ has rank $2$ and dimension $14$. Its Weyl group is the dihedral group $D_6$.
Since the root system is isomorphic to its dual, there is always a certain ambiguity when one wants to
describe the roots  explicitly.
We consider first  $\RR^3$ with standard basis $\{e_i\}$ and standard inner product. Then we can identify the (real) Cartan subalgebra
as $\ft_\RR = \left(e_1 + e_2 + e_3\right)^\perp\subset \RR^3$ with the induced inner product.
 The six short coroots are $\pm (e_i-e_j)$, having length $\sqrt{2}$, while the
six long coroots are $\pm (2e_i -e_j-e_k)$, having length $\sqrt{6}$.
Passing to the dual root system, we obtain  $\pm (e_i-e_j)$ for the long roots and $\pm \frac{1}{3}(2e_i -e_j-e_k)$ for the short roots.
Drawing the two root system one sees that a linear map $\mu$ mapping coroots to roots, is obtained by composing a scaling 
by a factor of $\frac{1}{\sqrt{3}}$ with rotation by $\frac{\pi}{2}$. Consequently, we can take
 ${\sf l}=\sqrt{3}\mu$. 
Explicitly, we can take ${\sf l}=\sqrt{3}\left(\begin{array}{rrr}
                                         0&-1&1\\
					  1&0&-1\\
					  -1&1&0\\
                                        \end{array}\right)$, 
where the matrix represents an automorphism of $\RR^3$, inducing the required rotation on $\ft_\RR$. Clearly, ${\sf l}\notin W$,
since $D_6$ does not contain rotations by $\frac{\pi}{2}$. At the same time, ${\sf l}^2$, rotation by $\pi$, belongs to the dihedral group.

We can get a more convenient description of $\ft$ and ${\sf l}$ if we use  the relation between $\fg_2$ and $\fs\fo_7(\CC)$
and realise the Cartan subalgebra as
\[
 \ft=
\left\{ \textrm{diag}\left(0,\lambda_1,\lambda_2,\lambda_3, -\lambda_1,-\lambda_2,-\lambda_3\right)\left| \sum_{i=1}^3\lambda_i =0,\ \lambda_i\in\CC\right.\right\} \subset \mend(\CC^7).
\]
Thus any $a\in \ft$  is determined by a matrix $A=\textrm{diag}(\lambda_1,\lambda_2,\lambda_3)\in \slt$. Its characteristic polynomial
is $\lambda\left(\lambda^6 - f \lambda^4 + \frac{f^2}{4}\lambda^2 -q\right)$, where
$f=\textrm{tr}A^2 = \sum_{i=1}^3\lambda_i^2$, $q=\det A^2=\lambda_1^2\lambda_2^2\lambda_3^3$.
We see then that the degrees of the generators of $\CC[\fg_2]^{G_2}$ are $d_1=2$, $d_2=6$.
If we use the above invariants as generators, and use the simple coroots $(e_1-e_2, 2e_2 - e_1-e_3)$
as a basis of $\ft$, we get an explicit expression for the adjoint quotient morphism:
\[
 \chi: \CC^2\to \CC^2,\ \chi(x,y)=(I_1(x,y), I_2(x,y))
\]
where 
      \begin{equation}\label{g2_invariants}
	\left| 
       \begin{array}{l}
	    I_1(x,y)= 2x^2 -  6xy +2y^2\\
	    I_2(x,y)= 4y^6 - 12 y^5 x + 13 y^4x^2 - 6 y^3x^3 + x^4y^2.
       \end{array}
	\right.
      \end{equation}
Since we know explicitly the roots of $\fg$ and have an expression for the invariants, we can compute directly that
      \begin{equation}\label{langlands_base}
       {\sf l}(f,q)= \left(f, -q + \frac{f^3}{54}\right)
      \end{equation}
and
      \begin{equation}
       \fD= \prod_{\alpha\in \cR}\alpha= 27 q\left(-q + \frac{f^3}{54}\right) = 27 q \check{q},
      \end{equation}
where we write $\check{q}:= -q + \frac{f^3}{54}$ following Hitchin.
Notice that the discriminant  $\fD= {\sf l}(\fD)$ is reducible and the two factors are the products over short and long roots, respectively.
This is a common feature for non simply-laced Lie algebras, since the Weyl group acts transitively on the set of roots of equal length.

We now recall briefly the numerical invariants for the generalised $G_2$ Hitchin system.
We fix a divisor $D$ on $X$ of $\deg D=\delta$, and assume, as in Section (\ref{dm_cubic}), that the square
of $L=K_X(D)$ is very ample.
Via the choice of $(I_1,I_2)$ the Hitchin base can be identified as
$\cB_{G_2,D}= H^0(L^2)\oplus H^0(L^6)$ and $\dim \cB_{G_2,D}= 14(g-1) + 8\delta$.
The group $G_2$ is simultaneously simply-connected and adjoint and
$\dim \Higgs_{G_2,D}^0=\dim G_2\deg L= 28(g-1) + 14\delta$.
In Section (\ref{dm_cubic}) we also introduced the subspace
$\cB_0= H^0(L^2(-D)\oplus L^6(-D))\subset\cB_{G_2,D}$ and one checks that 
$\dim \cB_0= 14(g-1)+ 6\delta$. 
In particular, $\cB_{G_2,D}/\cB_0$, the space of closures of symplectic leaves, has dimension $2\delta$.

    \begin{thm}
Let ${\sf l}:\cB_{G_2,D}\to \cB_{G_2,D}$ be the automorphism  (\ref{langlands_base}) and $c$ the Donagi--Markman cubic
for the generalised $G_2$ Hitchin system. Then ${\sf l}^\ast c=c$.     
    \end{thm}
\emph{Idea of proof:}
For a fixed point $o=(f_0, q_0)\in \cB_{G,D}$ the equation of the cameral curve $\widetilde{X}_o\subset \ft\ctimes L$ is $I_1=f_0$, 
$I_2=q_0$, where $(I_1,I_2)$ are given in (\ref{g2_invariants}).
By applying a linear automorphism of $\ft\simeq \CC^2$ we can rewrite these equations  as
\[
 \left|
      \begin{array}{l}
       \frac{3}{2}\left(x^2+ y^2\right)=f_0\\
	\frac{1}{16}\left(x^6 - 6 x^4y^2 + 9 x^2 y^4\right)=q_0\\
      \end{array}
\right.,
\]
which is equivalent to equation (26) in \cite{hitchin_G2}, namely, 
\[
 \left|
      \begin{array}{l}
       x^2+ y^2 =\frac{2}{3}f_0\\
	x^6 - f_0 x^4  + \frac{f_0^2}{4} x^2 =q_0\\
      \end{array}
\right. .
\]
The advantage of the latter description is that $\widetilde{X}_0$ is realised as a double cover of an intermediate curve
in a way which makes the action of ${\sf l}$ on the space of cameral curves explicit. The calculation then is completed as in
\S 6.5 of \cite{hitchin_G2}, but using the expression from Theorem \ref{thm_cubic} for the cubic.

\qed

\bibliographystyle{alpha}
\bibliography{biblio}

\end{document}